\RequirePackage{fix-cm}
\documentclass[]{svjour3}                     
\usepackage[letterpaper,margin=01.in]{geometry}
\smartqed  
\usepackage[utf8]{inputenc}
\usepackage{amssymb, amsmath}
\usepackage{graphicx}
\usepackage{graphicx}
\usepackage{subcaption}
\usepackage[shortlabels]{enumitem}
\usepackage{natbib}
 \journalname{Journal of Scientific Computing}
\begin{document}

\title{An efficient Haar wavelet method for the coupled non-linear transient PDE-ODEs system with discontinuous coefficients
}

\titlerunning{An efficient HWM for the coupled non-linear transient PDE-ODEs system with discontinuous coefficients}        

\author{B.V. Rathish Kumar* \and
        Meena Pargaei 
}


\institute{B.V. Rathish Kumar* \at
               Department of Mathematics and Statistics, Indian Institute of Technology, Kanpur, India \\
              \email{drbvrk11@gmail.com*}           
           \and
          Meena Pargaei \at
              Department of Mathematics and Statistics, Indian Institute of Technology, Kanpur, India \\
              Department of Mathematics, G.P.G.C. Champawat, Uttarakhand, India \\
              \email{meenamenu15@gmail.com}
}


\maketitle

\begin{abstract}
In this work, Haar wavelet method for the coupled non-linear transient PDE-ODEs system with Neumann boundary condition has been proposed. The capability of the method in handling multiple jump discontinuities in the coefficients and parameters of the transient and coupled PDE-ODEs system is brought out through a series of 1D/2D/3D test problems. The method is easy to implement, computationally efficient and compare well with conventional methods like finite element method. Convergence analysis of the method has been carried out and apriori error is found to be exponential order i.e. $\|u-u_H\|_X \sim o(2^{-j})$. ILU-GMRES is found to better accelerate the numerical solution convergence than other Krylov solvers for the class of transient non-linear PDE-ODE system.
\end{abstract}

\section{Introduction}
\label{intro}
Study of coupled PDE-ODEs systems has been an active research theme as such models arise in different fields of physical science, engineering, biology and economics. For instance such models are encountered while dealing with problems related to electro-magnetic coupling, coupling of chemical reaction and transport, mechanical coupling, propagation of electric waves in cardiac tissue etc. J. Doafouz et al. \cite{JD} model for a boost converter connected to a load via a transmission line, coupled the partial differential equation describing telegraph with the ODE representing the action of the converter. In \cite{SE}, E. Sabine presented the mathematical model of mitochondrial swelling phenomena via coupled PDE-ODE system. Such PDE-ODE coupled systems frequently arise in the field of cardiac electrophysiology \cite{Plonsey,keener,SRWJ}. Such PDE-ODE coupled systems are not easily solvable especially when they are non-linear and in addition they have discontinuities or sharp variations either in field variables or coefficients or parameters appearing in the governing mathematical model. This is especially the case in mathematical models representing the cardiac electrical activity in ischemic cardiac tissue. Such discontinuities are also encountered in solid mechanics while dealing with fractures, cracks, inclusions and voids. In fluid mechanics such sharp variations in field variables arise when dealing with shock waves, boundary layers etc.

Several specialized numerical methods based on different operator discretization procedures such as finite differences, finite elements, spectral elements, finite volume etc. have been used to solve such coupled systems with discontinuities or sharp variations. M. Vynnycky and Sarah L. Mitchell \cite{VS} proposed Keller box based finite difference method for parabolic PDE with discontinuous boundary conditions. 
Several finite element method based adaptive methods have been proposed to deal with jumps, kinks or sharp variations in PDE model. Andreas Bomto et al. in \cite{BRR} proposed adaptive finite element method for elliptic problems with discontinuous coefficients for dealing with diffusion through porous media, or rough surfaces, electromagnetic field propagation on heterogeneous media using distortion theory. C. Bernardi and R. Verfurth in \cite{BR} derived apriori and aposteriori errors for adaptive finite element method for elliptic equations with non-smooth coefficients for modeling multilayer fluids with different viscosities. In \cite{BD}, Andrea Bonito and Denis Devand introduced the adaptive finite element method for stokes problem with discontinuities viscosity. There are also meshfree approximation techniques which provide greater flexibility. C. Armando Duarte and J.T. Oden \cite{DJ} implemented h-p cloud method, fully exploiting the h, p and h-p adaptivity, devoid of conventional finite element grids on Poisson problems and 3-D elasticity problem. However these methods are not only limited by their adaptivity approach but may turn out to be very costly especially when dealing with the problems with evolving and moving features. In \cite{GFEM}, generalized or extended finite element method has been discussed especially to effectively handle the jumps and singularities but this again gets quite sensitive to indicator function choice and are computationally costly.

Recently, wavelet methods are getting much attention, due to their inherent adaptability to the complexities such as discontinuities, sharp variations etc. in solving the mathematical models from science and engineering. Its properties such as orthogonality, compact support, arbitrary regularity and high order vanishing moments are very attractive. Wavelet based Galerkin and collocation methods have been used to solve integral equations, ODEs, PDEs, fractional differential equations\cite{aziz,siraj,siraj2,aziz2,siraj3,siraj4,aziz3,siraj5,naldi,wu}. In \cite{WAMR} parallel wavelet adaptive method has been used in the domain decomposition framework to solve the compressible reacting flow with the shock wave in air. Oleg V. Vasilyev et al. \cite{OTD} have developed discontinuous interpolating wavelets for handling multiphase problems with discontinuities across the phase boundaries. Kai Schneider and Oleg V. Vasiley \cite{CFD} demonstrated the ability of adaptive wavelet basis in capturing coherent vortex structures in turbulent flows and thereby zoomed into their space-scale structure. They also present a hierarchy of wavelet-based turbulence models for carrying out complex industrial flow simulations. In \cite{KKK} Sebastion Kestler et al. have introduced a space-time adaptive wavelet Galerkin method for time periodic PDEs. In all these studies authors have focused on using either Daubechies class of continuous wavelets or second generation wavelets. Simplest of all the wavelets and easy to implement in a finite domain is the Haar wavelet.

Haar wavelets are piecewise constant functions which are orthogonal and have compact support. They also have scaling property. The ease of using a Haar wavelet has made it very popular in signal processing and image processing. Because of the discontinuity of Haar wavelets, their derivatives does not exists and hence cannot be directly used in approximating solution of a differential equation. Chen and Hsiao \cite{chen} have proposed the idea that the highest order derivative of the differential equation can be expanded into the Haar series and not the field variable function. Then on integration one can obtain lower order derivatives and the function too. 

Haar wavelet method has been used to solve the linear and non-linear differential equations and the eigenvalue problems \cite{bujurke,lepik,kannan,sheo,ev}. Largely these attempts are confined to 1D studies and this is specially the case with transient PDEs. Further no attempt has been made to test Haar wavelet capability in handling complexities like discontinuities, sharp variations etc. In this work, we propose the Haar wavelet method for the solution of coupled PDE-ODEs system with Neumann boundary condition and having discontinuities in either in source/coefficient function or in any parameter of the model. We will show the power of method in handling problem with multiple discontinuities like those which occur in cardiac electrophysiological models. We easily extend the method to higher dimensions and also to identify the desired linear solver and the influence of preconditioning in convergence acceleration. We also theoretically carry out the error analysis and show that the convergence is of exponential order i.e. $o(2^{-j})$. We compare the results against those from finite element method both for quality and computational complexity. The paper is organized as follows:

In section 2, we introduce Haar wavelet function, its properties and haar integrating functions. In the next section, one and two dimensional Haar wavelet method for non-linear PDE-ODEs system will be developed. Convergence analysis of the proposed method is presented in section 4. Numerical results and discussion for the several problems with different types of discontinuities are provided in section 5. Three dimensional Haar wavelet method is discussed in the next section and its convergence analysis has been established. Further the proposed method has been successfully demonstrated for few practical problems, with discontinuity from cardiac electrophysiology in ischemic cardiac tissue.

\section{Haar wavelets}
Let us consider the interval $x \in [A,B]$, A, B are finite real numbers. Define $M=2^J$, $J$ is 
the maximum level of resolution. This interval $[A,B]$ is equally divided into $2M$ subintervals such that the length of each subinterval is $\Delta x = (B-A)/2M$. Now, define the dilation and 
translation parameter $j=0,1,...,J$ and $k=0,1,...,m-1$ respectively, where, $m=2^j$. The wavelet number is given by $i=m+k+1$. 
Family of haar wavelets is defined as follows:

For $i=1$
\begin{align*}
h_i(x) = 
\begin{cases}
1 & \text{when } A \leq \ x < B\\    
0 & \text{otherwise}.
\end{cases}
\end{align*}

For $i \geq 2$
\begin{align}
\label{haar}
h_i(x) = 
\begin{cases}
1 & \text{when } \beta_1(i) \leq \ x < \beta_2(i)\\    
- 1 & \text{when } \beta_2(i) \leq \ x < \beta_3(i)\\  
0 & \text{otherwise},
\end{cases}
\end{align}
where, 
\begin{align*}
& \beta_1(i)= A+ 2k\zeta \Delta x, \hspace{0.5cm} \beta_2(i)= A+ (2k+1)\zeta \Delta x, \\
& \beta_3(i)= A+ 2(k+1)\zeta \Delta x,  \hspace{0.5cm} \zeta=M/m.
\end{align*}

Haar wavelets are orthogonal, since
\begin{align*}
\int_A^B h_p(x) h_q(x) dx = 
\begin{cases}
2^{-q}(B-A) & \text{when } p=q\\    
0 & \text{when } p \neq q.
\end{cases}
\end{align*}

For the solution of differential equation, we have to compute the integral
\begin{align}
p_{\alpha,i}(x)= \idotsint_A^B h_i(t)dt^{\alpha} = \frac{1}{(\alpha -1)!}
\int_A^x (x-t)^{(\alpha -1)} h_i(t) dt,
\end{align}

where $\alpha = 1, 2,..., n$ and $i=1,2,..., 2M$. 

For the case, $\alpha =0, p_{0,i}(x)=h_i(x)$.

This integral is calculated with the help of equation \eqref{haar}, defined as
\begin{align}
\label{HaarIntegral}
p_{\alpha ,i}(x) = 
\begin{cases}
0 & \text{when } x < \beta_1(i)\\    
\frac{1}{\alpha !}(x-\beta_1(i))^{\alpha} & \text{when } \beta_1(i) \leq \ x < \beta_2(i)\\ 
\frac{1}{\alpha !}[(x-\beta_1(i))^{\alpha}-2(x-\beta_2(i))^{\alpha}] & \text{when } \beta_2(i) \leq \ x < \beta_2(i)\\  
\frac{1}{\alpha !}[(x-\beta_1(i))^{\alpha}-2(x-\beta_2(i))^{\alpha}+(x-\beta_3(i))^{\alpha})] & \text{when }  x > \beta_3(i). 
\end{cases}
\end{align}

When $\alpha =1, 2$,
\begin{align}
p_{1,i}(B) = 
\begin{cases}
(B-A) & \text{when } i=1\\    
0 & \text{otherwise}.
\end{cases}
\end{align}

\begin{align*}
p_{2,i}(B) = 
\begin{cases}
(B-A)^2/2 & \text{when } i=1\\    
(B-A)^2/4m^2 & \text{otherwise}.
\end{cases}
\end{align*}

For the grid points $x_k=A+kh, k=0,1,...,2M$, $h=\Delta x$, collocation points are as follows:
\begin{align}
y_k= \frac{x_{k-1}+x_{k}}{2}, k=1,2,...,2M.
\end{align}

After this discretization, we define Haar matrix $H$, and Haar Integral matrices $P_1 , P_2$ of
size $2M \times 2M$ as $H(i,k)=h_i(y_k), P_1(i,k)=p_{1,i}(y_k) , P_2(i,k)=p_{2,i}(y_k)$.

\subsection{Function approximation}
Any function $f(x) \in L^2[0,1)$ can be approximated in terms of the Haar wavelet series as 
\begin{align*}
f(x)=\sum_{i=1}^{\infty}\alpha_i h_i(x),
\end{align*}
where the wavelet coefficients $\alpha_i$ are obtained by
\begin{align*}
\alpha_i= 2^j \int_0^1 f(x) h_i(x)dx.
\end{align*}

Since only the finite number of terms are taken for the computational purpose therefore, the function approximation $f$ is given by
\begin{align*}
f(x)=\sum_{i=1}^{2M}\alpha_i h_i(x).
\end{align*}
 where, $M=2^J$
\subsection{Functions consisting the coefficient having jump discontinuities}
Let us consider the domain $[0,1]$ and a function $f(x)$ defined on this domain given by
\begin{align}
f(x)=k(-x^3+1.1x^2-0.1x)
\end{align}
where, $k$ is a parameter.

\begin{figure}
	\vspace{-12em}
	\centering
	\begin{subfigure}[t]{0.5\textwidth}
		\centering
		\includegraphics[width=1.\textwidth]{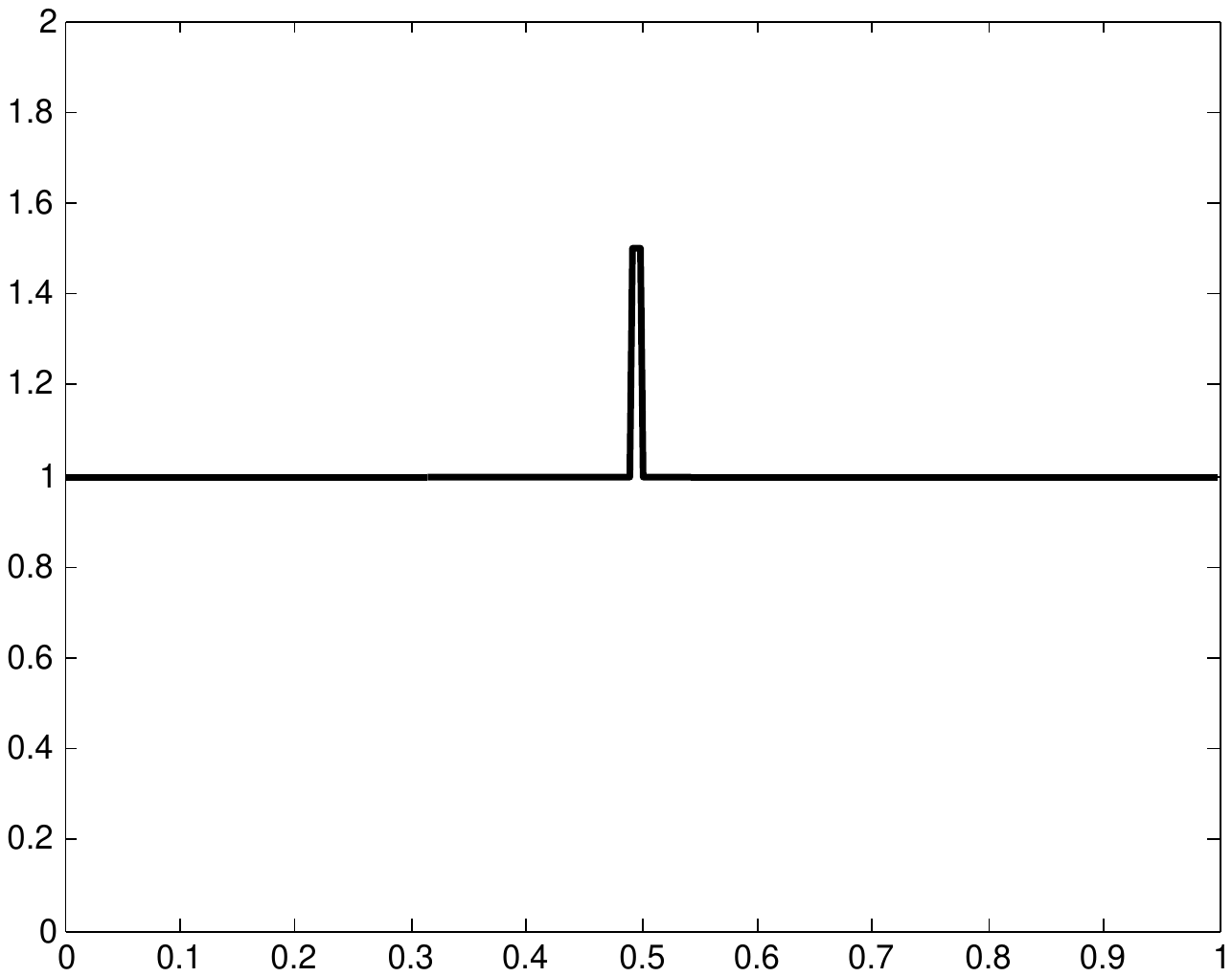}
		\vspace{-12em}
		\caption{}
		\label{1jump}
	\end{subfigure}\hfill
	\begin{subfigure}[t]{0.5\textwidth}
		\centering
		\includegraphics[width=1.\textwidth]{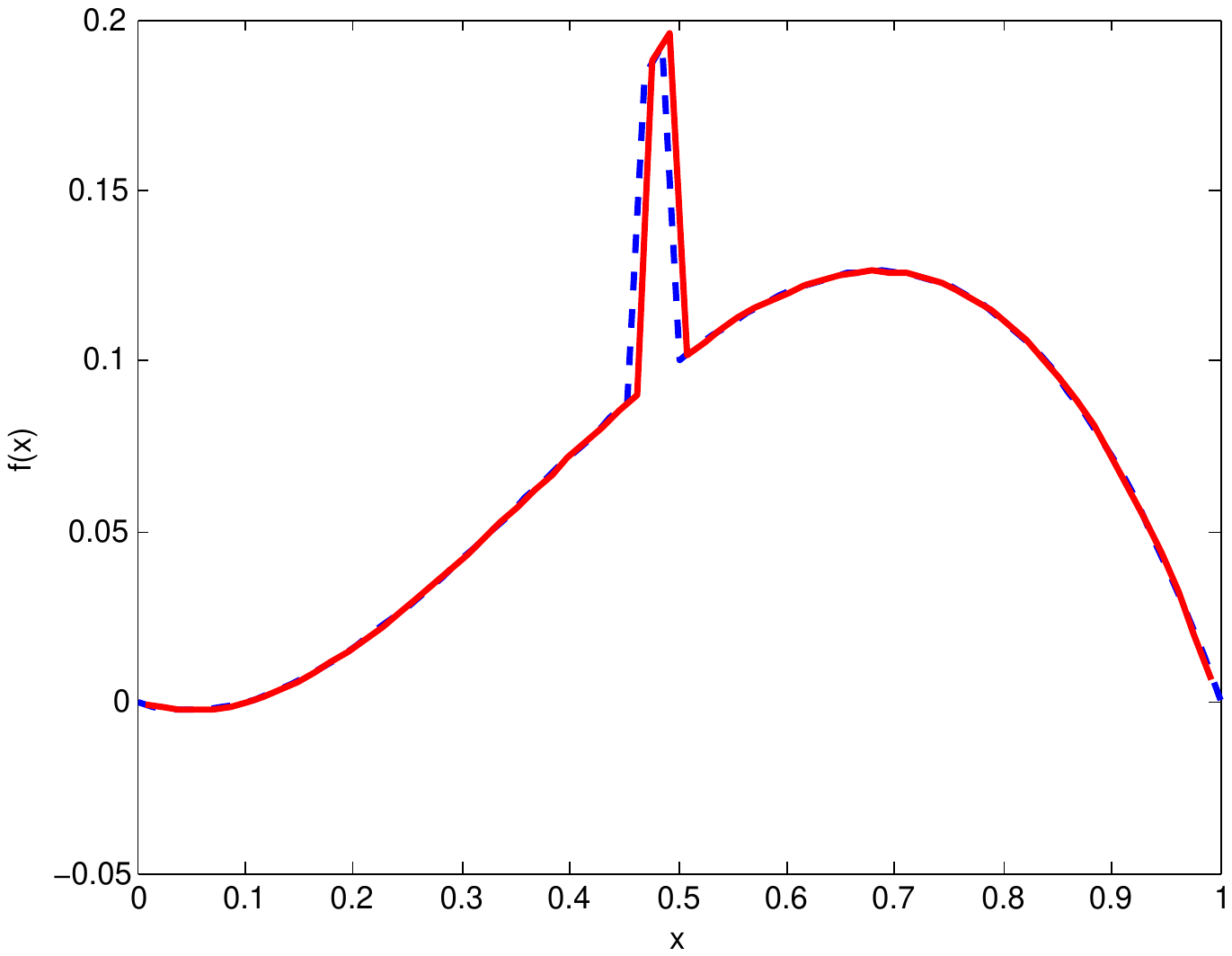}
		\vspace{-12em}
		\caption{}
		\label{f_1jump}
	\end{subfigure}\hfill 
	\vspace{-20em}    	
	\begin{subfigure}[t]{0.5\textwidth}
		\centering
		\includegraphics[width=1.\textwidth]{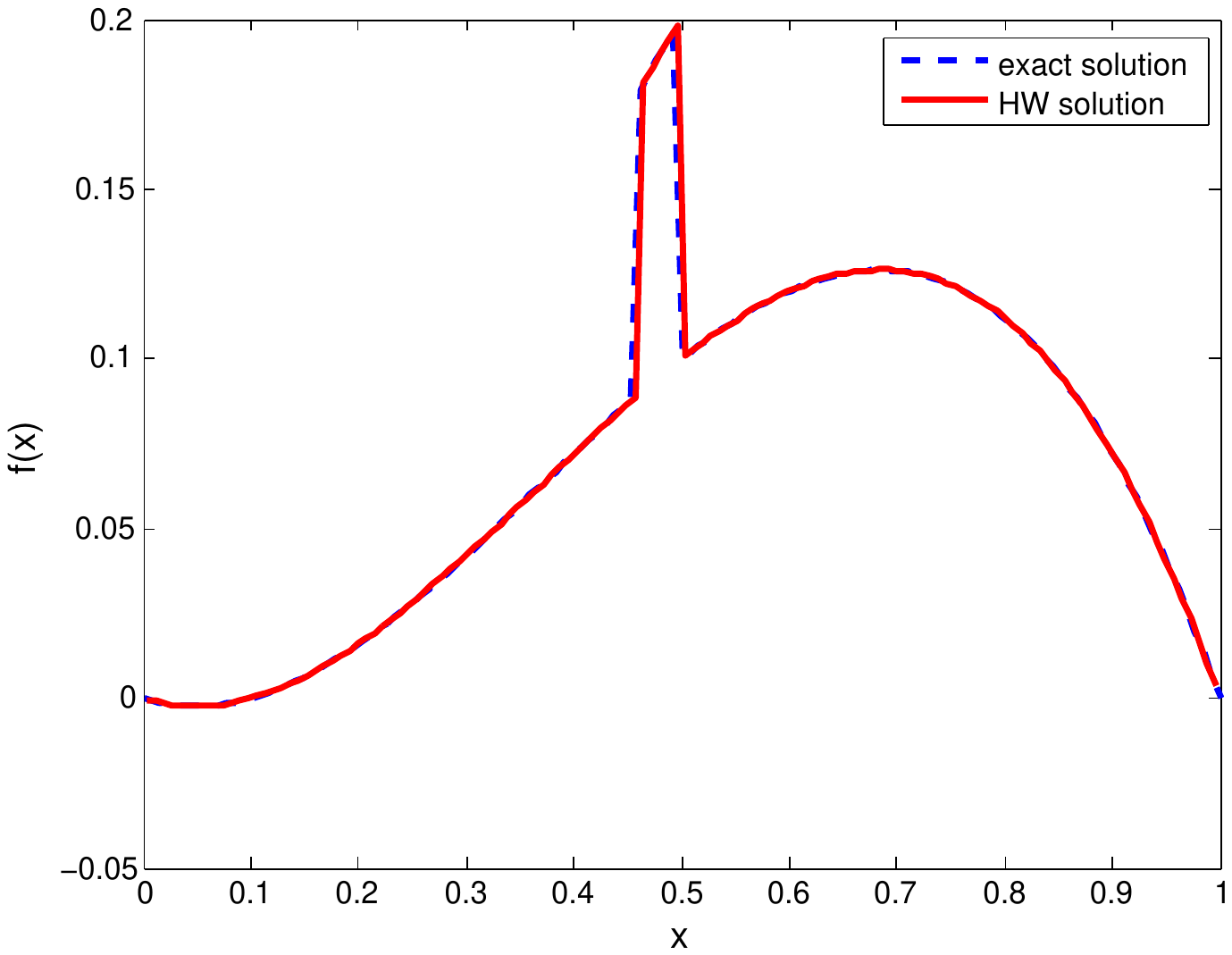}
		\vspace{-12em}
		\caption{}
		\label{f_1jump_J6}
	\end{subfigure}\hfill
	\begin{subfigure}[t]{0.5\textwidth}
		\centering
		\includegraphics[width=1.\textwidth]{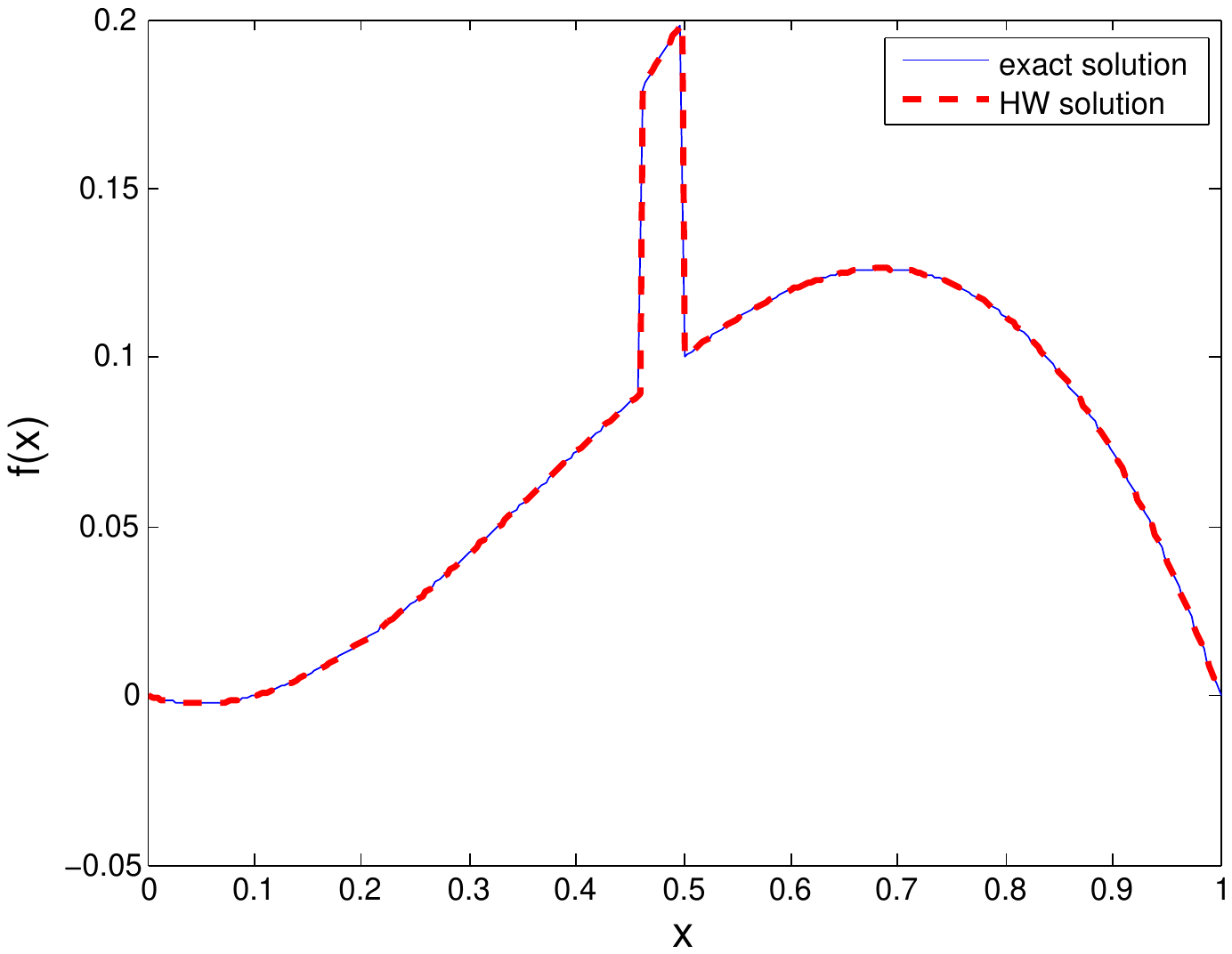}
		\vspace{-12em}
		\caption{}
		\label{f_1jump_J7}
	\end{subfigure}
	\vspace{-10em}
	\caption{\textbf{(a) Value parameter $k$, (b) Corresponding function $f(x)$, $J=5$, (c) $J=6$, (d) $J=7$.}}
	\label{FHN_jump}
\end{figure}

\begin{figure}
	\vspace{-12em}
	\centering
	\begin{subfigure}[t]{0.5\textwidth}
		\centering
		\includegraphics[width=1.\textwidth]{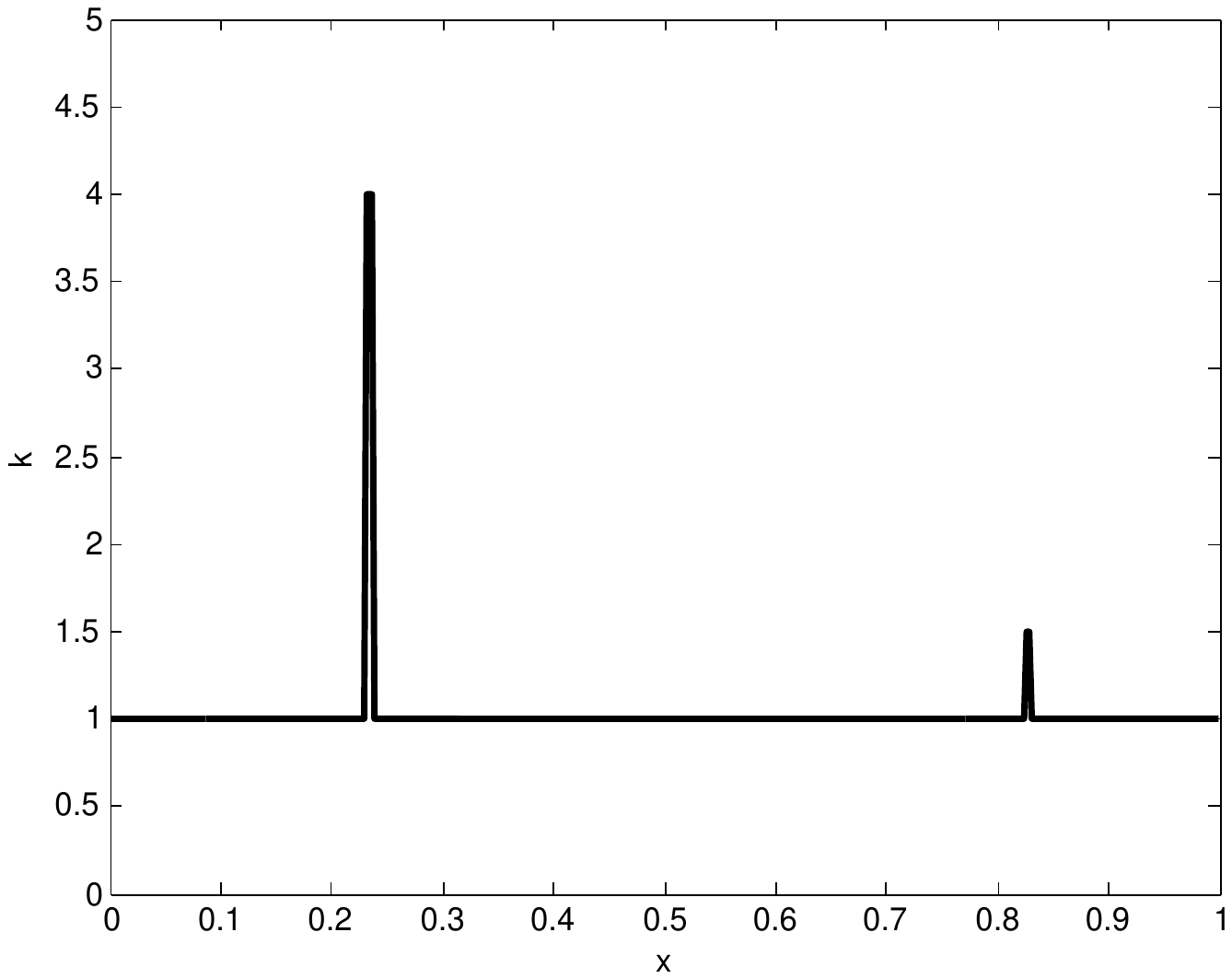}
		\vspace{-12em}
		\caption{}
		\label{2jump}
	\end{subfigure}\hfill
	\begin{subfigure}[t]{0.5\textwidth}
		\centering
		\includegraphics[width=1.\textwidth]{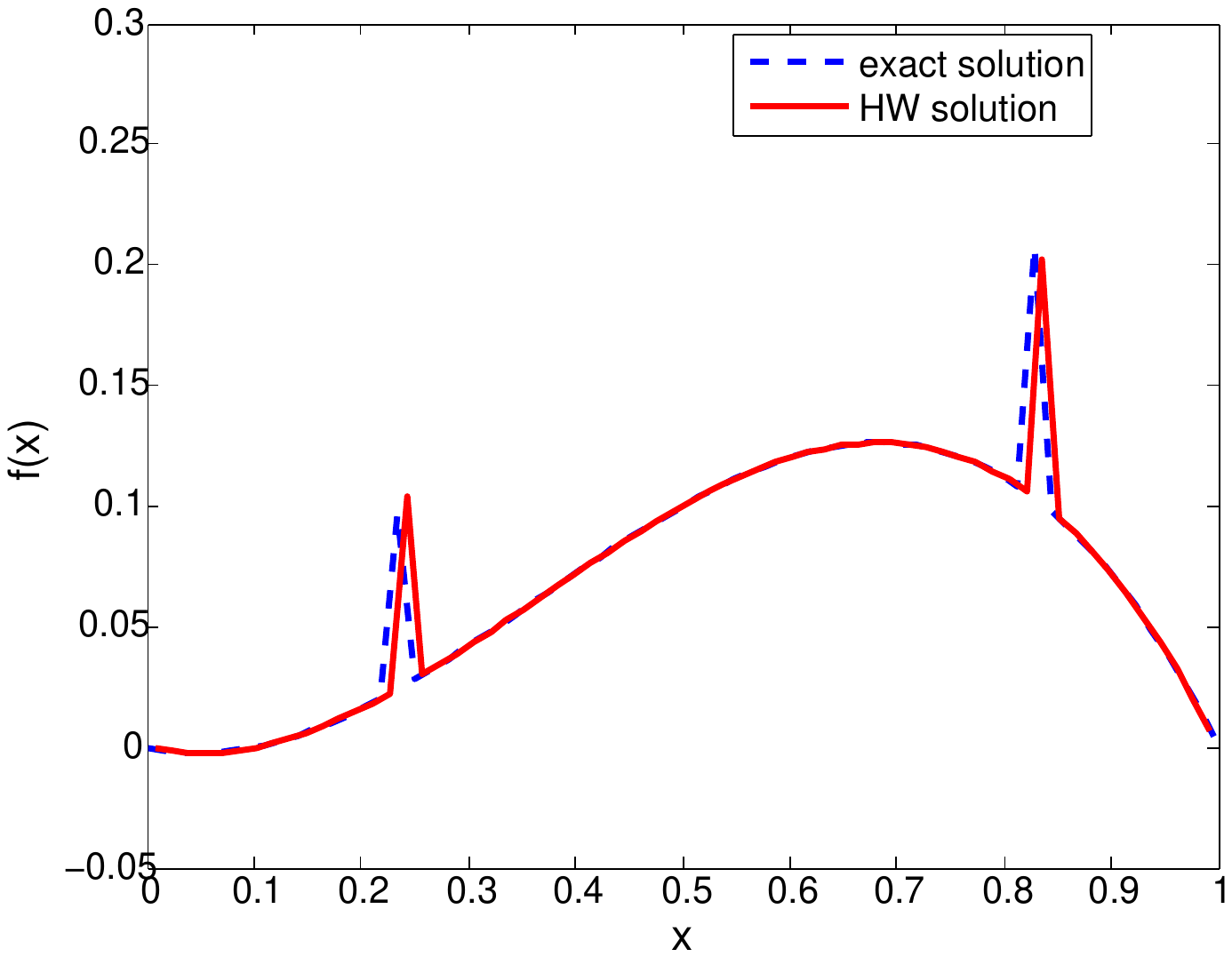}
		\vspace{-12em}
		\caption{}
		\label{f_2jump}
	\end{subfigure}\hfill 
	\vspace{-19em}  
	\begin{subfigure}[t]{0.5\textwidth}
		\centering
		\includegraphics[width=1.\textwidth]{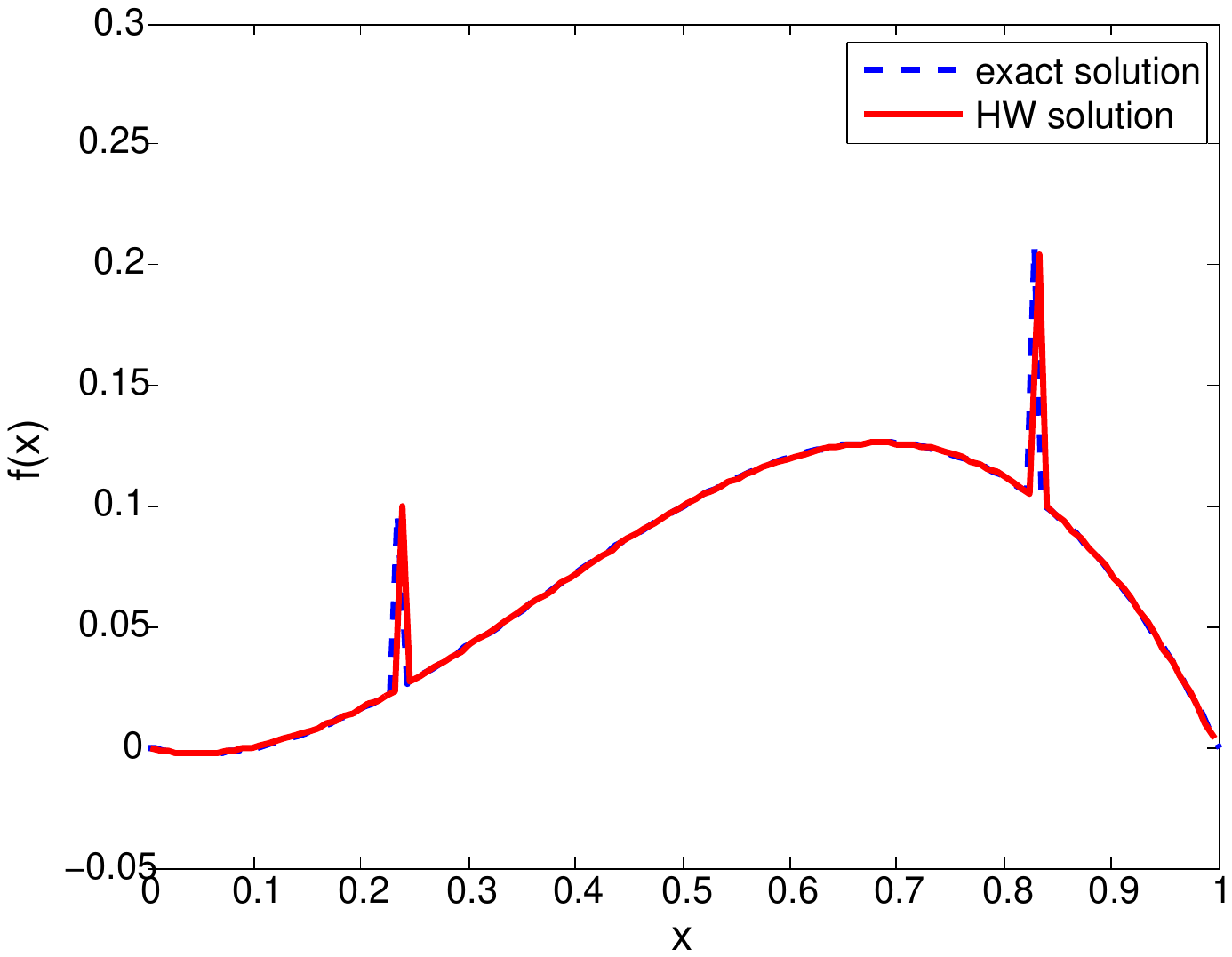}
		\vspace{-12em}
		\caption{}
		\label{f_2jump_J6}
	\end{subfigure}\hfill
	\begin{subfigure}[t]{0.5\textwidth}
		\centering
		\includegraphics[width=1.\textwidth]{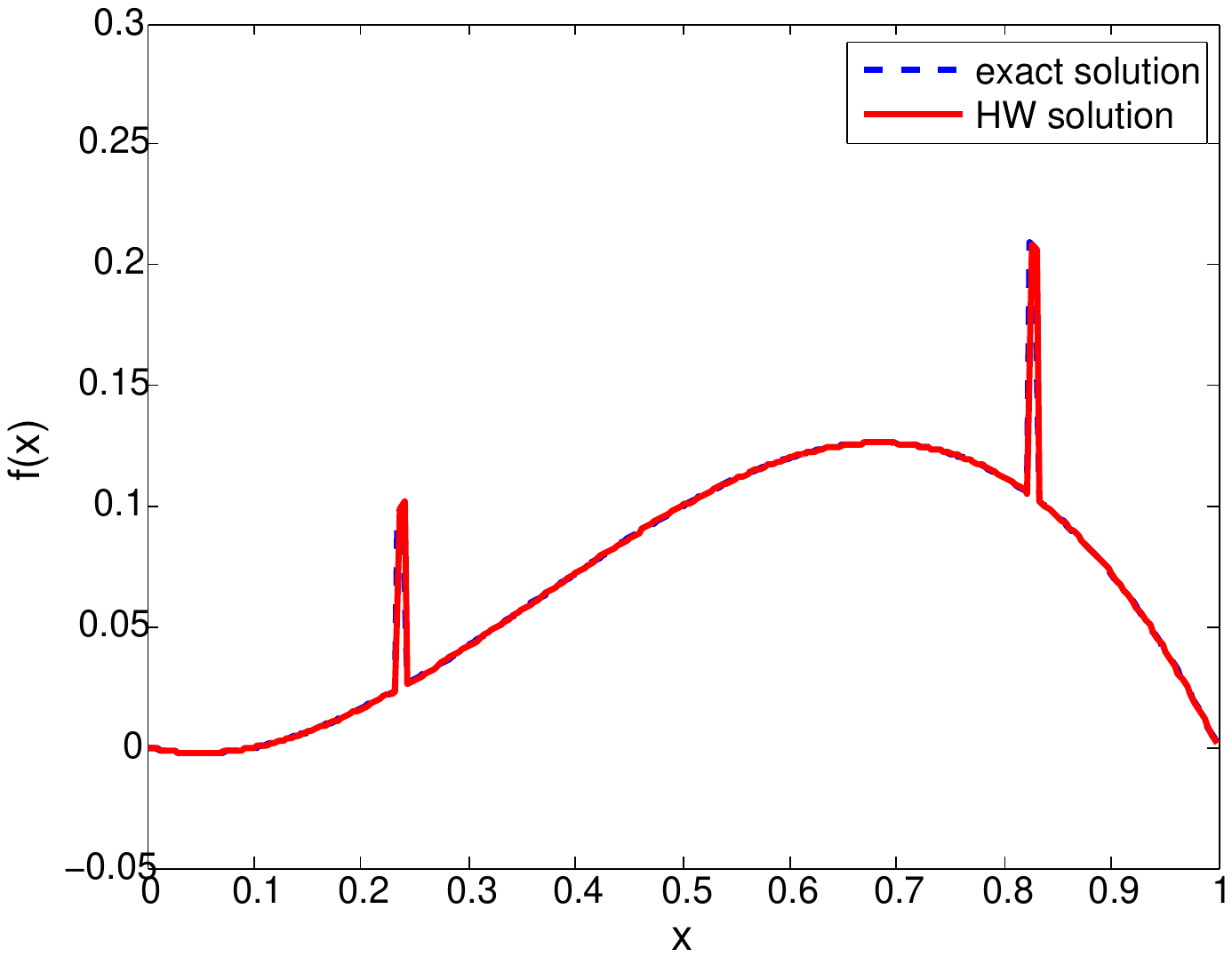}
		\vspace{-12em}
		\caption{}
		\label{f_2jump_J7}
	\end{subfigure}
	\vspace{-10em}
	\caption{\textbf{(a) Value parameter $k$, (b) Corresponding function $f(x)$ $J=5$ , (c) $J=6$, (d) $J=7$.}}
	\label{FHN_jump}
\end{figure}

Now, if the parameter $k$ of the above function $f(x)$ is discontinuities at one or multiple places of the domain, as given in
Fig. \ref{1jump} $(k=1.5 if 0.46<x<0.5, k=1$ otherwise) and Fig. \ref{2jump} $(k=4 if 0.2305<x<0.2383, k=2 if 0.8242<x<0.8302, k=1$ otherwise),  the corresponding function $f(x)$ will also be discontinuous, drawn in Fig. \ref{f_1jump}
, \ref{f_1jump_J6}, \ref{f_1jump_J7} \ref{f_2jump}, \ref{f_2jump_J6}, and, \ref{f_2jump_J7} respectively.

\section{Mathematical Model}
Let $\Omega \subset R^n$ then a general non-linear parabolic reaction-diffusion equation coupled with the system of ODEs, is given as follows: 
\begin{align*}
\epsilon\frac{\partial v}{\partial t}- div(D(x)\nabla v) + f(v,w)  &= 0 &  x \in \Omega,  0\leq t \leq T\\
\frac{\partial w_i}{\partial t}-g_i(v,w)&= 0, & x \in \Omega,  0\leq t \leq T\\
v(x,0)= v_0(x,0), \hspace{5mm} w_i(x,0)&=w_{i,0}(x,0), &x \in \Omega\\
n^T D(x) \nabla v &=0, & x \in \partial \Omega, 0\leq t \leq T,
\end{align*}
where, $i=1,2,...d$ and  $D(x)$ is the diffusion tensor.

Remark: For existence-uniqueness and stability of such coupled system one may refer to \cite{bourgault, veneroni}.



\subsection{Haar wavelet algorithm for the coupled non-linear PDE-ODEs system}
Consider the three dimensional non-linear parabolic reaction-diffusion equation coupled with a system of ODEs given as follows:
\begin{align}
\label{mv3d}
&\epsilon(x,y,z)\frac{\partial v}{\partial t}- div(D(x,y,z)\nabla v) + f(v,w)  = 0, & 0\leq x,y,z \leq 1,  0\leq t \leq T\\
\label{mw3d}
&\frac{\partial w}{\partial t}-g(v,w)= 0, & 0\leq x,y,z \leq 1,  0\leq t \leq T\\
&v(x,y,z,0)= v_0(x,y,z,0), \hspace{5mm} w(x,y,z,0)=w_0(x,y,z,0), & 0\leq x,y,z \leq 1
\end{align}
with Neumann boundary conditions in $v$. 

Let us write $\frac{\partial ^7 v}{\partial t \partial x^2 \partial y^2 \partial z^2}
\big (x,y,z,t\big )$ in terms of the Haar wavelet as follows:

\begin{align}
\label{approxv3d}
\frac{\partial ^7 v}{\partial t \partial x^2 \partial y^2 \partial z^2} (x,y,z,t)&
=\sum_{i,j,k=1}^{2M}\alpha_{i,j} h_i(x) h_j(y) h_k(z), \hspace{1cm}t\in [t_s, t_{s+1})\\
\label{approxw3d}
\frac{\partial w}{\partial t}(x,y,z,t) &= \sum_{l,m,n=1}^{2*M}\beta_{l,m,n} h_l(x) h_m(y)
h_n(z) , \hspace{0.5cm}t\in [t_s, t_{s+1}).
\end{align}
Integrating equation \eqref{approxv3d} w.r.t  $t$ from $t_s$ to $t$, we will get 
\begin{align}
\label{3dv6}
\frac{\partial ^6 v}{ \partial x^2 \partial y^2  \partial z^2}(x,y,z,t) &= (t-t_s)\sum_{i,j,k=1}^{2M}\alpha_{i,j,k} h_i(x) h_j(y)h_k(z)
+\frac{\partial ^6 v}{ \partial x^2 \partial y^2  \partial z^2} (x,y,z,t_s), \hspace{0.5cm}t\in [t_s, t_{s+1}).
\end{align}
Now, Integrate equation \eqref{3dv6} twice w.r.t  $x$ from $0$ to $x$ also using the
boundary conditions, we will obtain the following
\begin{align}
\label{3dv4}
\nonumber
\frac{\partial ^4 v}{ \partial y^2 \partial z^2}(x,y,z,t) &= (t-t_s)\sum_{i,j,k=1}^{2M}
\alpha_{i,j,k} P_{2,i}(x) h_j(y)h_k(z)
+\frac{\partial ^4 v}{ \partial y^2  \partial z^2} (x,y,z,t_s) -\frac{\partial ^4 v}
{\partial y^2  \partial z^2} (0,y,z,t_s) \\ &+ \frac{\partial ^4 v}{ \partial y^2  
	\partial z^2} (0,y,z,t), \hspace{0.5cm}t\in [t_s, t_{s+1}).
\end{align}
Now, Integrate equation \eqref{3dv4} twice w.r.t  $y$ from $0$ to $y$ also using the 
boundary conditions, we get
\begin{align}
\label{3dvz2}
\nonumber
\frac{\partial ^2 v}{\partial z^2}(x,y,z,t) &= (t-t_s)\sum_{i,j,k=1}^{2M}\alpha_{i,j,k} 
P_{2,i}(x) P_{2,j}(y)h_k(z) + 
\frac{\partial ^2 v}{\partial z^2}(x,y,z,t_s)-\frac{\partial ^2 v}{\partial z^2}(x,0,z,t_s)
-\frac{\partial ^2 v}{\partial z^2}
(0,y,z,t_s)\\
&+ \frac{\partial ^2 v}{\partial z^2}(0,0,z,t_s)+ \frac{\partial ^2 v}{\partial z^2}
(0,y,z,t)- \frac{\partial ^2 v}{\partial z^2}(0,0,z,t)+ \frac{\partial ^2 v}{\partial z^2}(x,0,z,t),
\hspace{0.5cm}t\in [t_s, t_{s+1}).
\end{align}
Similarly, Integrate equation \eqref{3dv4} twice w.r.t  $z$ from $0$ to $z$ also using the boundary conditions, we get	
\begin{align}
\label{3dvy2}
\nonumber
\frac{\partial ^2 v}{\partial y^2}(x,y,z,t) &= (t-t_s)\sum_{i,j,k=1}^{2M}\alpha_{i,j,k} 
P_{2,i}(x)h_j(y) P_{2,k}(z) + 
\frac{\partial ^2 v}{\partial y^2}(x,y,z,t_s)-\frac{\partial ^2 v}{\partial y^2}(x,y,0,t_s)-
\frac{\partial ^2 v}{\partial y^2}
(0,y,z,t_s)\\
&+ \frac{\partial ^2 v}{\partial y^2}(0,y,0,t_s)+ \frac{\partial ^2 v}{\partial y^2}(0,y,z,t)
- \frac{\partial ^2 v}{\partial y^2}
(0,y,0,t)+ \frac{\partial ^2 v}{\partial y^2}(x,y,0,t),
\hspace{0.5cm}t\in [t_s, t_{s+1}).
\end{align}
Again, Integrating \eqref{3dv6} twice w.r.t  $z$ from $0$ to $z$ and then twice w.r.t 
$y$ from $0$ to $y$ also using the boundary conditions, we get
\begin{align}
\label{3dvx2}
\nonumber
\frac{\partial ^2 v}{\partial x^2}(x,y,z,t) &= (t-t_s)\sum_{i,j,k=1}^{2M}\alpha_{i,j,k} h_i(x) P_{2,j}(y) P_{2,k}(z) + 
\frac{\partial ^2 v}{\partial x^2}(x,y,z,t_s)-\frac{\partial ^2 v}{\partial x^2}(x,0,z,t_s)-\frac{\partial ^2 v}{\partial x^2}(x,y,0,t_s)\\
&+ \frac{\partial ^2 v}{\partial x^2}(x,0,0,t_s)+ \frac{\partial ^2 v}{\partial x^2}(x,0,z,t)- 
\frac{\partial ^2 v}{\partial x^2}(x,0,0,t)+ \frac{\partial ^2 v}{\partial x^2}(x,y,0,t),
\hspace{0.5cm}t\in [t_s, t_{s+1}).
\end{align}
Now, Integrate equation \eqref{approxv3d} twice w.r.t. x , y and then z also using boundary conditions, we will obtain
\begin{align}
\label{3dvt}
\nonumber
\frac{\partial v}{\partial t}(x,y,z,t) &= \sum_{i,j,k=1}^{2M}\alpha_{i,j,k} P_{2,i}(x) P_{2,j}(y)
P_{2,k}(z) + \frac{\partial v}{\partial t}(0,y,z,t)- \frac{\partial v}{\partial t}(0,y,0,t)
-\frac{\partial v}{\partial t}(0,0,z,t) + \frac{\partial v}{\partial t}(0,0,0,t)\\
& + \frac{\partial v}{\partial t}(x,0,z,t)- \frac{\partial v}{\partial t}(x,0,0,t) + 
\frac{\partial v}{\partial t}(x,y,0,t), \hspace{0.5cm}t\in [t_s, t_{s+1}).
\end{align}
Now, Integrating the above equation \eqref{3dvt} w.r.t. $t$ from $t_s$ to $t$, we will get 
\begin{align}
\label{v3d}
\nonumber
v(x,y,z,t) &= (t-t_s) \sum_{i,j,k=1}^{2M}\alpha_{i,j,k} P_{2,i}(x) P_{2,j}(y) P_{2,k}(z) + 
v(x,y,z,t_s)+ v(0,y,z,t) -v(0,y,z,t_s) - v(0,y,0,t)\\ \nonumber
&+v(0,y,0,t_s)- v(0,0,z,t) + v(0,0,z,t_s) + v(0,0,0,t)- v(0,0,0,t_s) + v(x,0,z,t) - v(x,0,z,t_s)\\
& - v(x,0,0,t) + v(x,0,0,t_s) + v(x,y,0,t) -  v(x,y,0,t_s).
\end{align}
Again, Integrate \eqref{approxw3d} w.r.t $t$ from $t_s$ to $t$, we acquire
\begin{align}
\label{w3d}
w(x,y,z,t) = (t-t_s)\sum_{l,m,n=1}^{2M}\beta_{l,m,n} h_{l}(x) h_{m}(y) h_{n}(z)+ w(x,y,z,t_s)
\end{align}
To find the solution at the collocation points, we have to discretized the equation \eqref{mv3d} and \eqref{mw3d} when $t \rightarrow t_{s+1}$.

The discrete form is as follows:
\begin{align}
\nonumber
\label{discmv3d}
& \epsilon(x_{k_1},y_{k_2},z_{k_3})\frac{\partial v}{\partial t}(x_{k_1},y_{k_2},z_{k_3},t_{s+1})- \bigg[\sigma _l(x_{k_1},y_{k_2},z_{k_3})\frac{\partial ^2 v}
{\partial x^2}(x_{k_1},y_{k_2},z_{k_3},t_{s+1})+ 
\sigma _{l,x}(x_{k_1},y_{k_2},z_{k_3})\frac{\partial  v}{\partial x}(x_{k_1},y_{k_2},z_{k_3},t_{s+1}) \\ \nonumber
&+(\sigma _t(x_{k_1},y_{k_2},z_{k_3})\frac{\partial ^2 v}{\partial y^2}(x_{k_1},y_{k_2},z_{k_3},t_{s+1})
+\sigma _{t,y}(x_{k_1},y_{k_2},z_{k_3})\frac{\partial  v}{\partial y}(x_{k_1},y_{k_2},z_{k_3},t_{s+1}) 
+\sigma _t(x_{k_1},y_{k_2},z_{k_3})\\ \nonumber
&\frac{\partial ^2 v}{\partial z^2}(x_{k_1},y_{k_2},z_{k_3},t_{s+1})
+\sigma _{t,z}(x_{k_1},y_{k_2},z_{k_3})\frac{\partial  v}{\partial z}(x_{k_1},y_{k_2},z_{k_3},
t_{s+1})\bigg ]
+ f(v(x_{k_1},y_{k_2},z_{k_3},t_{s+1}),w(x_{k_1},y_{k_2},z_{k_3},\\ &t_{s+1})) = 0,\\
\label{discmw3d}
& \frac{\partial w}{\partial t}(x_{k_1},y_{k_2},z_{k_3},t_{s+1})=g(v(x_{k_1},y_{k_2},z_{k_3},
t_{s+1}),w(x_{k_1},y_{k_2},z_{k_3},t_{s+1}).
\end{align}
Using \eqref{approxw3d} at the grid points in \eqref{discmw3d} and linearize the non-linear terms 
by treating it explicitly, we obtain the following
\begin{align*}
&\sum_{l,m,n=1}^{2M}\beta_{l,m,n} h_l(x_{k_1}) h_m(y_{k_2}) h_n(z_{k_3}) 
= g(v(x_{k_1},y_{k_2},z_{k_3},t_{s}),w(x_{k_1},y_{k_2},z_{k_3},t_{s})).
\end{align*}
Matrix system of the above equation is given by,
\begin{align}
\label{matw3d}
H_l H_m H_n \beta=c,
\end{align}
where $H_l, H_m, H_n$ are the Haar matrices and $c^t=(c_{k_1, k_2, k_3})$, which is given by,
\begin{align}
\label{wc3d}
c_{{k_1, k_2, k_3}}= g(v(x_{k_1},y_{k_2},z_{k_3},t_{s}),w(x_{k_1},y_{k_2},z_{k_3},t_{s})).
\end{align}	
Now, at each time step we will calculate the wavelet coefficient $\beta$ and then from \eqref{w3d}
at the collocation points we will calculate the solution $w$. So, now we will use this $w$ to calculate the solution $v$.

Again, Calculate equations \eqref{3dvy2}, \eqref{3dvx2} and \eqref{3dvt} at the collocation points and substitute in \eqref{discmv3d} and \eqref{discmw3d} and treat non-linear terms explicitly in $v$, we get the following equation in matrix form at time $t_{s+1}$:
\begin{align}
\label{matv3d}
K\alpha=b,
\end{align}
where $ K=(k_{ij})$ ia matrix of size $8M^3 \times 8M^3$ and $b$ is a column vector of size $8M^3 
\times 1$.

Now from the above equation we will calculate the wavelet coefficient $\alpha$ and the obtain the 
solution $v$ with the use of calculated $w$, at the desired time step.

\section{Convergence Analysis}

\begin{lemma}  
	If $v(x,y,z)$ and $w(x,y,z)$ are Lipschitz continuous on domain $[0,1]^3$, then the wavelet coefficients
	$a_{i_1,i_2,i_3}$, $b_{l_1,l_2,l_3}$ corresponding to $v$ and $w$ satisfy the inequality 
	\begin{align}
	& |a_{i_1,i_2,i_3}| \leq \frac{L}{8m^4},\\
	& |b_{l_1,l_2,l_3}| \leq \frac{L}{8m^4},
	\end{align}
	where $L>0$ depends on the Lipschitz constant and the coefficients are defined as
	\begin{align}
	\label{coefa}
	& a_{i_1,i_2,i_3} = \int_0^1 \int_0^1 v(x,y,z) h_{i_1}(x) h_{i_2}(y) h_{i_3}(z)dx dy dz,\\
	\label{coefb}
	& b_{l_1,l_2,l_3} = \int_0^1 \int_0^1 w(x,y,z) h_{l_1}(x) h_{l_2}(y) h_{l_3}(z)dx dy dz,
	\end{align}
\end{lemma}
Proof: We can redefine the definition of the coefficient $a_{i_1,i_2,i_3}$ in terms of the inner product as follows 
\begin{align*}
a_{l_1,l_2,l_3} =\int_0^1 \int_0^1 v(x,y,z) h_{i_1}(x) h_{i_2}(y) h_{i_3}(z)dx dy dz
= \big < h_{i_1}(x), \big < h_{i_2}(y), \big < v(x,y,z),h_{i_3}(z) \big > \big> \big >,
\end{align*}
where \big< . \big> is defined as the inner product.
\begin{align*}
\big < v(x,y,z),h_{i_3}(z) \big > &= \int_0^1 v(x,y,z)h_{i_3}(z) dz.
\end{align*}
Using the definition of Haar wavelet, we get
\begin{align*}
\big < v(x,y,z),h_{i_3}(z) \big > &= \int_{\frac{k}{m}}^{\frac{k+0.5}{m}} v(x,y,z) dz
- \int_{\frac{k+0.5}{m}}^{\frac{k+1}{m}} v(x,y,z) dz.
\end{align*}

Now, applying the Mean Value theorem for both the integrals, we find an $z_1 \in 
\bigg(\frac{k}{m},\frac{k+0.5}{m} \bigg)$ and an $z_2 \in 
\bigg(\frac{k+0.5}{m},\frac{k+1}{m} \bigg)$ such that
\begin{align*}
\big < v(x,y,z),h_{i_3}(z) \big > &= \bigg(\frac{k+0.5}{m}-\frac{k}{m} \bigg)v(x,y,z_1)
-\bigg(\frac{k+1}{m}-\frac{k+0.5}{m} \bigg)v(x,y,z_2)\\
&=\frac{1}{2m} [v(x,y,z_1)-v(x,y,z_2)].
\end{align*}
Now, 
\begin{align*}
\big < h_{i_2}(y), \big < v(x,y,z),h_{i_3}(z) \big > \big> = \frac{1}{2m}\bigg\{ \int_0^1v(x,y,z_1)h_j(y)dy
- \int_0^1 v(x,y,z_2)h_j(y)dy \bigg\}.
\end{align*}
Again, using the definition of Haar wavelet in both the integral, we obtain

\begin{align*}
\big < h_{i_2}(y), \big < v(x,y,z),h_{i_3}(z) \big > \big> = \frac{1}{2m}\bigg[ \bigg\{ \int_{\frac{k}{m}}^{
	\frac{k+0.5}{m}}v(x,y,z_1)dy - \int_{\frac{k+0.5}{m}}^{\frac{k+1}{m}}v(x,y,z_1) dy \bigg\}\\
-\bigg\{ \int_{\frac{k}{m}}^{\frac{k+0.5}{m}}v(x,y,z_2)dy - 
\int_{\frac{k+0.5}{m}}^{\frac{k+1}{m}}v(x,y,z_2) dy \bigg\}\bigg].
\end{align*}
Now, applying the Mean Value Theorem for all the four integrals, we find $y_1 ,y_2 \in 
\bigg(\frac{k}{m},\frac{k+0.5}{m} \bigg)$ and $y_3, y_4 \in 
\bigg(\frac{k+0.5}{m},\frac{k+1}{m} \bigg)$ such that

\begin{align*}
\big < h_{i_2}(y), \big < v(x,y,z),h_{i_3}(z) \big > \big> = \frac{1}{4m^2}\bigg[v(x,y_1,z_1)
-v(x,y_2,z_1)-v(x,y_3,z_2)+ v(x,y_4,z_2) \bigg].
\end{align*}
Hence,
\begin{align*}
a_{i_1,i_2,i_3} &= \big < h_{i_1}(x), \big < h_{i_2}(y), \big < v(x,y,z),h_{i_3}(z) \big > \big> \big >\\
&= \frac{1}{4m^2}\bigg\{ \int_0^1v(x,y_1,z_1)h_{i_1}(x)dx - \int_0^1v(x,y_2,z_1)h_{i_1}(x)dx
- \int_0^1 v(x,y_3,z_2)h_{i_2}(x)dx + \int_0^1 v(x,y_4,z_2)h_{i_2}(x)dx \bigg\}
\end{align*}
Using the definition of Haar wavelet in all the integral and then the Mean Value Theorem, we
have we find $x_1, x_2, x_3, x_4 \in \bigg(\frac{k}{m},\frac{k+0.5}{m} \bigg)$ and
$x_5, x_6, x_4 , x_8 \in \bigg(\frac{k+0.5}{m},\frac{k+1}{m} \bigg)$ such that

\begin{align*}
|a_{i_1,i_2,i_3}| \leq \frac{1}{8m^3}\bigg[ |v(x_1,y_1,z_1)-v(x_2,y_1,z_1)|+|v(x_3,y_2,z_1)-v(x_4,y_2,z_1)|\\
+|v(x_5,y_3,z_2)-v(x_6,y_3,z_2)|+ |v(x_7,y_4,z_2)-v(x_8,y_4,z_2)| \bigg].
\end{align*}

Now, using the Lipschitz continuity of $v$, we get 
\begin{align*}
|a_{i_1,i_2,i_3}| \leq \frac{1}{8m^3} \{L_1 |x_1-x_2| + L_2 |x_3-x_4| + L_3 |x_5-x_6| + L_4 |x_7-x_8|\},
\end{align*}
where, $L_1, L_2, L_3, L_4$ are the Lipschitz constants.

Choose $L=max(L_1, L_2, L_3, L_4)$ and from the definition of Haar wavelet we can conclude that
$|x_i-x_j|\leq \frac{1}{2m}$, we obtain
\begin{align*}
|a_{i_1,i_2,i_3}| &\leq \frac{1}{8m^3} \bigg(\frac{4L}{2m}\bigg)\\
&=\frac{L}{8m^4}.
\end{align*}
Similarly, the wavelet coefficient $b_{l_1,l_2,l_3}$ can also be determined. \\ \\
$\bullet$ Now, introducing the norm as follows:
\begin{align}
\label{norm}
\parallel u \parallel_X = \parallel v \parallel_2 + \sum_{k=1}^{d}\parallel w^k \parallel_2,
\end{align}

where $ \parallel . \parallel_2 $ is the $L^2([0,1]^3)$-norm. 

Let  $ u =\begin{bmatrix} v & w \end{bmatrix}$ be the exact solution of the problem and
$u_H =\begin{bmatrix} v_H & w_H \end{bmatrix}$ be the solution approximated by the Haar wavelets. 
The error $(E=u-u_H)$ is given as follows:

\begin{theorem}
	Let $u$ be the exact solution of the problem and $u_H$ be the solution approximated by the Haar wavelets, then
	\begin{align}
	\parallel E \parallel_X ^2 = \parallel u(x,y,t_{n+1})-u_H(x,y,t_{n+1}) \parallel_X ^2
	\leq C dt^2 \bigg[ \frac{K_1 K_2 K_3}{m^2} + \frac{d}{m^8} \bigg], 
	\end{align}
	which implies that $\parallel E \parallel_X ^2 \to 0$ as $dt \to 0$ and $m \to \infty$.
\end{theorem}
Proof.Calculate the error using equations \eqref{v3d} and \eqref{w3d} for the Haar wavelet solution, we obtain the following
\begin{align*}
\parallel E \parallel_X ^2 &=  \parallel dt \sum_{i_1,i_2,i_3=2M+1}^{\infty} a_{i_1,1_2,i_3} P_{2,i_1}(x) P_{2,i_2}(y)  P_{2,i_3}(z)\parallel_2 ^2 +  \sum_{w=1}^d \parallel dt \sum_{l_1,l_2,l_3=2M+1}^{\infty} b_{l_1,l_2,l_3}^w H_{l_1}(x) H_{l_2}(y) H_{l_3}(z) \parallel_2 ^2 \\
& = dt^2 \int_{0}^{1} \int_{0}^{1} |\sum_{i_1,i_2,i_3=2M+1}^{\infty} a_{i_1,1_2,i_3} P_{2,i_1}(x) P_{2,i_2}(y) P_{2,i_3}(z) \sum_{p,q,r=2M+1}^{\infty} a_{p,q,r} P_{2,p}(x) P_{2,q}(y)P_{2,r}(z)| \\
&+ \sum_{w=1}^d \parallel dt^2 \int_{0}^{1} \int_{0}^{1} | \sum_{l_1,l_2,l_3=2M+1}^{\infty} b_{l_1,l_2,l_3}^w H_{l_1}(x) H_{l_2}(y)H_{l_3}(z) \sum_{s,t,h=2M+1}^{\infty} b_{s,t,h}^k H_{s}(x) H_{t}(y) H_{h}(z)|,
\end{align*}
\begin{align*}
& \leq dt^2 \sum_{i_1,i_2,i_3=2M+1}^{\infty} |a_{i_1,i_2,i_3} |  \sum_{p,q,r=2M+1}^{\infty} |a_{p,q,r} | \bigg(\int_{0}^{1}|P_{2,i_1}(x) P_{2,p}(x)|dx \bigg)  \bigg(\int_{0}^{1}|P_{2,i_2}(y) P_{2,q}(y)|dy\bigg)\\ &\bigg(\int_{0}^{1}|P_{2,i_3}(z) P_{2,r}(z)|dz\bigg)+ dt^2 \sum_{w=1}^d \bigg[\sum_{l_1,l_2,l_3=2M+1}^{\infty} |b_{l_1,l_2,l_3}^w| \sum_{s,t,h=2M+1}^{\infty} |b_{s,t,h}^w| \bigg( \int_{0}^{1} |H_{l_1}(x) H_{s}(x)|dx \bigg)\\ &\bigg(\int_{0}^{1} |H_{l_2}(y) H_{t}(y)|dy\bigg) \bigg(\int_{0}^{1} |H_{l_3}(z) H_{h}(z)|dz\bigg) \bigg]\\
&= dt^2 \sum_{i_1,i_2,i_3=2M+1}^{\infty} |a_{i_1,i_2,i_3}| \sum_{p,q,r=2M+1}^{\infty} |a_{p,q,r}| K_{i_1,p}K_{i_2,q} K_{i_3,r} + dt^2 \sum_{w=1}^d \bigg[\sum_{l_1,l_2,l_3=2M+1}^{\infty} |b_{l_1,l_2,l_3}^w| \sum_{p,q,r=2M+1}^{\infty} |b_{p,q,r}^k| L_{l_1,s} \\
&L_{l_2,t} L_{l_3,h} \bigg],\\
&= I_1+I_2,
\end{align*}
where,
\begin{align*}
&K_{i_1,p}=\int_{0}^{1}|P_{2,i_1}(x) P_{2,p}(x)|dx,\\  &K_{i_2,q}=\int_{0}^{1}|P_{2,i_2}(y) P_{2,q}(y)|dy,\\ 
&K_{i_3,r}=\int_{0}^{1}|P_{2,i_3}(z) P_{2,r}(z)|dz,
\end{align*}
and 
\begin{align*}
&L_{l_1,s}= \int_{0}^{1} |H_{l_1}(x) H_{p}(x)|dx, \\
&L_{l_2,t}= \int_{0}^{1} |H_{l_2}(y) H_{q}(y)|dy,\\
&L_{l_3,h}= \int_{0}^{1} |H_{l_3}(z) H_{h}(z)|dz.
\end{align*}

Let $K_{i_1} = \sup_r K_{i_1,p}$, $K_{i_2} = \sup_s K_{i_2,q}$, and,  $K_{i_3} = \sup_s K_{i_3,r}$ first term $(I_1)$ becomes as 

\begin{align*}
I_1 & \leq  dt^2 \sum_{i_1,i_2,i_3=2M+1}^{\infty} |a_{i_1,i_2,i_3}|  K_{i_1}K_{i_2}K_{i_3} \sum_{p,q,r=2M+1}^{\infty} |a_{p,q,r}| \\
& \leq  \frac{Ldt^2}{8} \sum_{i_1,i_2,i_3=2M+1}^{\infty} |a_{i_1,i_2,i_3}|  K_{i_1}K_{i_2}K_{i_3} \sum_{j=J+1}^ \infty \sum_{p,q,r}^{2^{j}-1} \frac{1}{m^4}, \hspace{2mm} \text{using Lemma 3} \\
& \leq  \frac{Ldt^2}{8} \sum_{i_1,i_2,i_3=2M+1}^{\infty} |a_{i_1,i_2,i_3}|  K_{i_1}K_{i_2}K_{i_3} \frac{1}{m} \\
\end{align*}

Now, let $K_1 = \sup_{i_1} K_{i_1}$, $K_2 = \sup_{i_2} K_{i_2}$, and, $K_3 = \sup_{i_3} K_{i_3}$ we obtain
\begin{align*}
I_1 & \leq  \frac{Ldt^2}{8}  K_1 K_2 K_3 \sum_{i_1,i_2,i_3=2M+1}^{\infty} |a_{i_1,i_2,i_3}|  \frac{1}{m} \\
& \leq  \frac{L^2dt^2}{64}  K_1 K_2 K_3 \sum_{j=J+1}^ \infty \sum_{i_1,i_2,i_3=0}^{2^{j}-1} \frac{1}{m^5}, \hspace{2mm} \text{using Lemma 3} \\
& \leq  \frac{L^2dt^2 K_1 K_2 K_3}{64m^2} \\
\end{align*}

Again, we know by the definition of Haar wavelet ,  
\begin{align*}
\int_{0}^{1} H_{l_1}(x) H_{p}(x)dx =  2^{-j} =\frac{1}{m}
\end{align*}

\begin{align*}
I_2 &= (dt)^2 \sum_{w=1}^d \bigg[\sum_{l_1,l_2,l_3=2M+1}^{\infty} |b_{l_1,l_2,l_3}^{q1}|^2 \frac{1}{m^3} \bigg]\\
& \leq \frac{L^2dt^2}{64}\sum_{w=1}^d \bigg(\frac{1}{m^4}\bigg)^2, \hspace{2mm} \text{using Lemma 3}\\
& \leq \frac{dL^2dt^2}{16m^8}
\end{align*}

Therefore,
\begin{align*}
\parallel E \parallel_X ^2 & \leq \frac{C^2dt^2 K_1 K_2K_3}{64m^2} + \frac{dC^2dt^2}{64m^8} \\
& \leq C dt^2 \bigg[ \frac{K_1 K_2 K_3}{m^2} + \frac{d}{m^8} \bigg]
\end{align*}

\paragraph{{Algorithm to calculate the Haar wavelet solution}}

\begin{itemize}
	\item  Step1: Calculate the Haar matrix H and the Haar Integral matrix $P_1$ and $P_2$ by using the equations
	\eqref{haar} and \eqref{HaarIntegral} respectively. 
	
	\item Step2: Compute c using the equation \eqref{wc3d} and then solve the matrix
	equation $H \alpha = c$ to find the wavelet coefficient $\beta$. 
	
	\item Step3: Use the wavelet coefficient $\beta$ in equation \eqref{w3d} and the given initial values of
	$v$ and $w$ at the to compute the gating variable $w$. 
	
	\item Step4: Compute the matrix $K$ and the vector $b$. Use the above calculated gating variables $w$ and solve the matrix equation $K \alpha = b$ to find the wavelet coefficient $\alpha$ and then use $\alpha$
	in \eqref{v3d} to compute the solution $v$.
\end{itemize}

\textbf{Remark:} (1) Theorem 1 implies the stability of the solution and theorem 2 implies that the truncation error goes to zero and hence in conjunction they implies the convergence of the solution.

(2) If there is discontinuity in the parameter of the function $f(v,w)$, $g(v,w)$ or the coefficient $\epsilon$, $D(x,y)$ , it can be handled using Haar wavelet interpolation.

\section{Numerical Result and Discussions}
To demonstrate the power of Haar wavelet method we presented a few examples which are related to the pathological case studies in cardiac electrical activity.
We solve all the examples using the above developed Haar wavelet method and calculate the absolute error also. Grid validation test
or resolution level test has been done for all the problems and here we are presenting for some of the problems. From grid 
validation we observe that resolution level $J=4$ in two dimension is good enough to calculate the solution.
We also compare the Haar wavelet solution with the solution obtained using the finite element method for almost all the problems. 
We also calculate the CPU time for Haar wavelet method and the finite element method.
For all the examples, we use the GMRES solver to solve the linear system of equations obtained from the discrete problem as it is 
more appropriate in terms of the CPU time in comparison to the other solvers like, CGS, BICG and BICGSTAB in terms of the CPU time (see Table \eqref{cputime_pcg}).

\paragraph{}
\textbf{Example $1$.} We consider the one dimensional reaction-diffusion system coupled with the ODE having homogeneous Neumann boundary. Such type of coupled systems are important because it models the problems related to the field of cardiac electrical activity of any species. We will also discuss the problem with jump discontinuity in the parameter or coefficients which is related to the problem of ischemia in cardiac tissue. 
\begin{align*}
\epsilon \frac{\partial v}{\partial t}- d\frac{d}{d x}(\frac{d v}{d x}) + kv(v-0.1)(1-v)-kw  &= I_{app}, & 0\leq x \leq 1,  0\leq t \leq T\\
\frac{\partial w}{\partial t}&=v-2w, & 0\leq x \leq 1,  0\leq t \leq T\\
v(x,0)= 0.2, \hspace{5mm} w(x,0)&=0.2, & 0\leq x \leq 1\\
\frac{d v}{d x}(0,t)=0,\hspace{0.5cm} \frac{d v}{d x}(1,t)=0 &=0, & 0\leq t \leq T\\
\end{align*}
where, $d=0.005$ and $\epsilon=0.01$. 

\begin{figure}[h] 
	\centering
	\vspace{-18em}
	\includegraphics[width=0.9\textwidth]{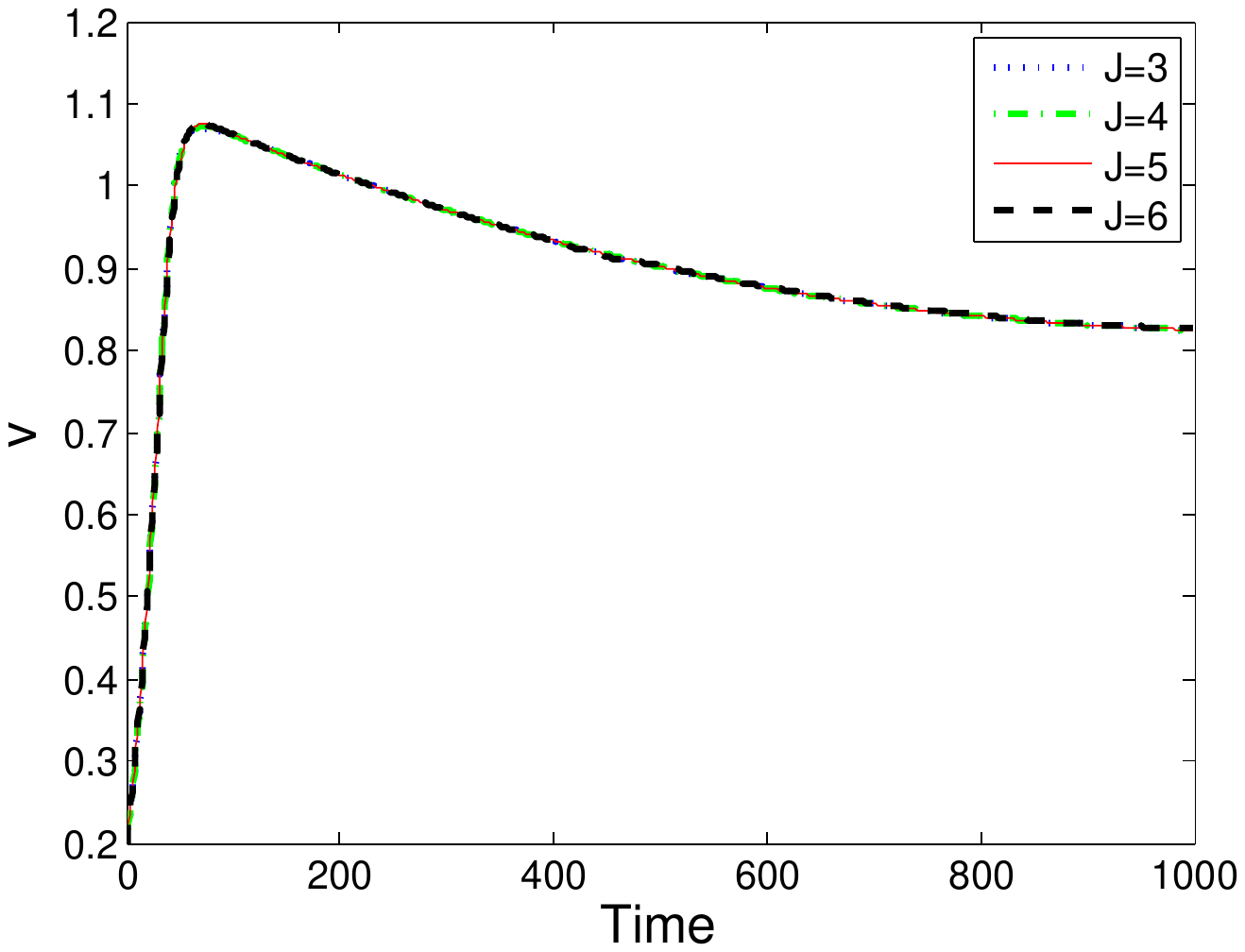}
	\vspace{-18em}
	\caption{\textbf{Resolution level test for Haar wavelet solution $v$}}
	\label{gridValid_fhn1D}
\end{figure}

Resolution level test for the proposed Haar wavelet method has been presented in Fig. \ref{gridValid_fhn1D}.

\begin{figure}[h] 
	\centering
	\vspace{-13em}
	\includegraphics[width=0.6\textwidth]{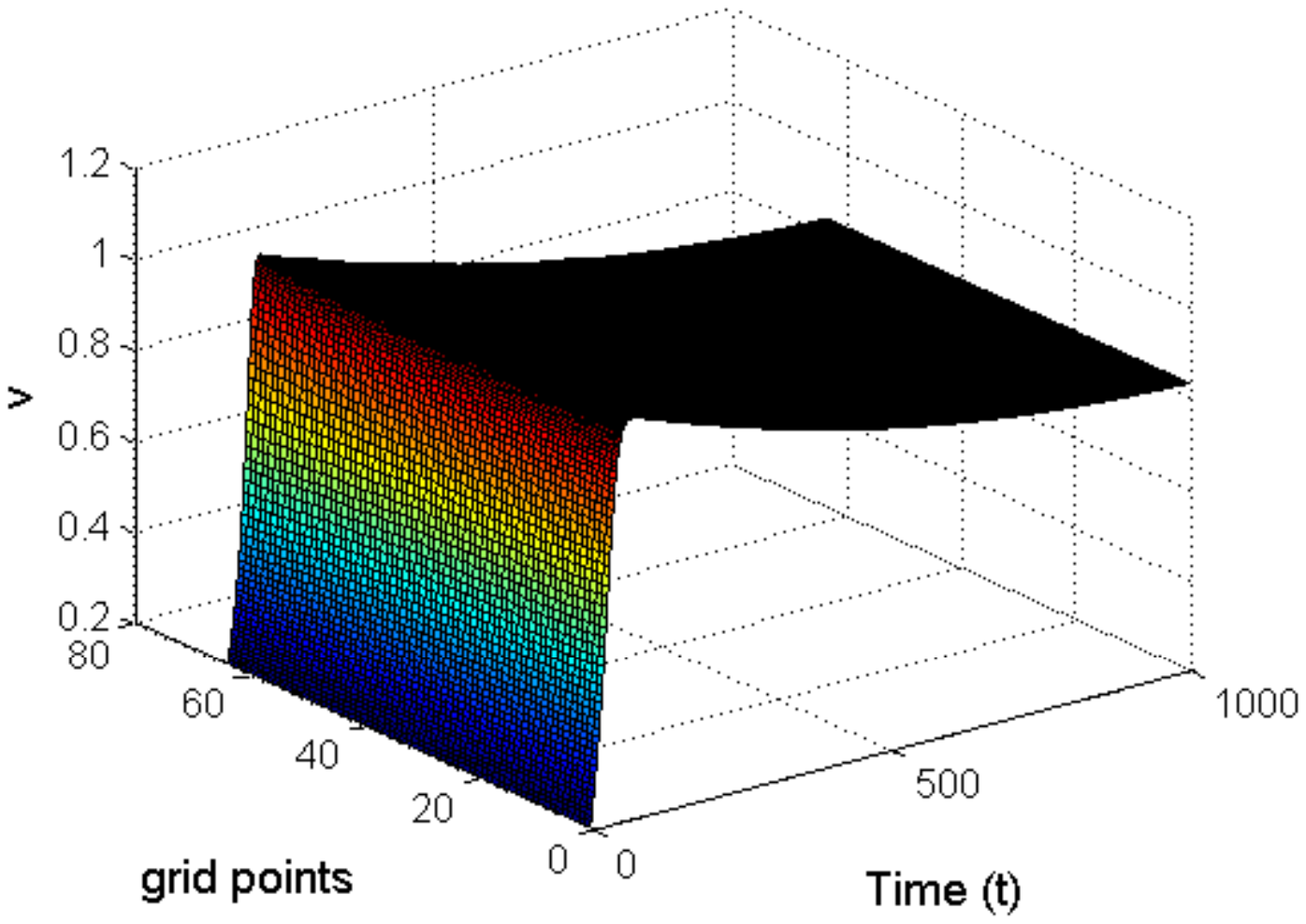}
	\vspace{-9em}
	\caption{\textbf{(for example1) Haar wavelet solution for $v$. $I_{app}=0.3$.}}
	\label{HW}
\end{figure}

\begin{figure}[h]
	\vspace{-10em}
	\centering
	\includegraphics[width=0.6\textwidth]{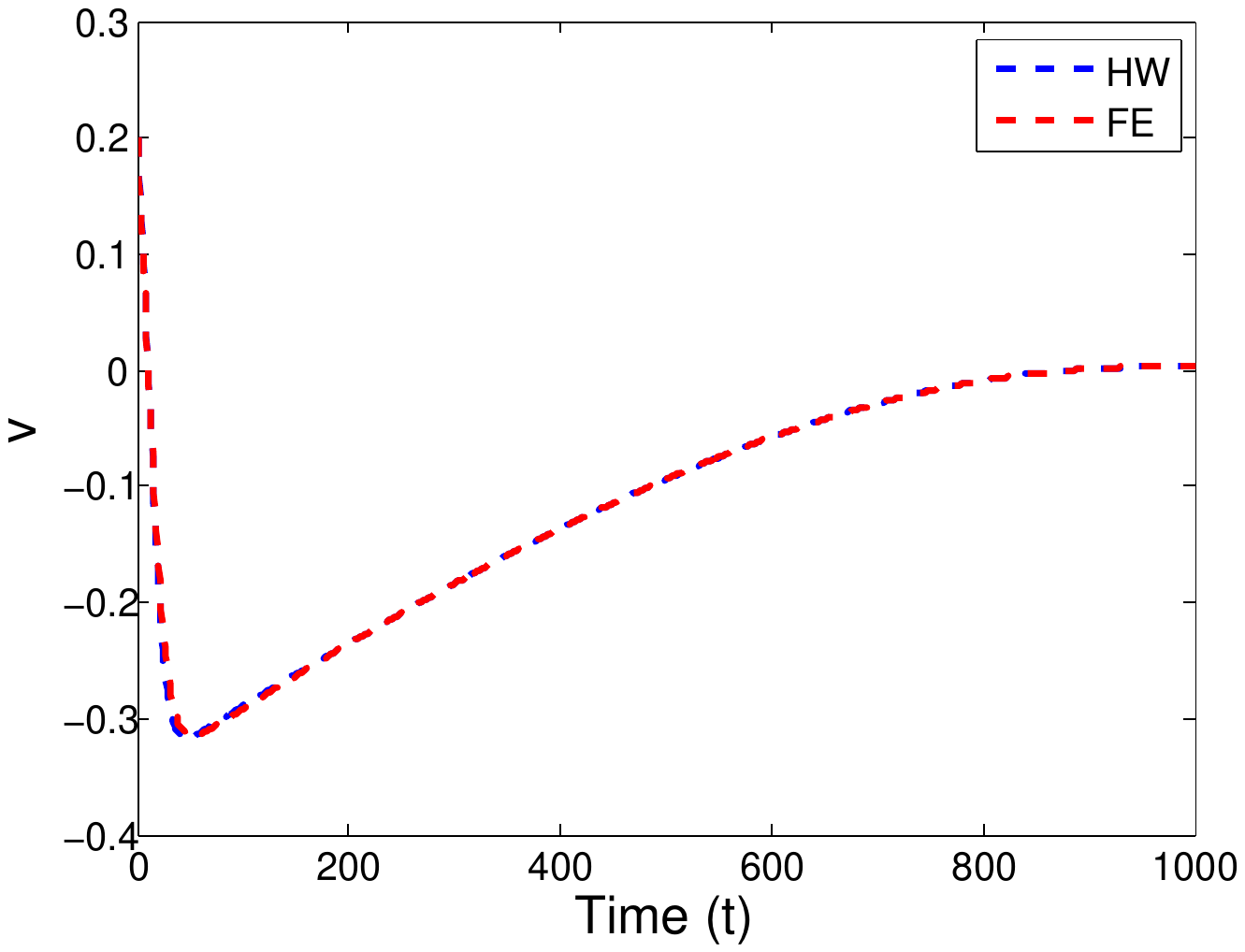}
	\vspace{-10em}
	\caption{\textbf{FE and HW solution at x=0.5. $I_{app}=0$, J=5, $dt=10^{-3}$ and T=1.}}
	\label{fhncomparison}
\end{figure}

\paragraph{}
Pointwise absolute error for different time steps is shown in Table \ref{table1}.
Solution at $dt=10^{-5}, J=8$ is taken as the reference solution. From the Table \ref{table1},
it can be seen clearly that absolute error decreases significantly with the smaller time step size.
The Haar wavelet solution for $v$ at the grid points is presented in Fig. \ref{HW}. 
We also solved the same problem using the linear finite element in space and implicit-explicit (IE) Euler method in time. It is observed that the Haar wavelet solution matches with the finite element solution (Fig. \ref{fhncomparison}). 

\begin{table}[h]
	\begin{center}
		\begin{tabular}{ |c  c  c | } 
			\hline
			x &  absolute error & \\
			\hline 
			& $dt=10^{-3}$ & $dt =10^{-4}$\\
			\hline
			0.1406 & $6.4 \times10^{-3}$ & $6.2 \times 10^{-4}$\\ 
			
			0.2656 & $6.9 \times10^{-3}$ & $6.085 \times 10^{-4}$\\ 
			
			0.3594 & $8.1 \times10^{-3}$ & $5.464 \times 10^{-4}$\\ 
			
			0.4531 & $8.8 \times10^{-3}$ & $5.4 \times 10^{-4}$\\ 
			\hline
		\end{tabular}
		\caption{(For example 1)Absolute Error for $v$ at different points of the domain when $J=8$ and T=1.}
		\label{table1}
	\end{center}
\end{table} 

From Fig. \ref{fhncomparison} we can see that when $I_{app}=0$, variable $v$ goes to the resting state without 
shooting up. But as we provide $I_{app}=0.3$, the solution $v$ shoot up from the initial state and then it goes down slowly to 
the stabilizing state. This type of problems models pathological phenomena like FitzHugh Nagumo (FHN) model in cardiac electrical activity, which is the reduced ionic model of animal. 

\begin{figure}[h]
	\vspace{-15em}
	\centering
	\begin{subfigure}[t]{0.5\textwidth}
		\centering
		\includegraphics[width=1.2\textwidth]{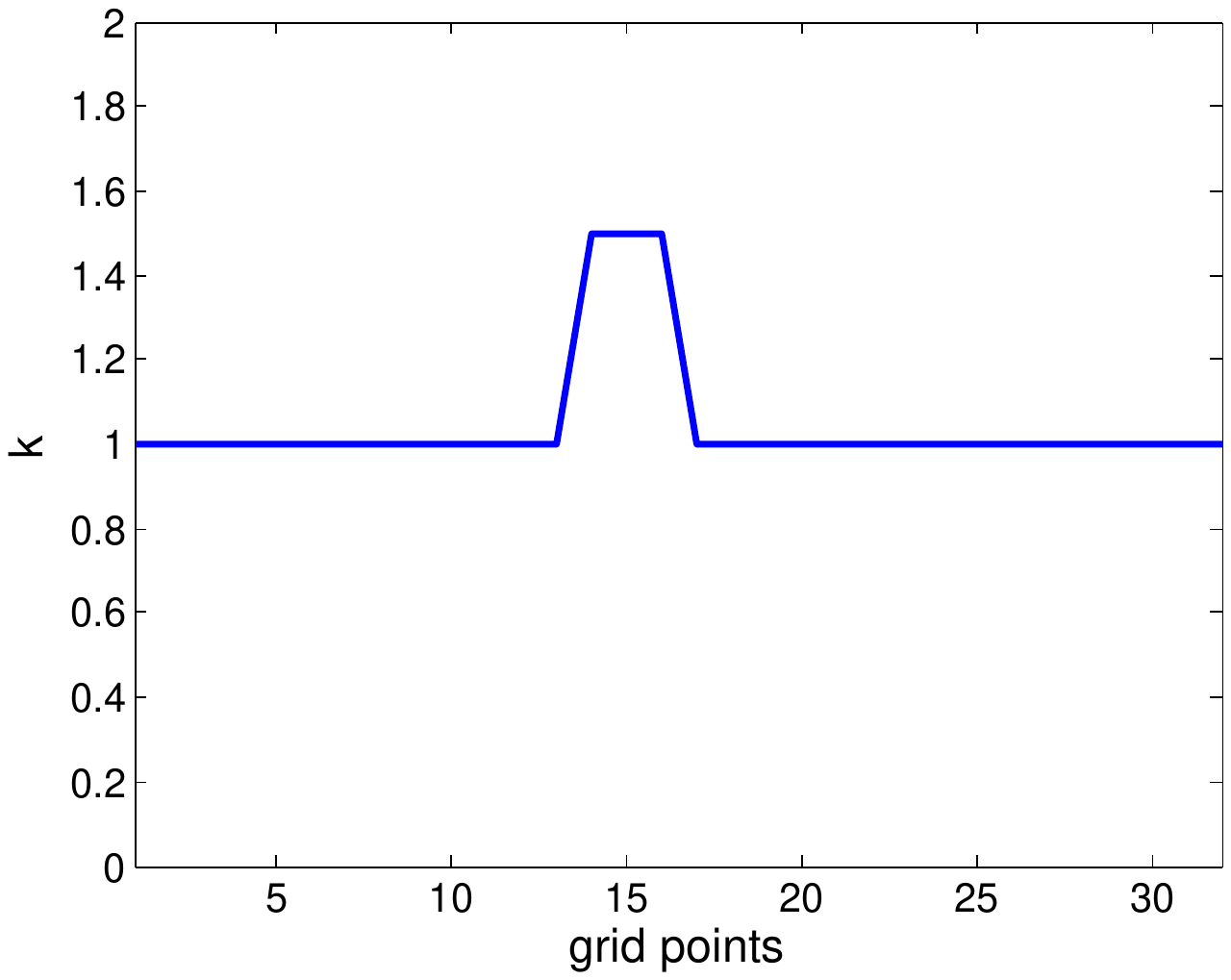}
		\vspace{-13em}
		\caption{}
		\label{k_1djump}
	\end{subfigure}\hfill
	\begin{subfigure}[t]{0.5\textwidth}
		\centering
		\includegraphics[width=1.2\textwidth]{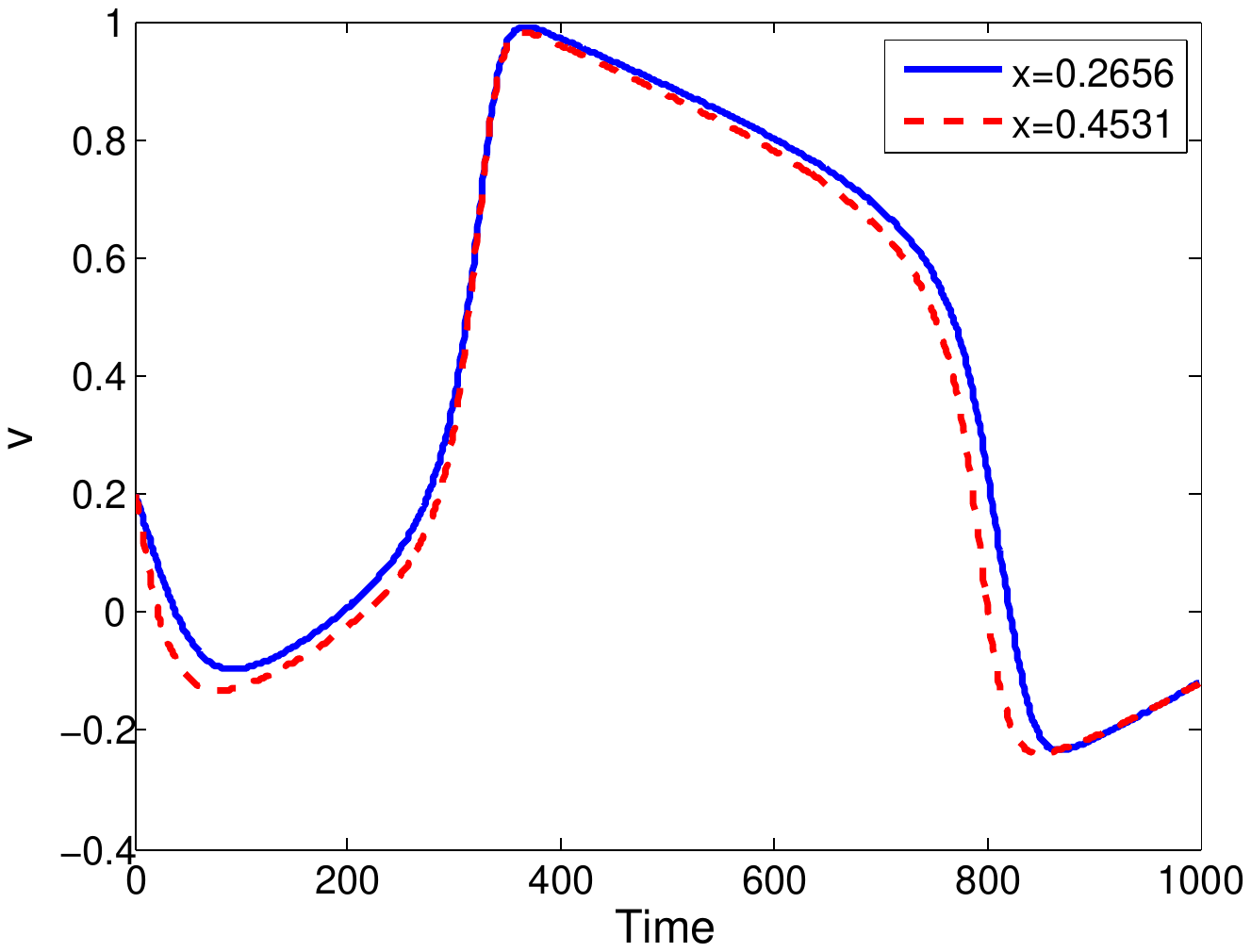}
		\vspace{-13em}
		\caption{}
		\label{AP_fhn_1d_jump_k}
	\end{subfigure}
	\vspace{-10em}
	\caption{\textbf{(a) parameter $k$, (b) Corresponding Haar wavelet solution $v$ at the points falling inside and outside the jump region. $I_{app}=0.15$}}
	\label{FHN_jump_k_1d}
\end{figure}

Next, we consider the same problem but with a jump discontinuity in the parameter $k$, drawn in Fig. \ref{k_1djump} which creates the local ischemia in a cardiac tissue. 
We calculate the solution $v$ at the points, as a sample, in Fig \eqref{AP_fhn_1d_jump_k} we present the solution corresponding to the points $x=0.2656$ and $x=0.4531$, which lie outside and inside the jump region respectively. The absence of any spurious oscillation in the solution presented in Fig. \ref{AP_fhn_1d_jump_k} indicated that Haar wavelets efficiently handle this jump discontinuity.

\paragraph{Example $2$} We will consider the two dimensional problems which are important in the field of cardiac electrophysiology. We will calculate the transmembrane potential in the heart which is important to know the functionality of the heart. We will solve the problem having jump discontinuity in any parameter or coefficients using the Haar wavelets.
\begin{align*}
&\epsilon \frac{\partial v}{\partial t}- div(D\nabla v) + kv(v-0.1)(1-v)-kw = I_{app},&  0\leq x,y \leq 1,  0\leq t \leq T\\
&\frac{\partial w}{\partial t} = v-2w, & 0\leq x,y \leq 1,  0\leq t \leq T\\
&v(x,0)= 0.2, \hspace{5mm} w(x,0)=0.2, & 0\leq x,y \leq 1\\
&\frac{d v}{d x}(0,y,t)=0,& 0\leq t \leq T\\
&\frac{d v}{d x}(1,y,t) =0, & 0\leq t \leq T\\
&\frac{d v}{d y}(x,0,t)=0,& 0\leq t \leq T\\
&\frac{d v}{d y}(x,1,t) =0,& 0\leq t \leq T
\end{align*}
where $ D = \begin{bmatrix} 1.2e{-3} & 0 \\ 0 &  2.5562e{-4}\end{bmatrix}$, $\epsilon=0.01$ and $I_{app}=0.15$.

\begin{figure}[h] 
	\centering
	\vspace{-14em}
	\includegraphics[width=0.6\textwidth]{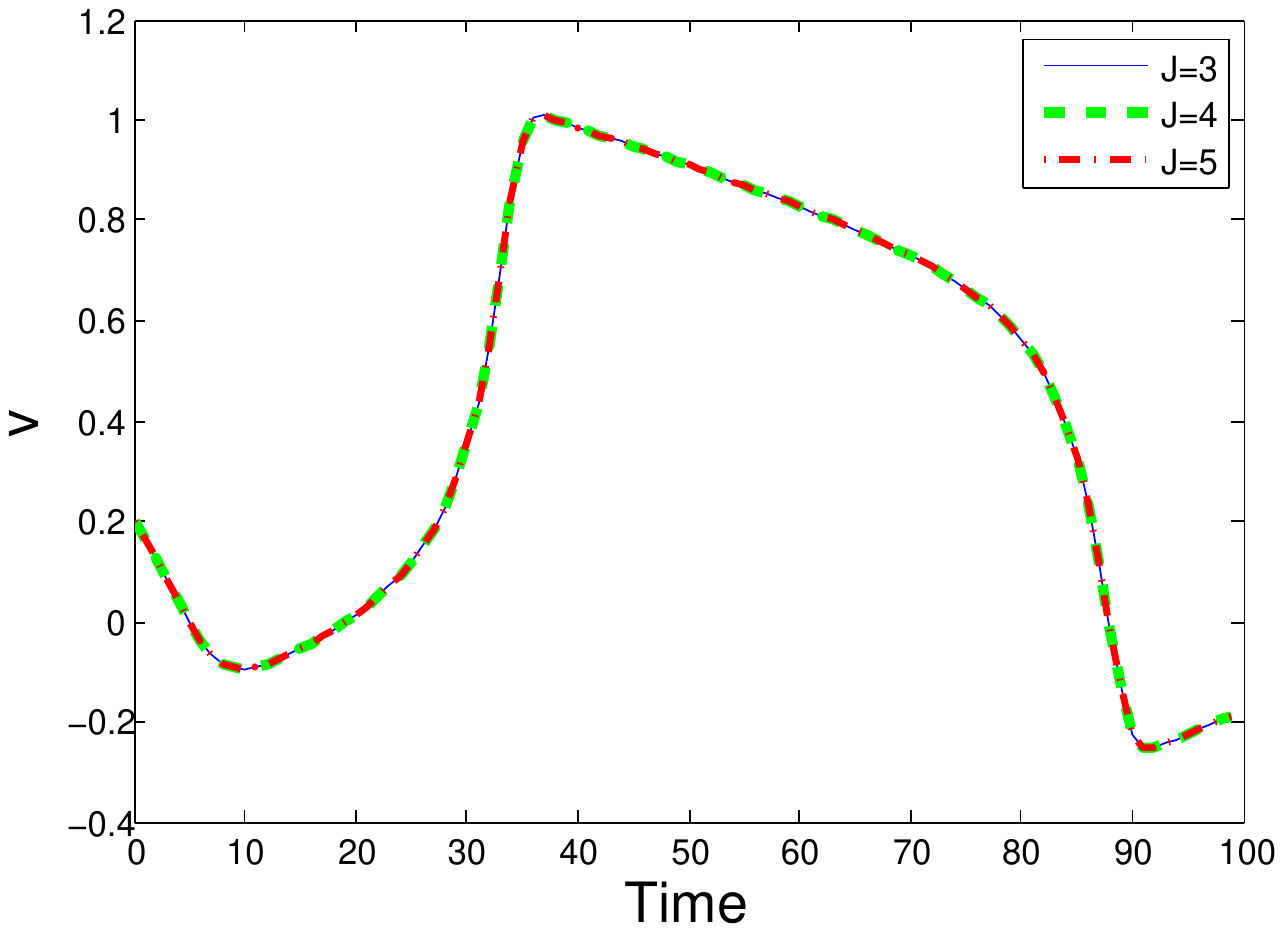}
	\vspace{-13em}
	\caption{\textbf{ Resolution level test for Haar wavelet solution $v$}}
	\label{gridValid_fhn2D}
\end{figure}

\paragraph{}
First of all the grid validation of the proposed algorithm for this problem is presented in Fig. \ref{gridValid_fhn2D} which clearly shows the accuracy of the solution at the different resolution level. So, resolution level $J=4$ is good enough to calculate the results.

Absolute error for different time steps is shown in Table \ref{table3}.
Solution at $dt=10^{-5}, J=6$ is taken as the reference solution. From the Table \ref{table3}
it can be seen clearly that absolute error decreases seriously with the smaller time step size.

\begin{figure}[h]
	\vspace{-13em}
	\centering
	\begin{subfigure}[t]{0.5\textwidth}
		\centering
		\includegraphics[width=1.2\textwidth]{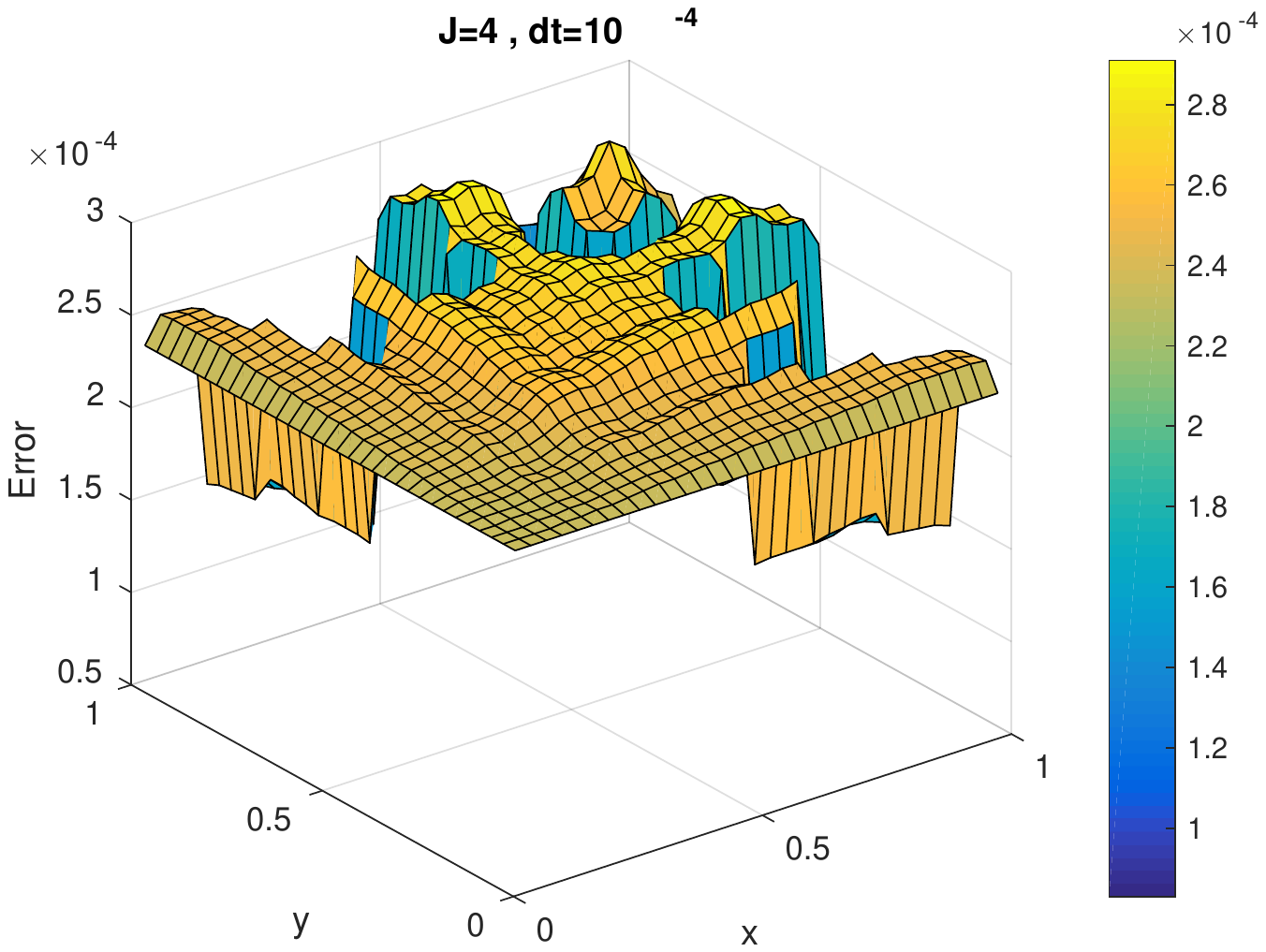}
		\vspace{-12em}
		\caption{}
		\label{error2d_fhn}
	\end{subfigure}\hfill
	\begin{subfigure}[t]{0.5\textwidth}
		\centering
		\includegraphics[width=1.2\textwidth]{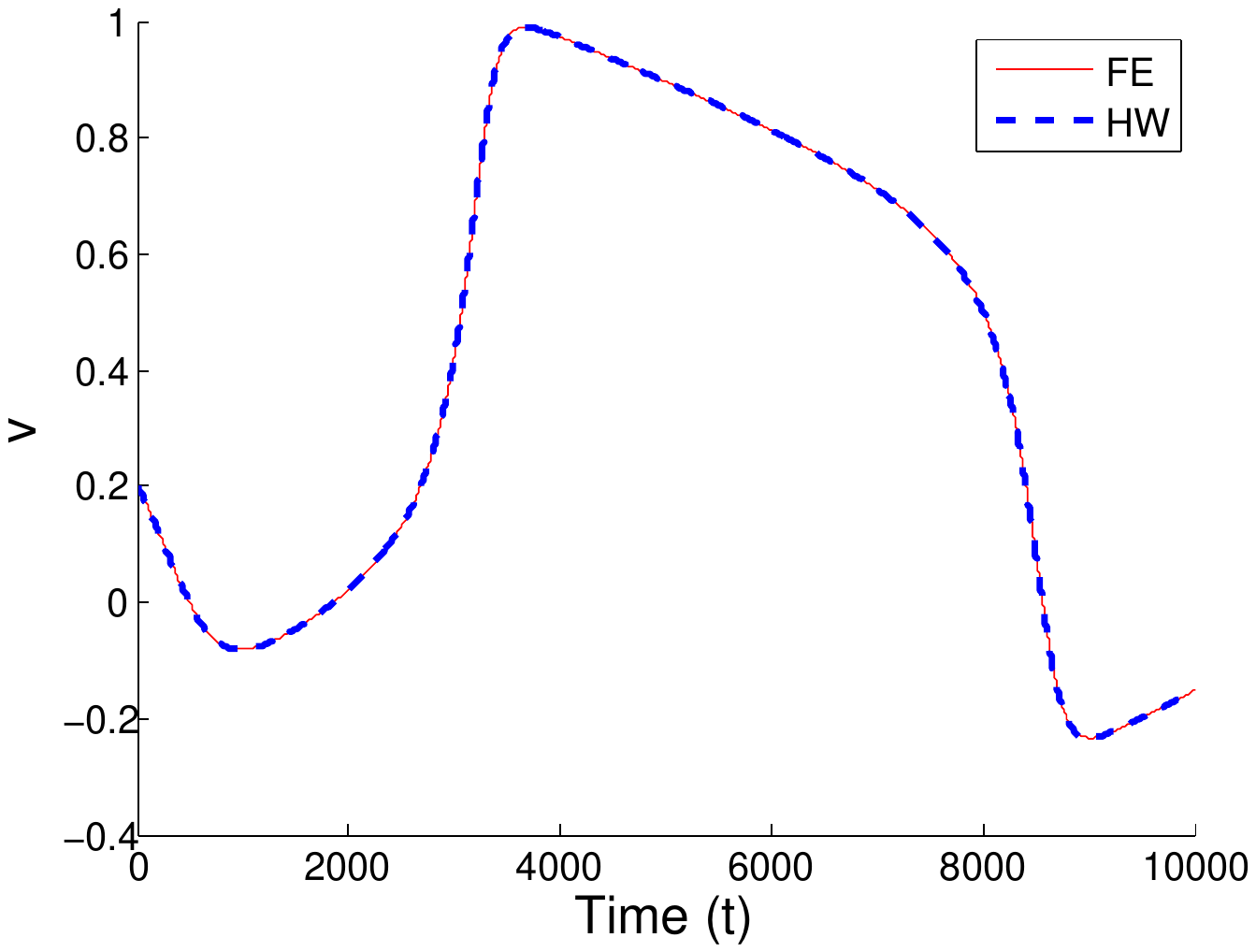}
		\vspace{-12em}
		\caption{}
		\label{fe_hw_2d_compariosn}
	\end{subfigure}
	\vspace{-11em}
	\caption{\textbf{(a) Numerical error at J=6 and T=0.5 and $dt=10^{-4}$, (b) FE and HW solution comparison at J=5 and T=1 and $dt=10^{-4}$, $k=1$.}}
	\label{2D_FHN}
\end{figure}

\begin{figure}
	\vspace{-10em}
	\centering
	\begin{subfigure}[t]{0.5\textwidth}
		\centering
		\includegraphics[width=1\textwidth]{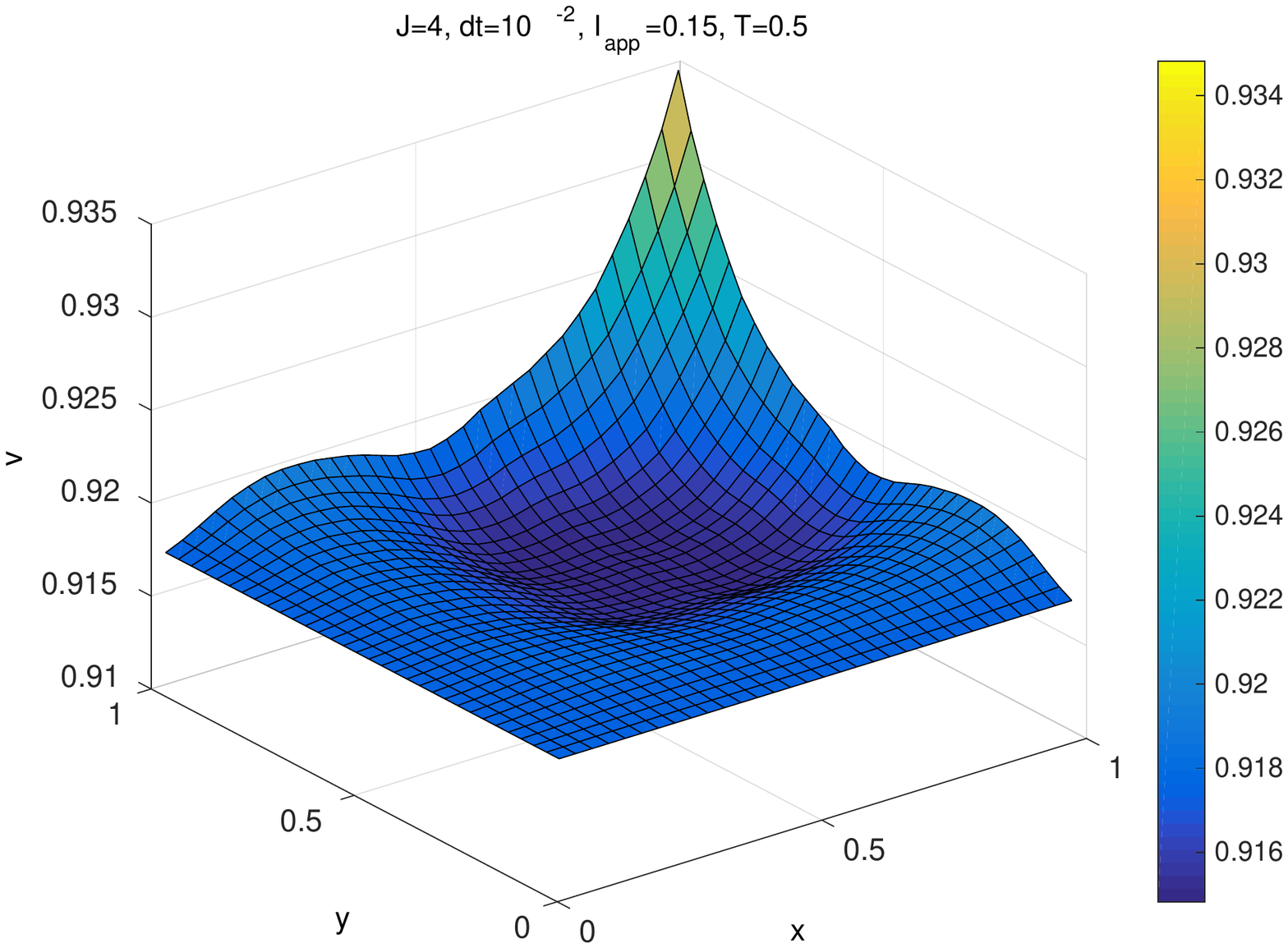}
		\vspace{-12em}
		\caption{}
		\label{surf2d_xy}
	\end{subfigure}\hfill
	\begin{subfigure}[t]{0.5\textwidth}
		\centering
		\includegraphics[width=1\textwidth]{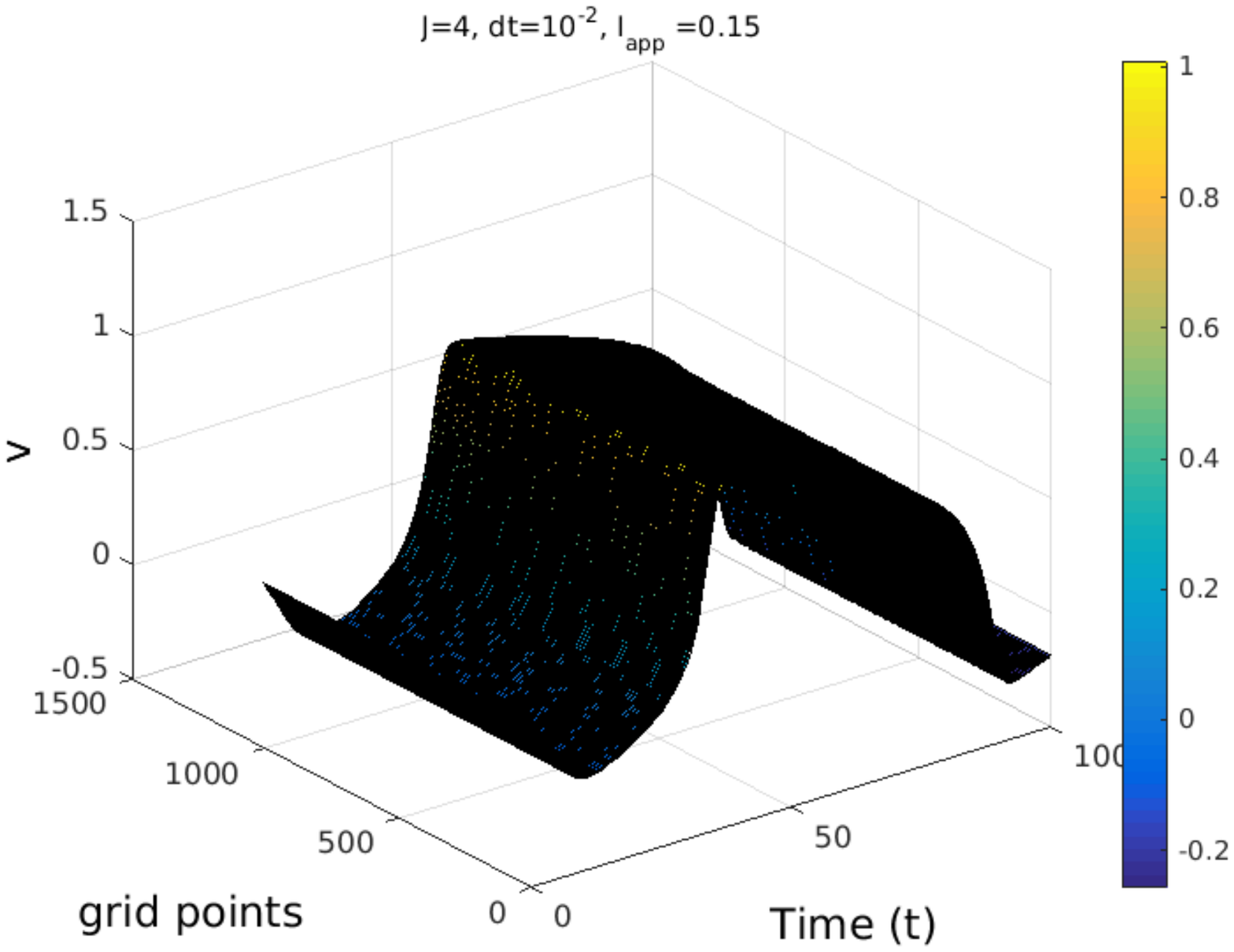}
		\vspace{-12em}
		\caption{}
		\label{surf2d_fhn}
	\end{subfigure}
	\vspace{-10em}
	\caption{\textbf{(a) Surface plot, Haar wavelet solution in 2D, (b)Haar wavelet solution in 2D}}
	\label{FHN_surf2d}
\end{figure}

Pointwise absolute error of the Haar wavelet solution $v$ is presented in Fig.\ref{error2d_fhn}.
Haar wavelet solution in surface is shown in Fig. \ref{surf2d_xy} and Fig.
\ref{surf2d_fhn} respectively. We also solved the same problem using the linear finite element in space and implicit-explicit (IE) Euler method in time. Haar wavelet solution matches with the finite element solution as shown in Fig. \ref{fe_hw_2d_compariosn}.

\begin{table}[h]
	\begin{center}
		\begin{tabular}{ |c  c  c | } 
			\hline
			& Maximum absolute error & \\
			\hline 
			$dt=10^{-2}$ & $dt=10^{-3}$ & $dt =10^{-4}$\\
			\hline
			$1.7 \times10^{-2}$ & $1.9 \times10^{-3}$ & $8.3 \times 10^{-5}$\\ 
			\hline
		\end{tabular}
		\caption{Maximum Absolute Error (for example 2) of $v$ when $J=6$ and T=0.5, k=1.}
		\label{table3}
	\end{center}
\end{table} 

\begin{figure}[h]
	\vspace{-12em}
	\centering
	\begin{subfigure}[t]{0.5\textwidth}
		\centering
		\includegraphics[width=1.2\textwidth]{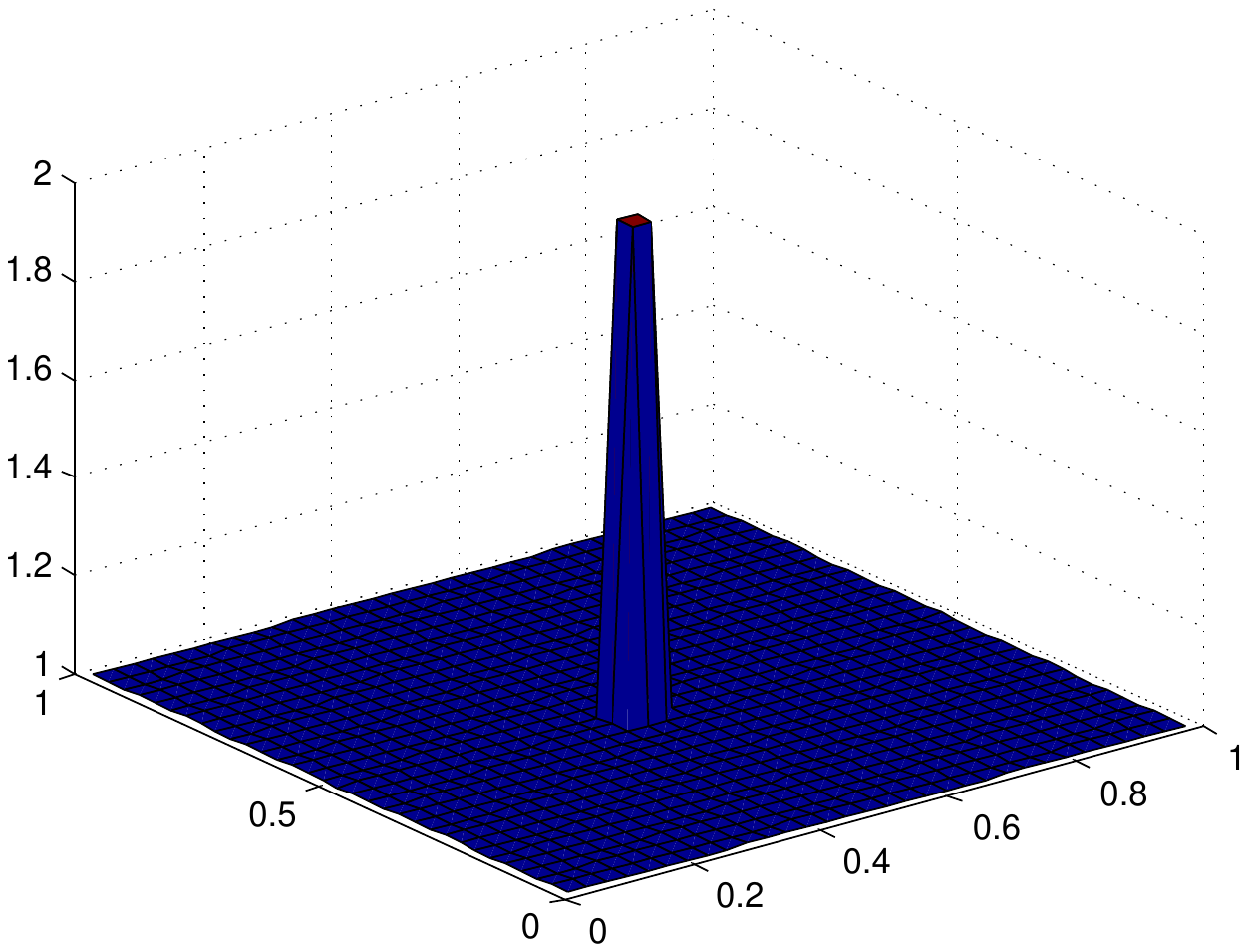}
		\vspace{-12em}
		\caption{}
		\label{fhn_k}
	\end{subfigure}\hfill
	\begin{subfigure}[t]{0.5\textwidth}
		\centering
		\includegraphics[width=1.2\textwidth]{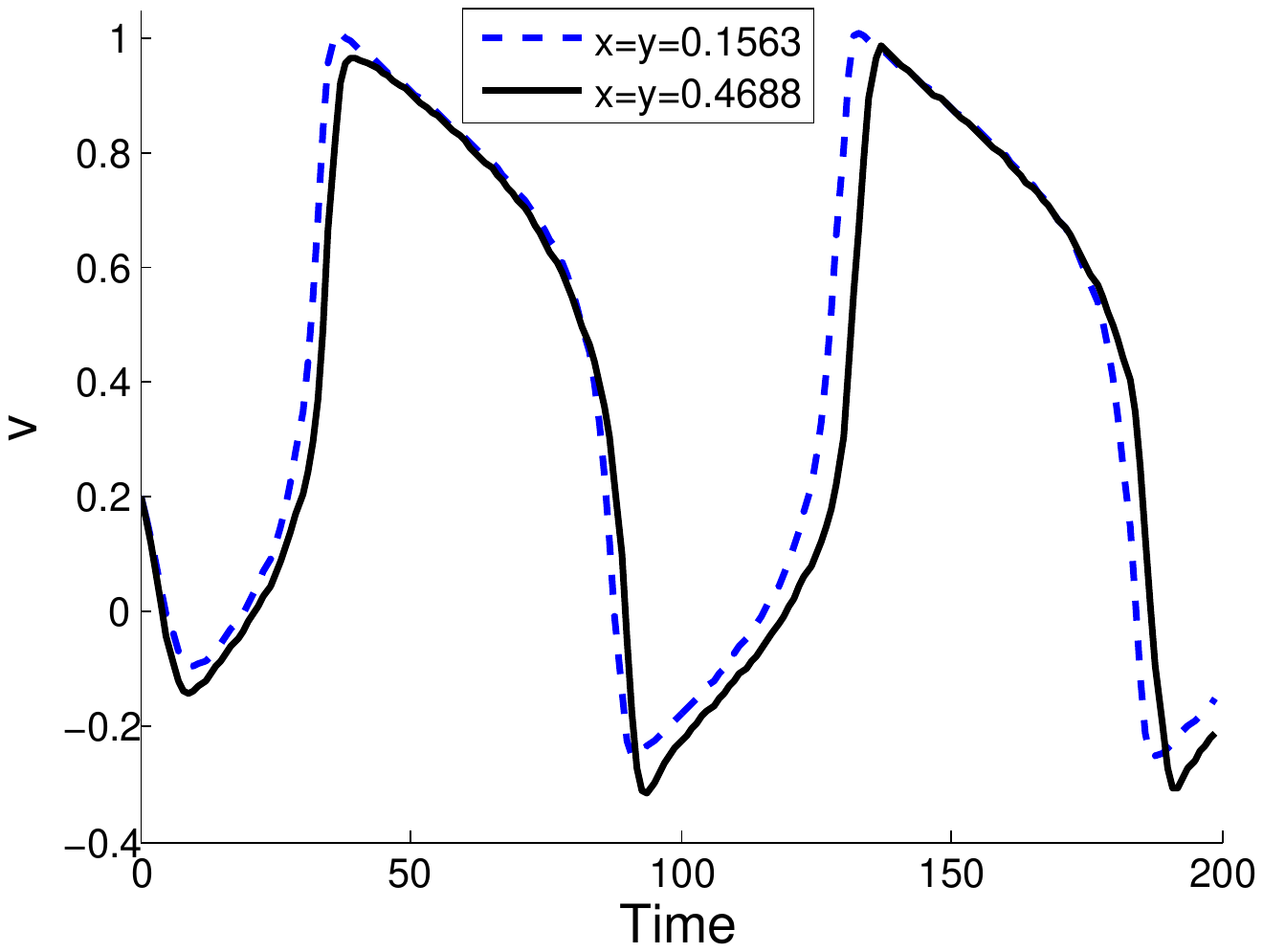}
		\vspace{-12em}
		\caption{}
		\label{AP_FHN}
	\end{subfigure}
	\vspace{-10em}
	\caption{\textbf{(a) $k$ parameter value, (b) Solution for $v$ with jump discontinuity in $k$ parameter.}}
	\label{FHN2d}
\end{figure}

Now, if the parameter $k$ in the above equation have jump discontinuity at some places in the domain, as shown in Fig. \ref{fhn_k}. This jump parameter models the local ischemia in the heart tissue. It is important to calculate the transmembrane potential with this ischemia to know the change in electrical activity of the heart.  

The solution $v$ corresponding to this value of the parameter $k$ at the points 
(0.1563,0.1563) (lies outside the jump region) and (0.4688,0.4688) (lies within the jump region of the domain) are drawn in Fig. \ref{AP_FHN}. There is no oscillation in the solution with discontinuity is noticed also the behavior of the solution clearly depict that Haar wavelets handle the jump discontinuity in the parameter $k$.

\begin{figure}[h]         
	\vspace{-12em}         
	\centering
	\begin{subfigure}[t]{0.5\textwidth}
		\centering
		\includegraphics[width=1.2\textwidth]{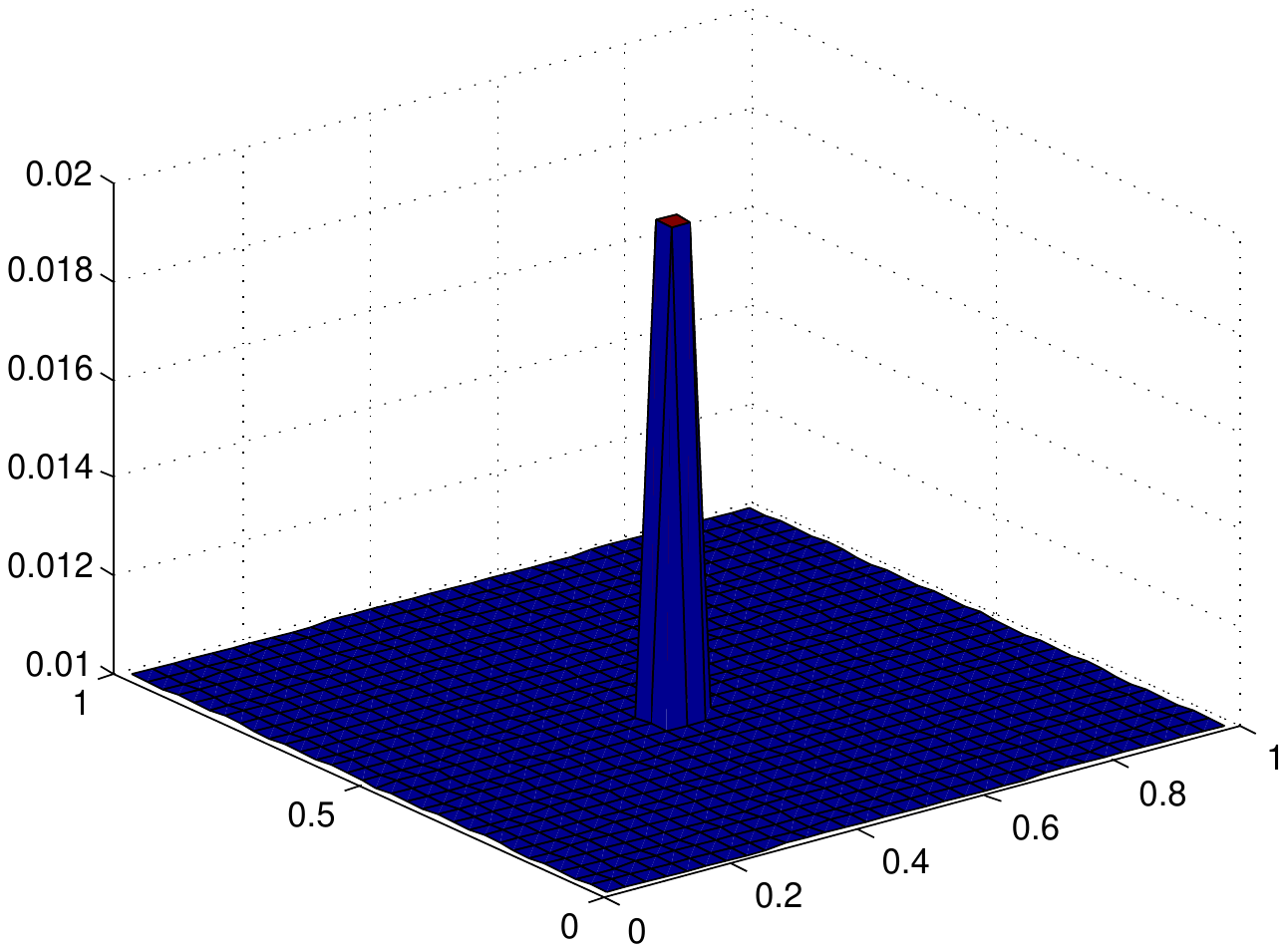}
		\vspace{-12em}
		\caption{}
		\label{fhn_eps}
	\end{subfigure}\hfill
	\begin{subfigure}[t]{0.5\textwidth}
		\centering
		\includegraphics[width=1.2\textwidth]{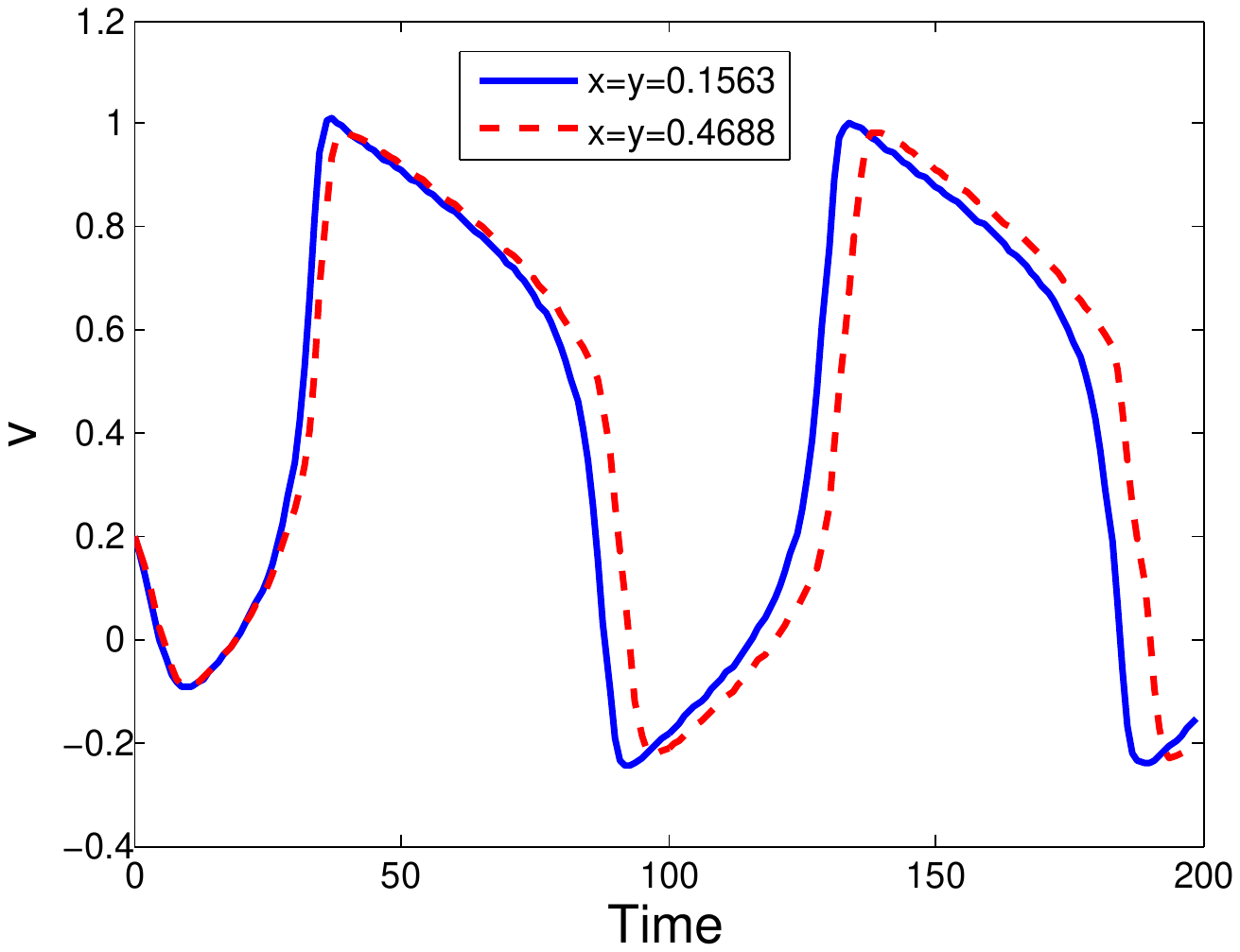}
		\vspace{-12em}
		\caption{}
		\label{AP_fhn_diff_eps}
	\end{subfigure}
	\vspace{-10em}
	\caption{\textbf{ (a) $\epsilon$ parameter value, (b) Corresponding solution $v$.}}
	\label{FHN}
\end{figure}

Again, if the discontinuity is in the parameter $\epsilon$ of the model equation, like drawn in Fig. \ref{fhn_eps}. 
The solution $v$ corresponding to this value of the parameter $\epsilon$ at the points 
(0.1563,0.1563) and (0.4688,0.4688) lies outside and inside the jump region respectively in the domain are presented
in Fig. \ref{AP_fhn_diff_eps}. Behavior of the solution with discontinuity in this parameter conclude that the $\epsilon$ discontinuity is automatically handled by the Haar wavelets.  

\paragraph{Example $3$} We will consider the example which is important in the field of cardiac electrophysiology. This model includes only inward and outward current. It contains four time constants, for the four phases of
cardiac action potential: initial, plateau, decay and recovery. We will solve the problem having multiple jump discontinuity in any of these parameter using the Haar wavelets.
\begin{align*}
&\frac{\partial v}{\partial t}- div(D\nabla v) + f(v,w) = I_{app},&  0\leq x,y \leq 1,  0\leq t \leq T\\
&\frac{\partial w}{\partial t} = G(v,w), & 0\leq x,y \leq 1,  0\leq t \leq T\\
&v(x,0)= 0.2, \hspace{5mm} w(x,0)=0.2, & 0\leq x,y \leq 1\\
&\frac{d v}{d x}(0,y,t)=0,&0\leq t \leq T\\
&\frac{d v}{d x}(1,y,t) =0, & 0\leq t \leq T\\
&\frac{d v}{d y}(x,0,t)=0,& 0\leq t \leq T\\
&\frac{d v}{d y}(x,1,t) =0,& 0\leq t \leq T
\end{align*}
where,
\begin{align*}
&f(v,w)=-\frac{w}{\tau_{in}}v^2(v-1)-\frac{v}{\tau_{out}}, \\		
&G(u,w)=  \begin{cases} 
\frac{1-w}{\tau_{open}} &  v\leq u_{gate},\\
\frac{-w}{\tau_{close}} & v>u_{gate}.\\
\end{cases}
\end{align*}

$ D = \begin{bmatrix} 1.2e{-3} & 0 \\ 0 &  2.5562e{-4}\end{bmatrix}$, and
$I_{app}=\begin{cases} 
20 &  Nt\leq 100,\\
0 &   o.w.\\
\end{cases}$.		

\begin{figure}
	\centering
	\vspace{-23em}
	\includegraphics[width=1.\textwidth]{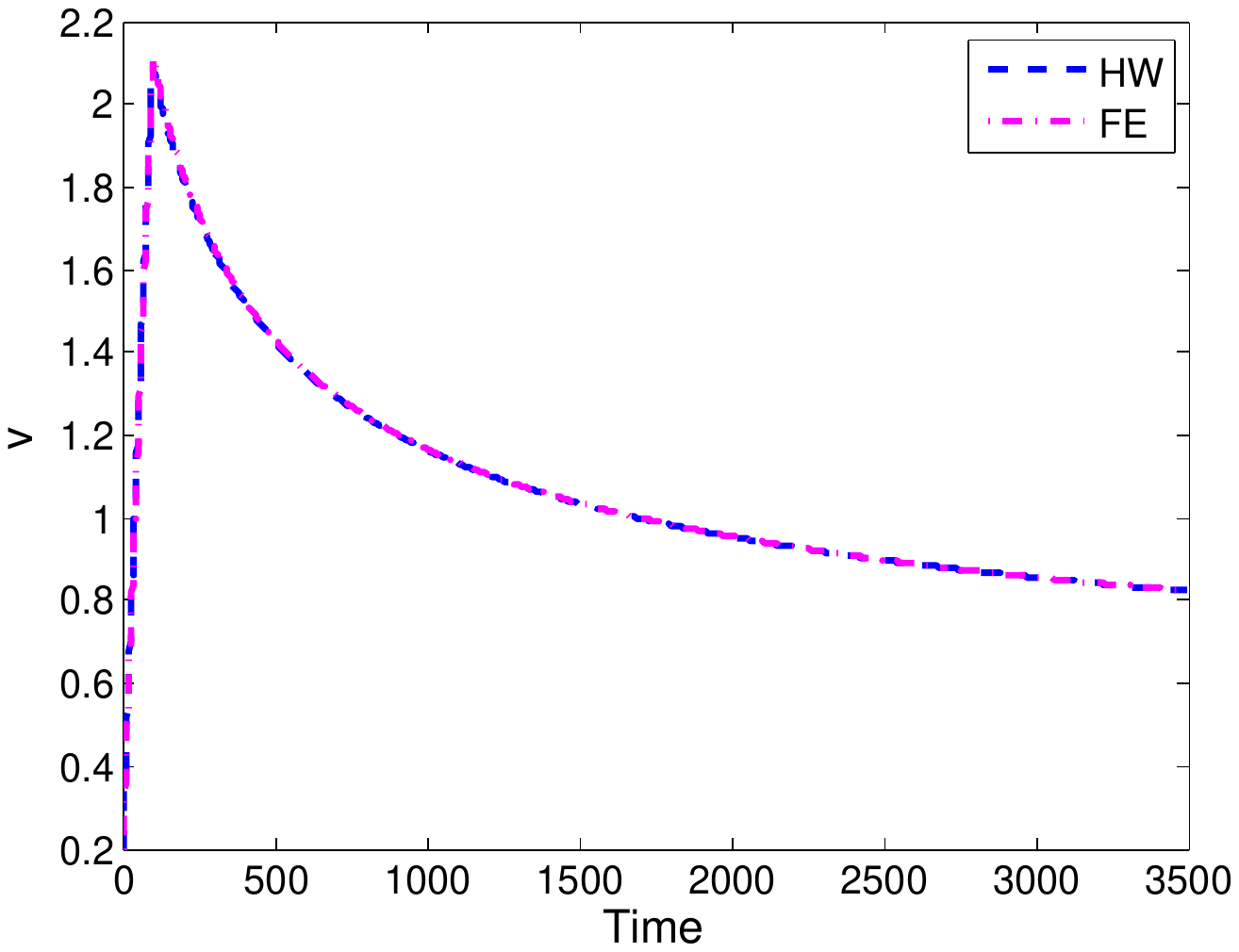}
	\vspace{-22em}
	\caption{\textbf{FE and HW solution comparison at J=4 and T=3.5 and $dt=10^{-3}$.}}
	\label{ms_fe_comparison_2d}
\end{figure}

The solution of this model with the above proposed HW matches with the FE solution, shown in Fig. \ref{ms_fe_comparison_2d}. 

\begin{figure} 
	\vspace{-14em}
	\centering
	\begin{subfigure}[t]{0.5\textwidth}
		\centering
		\includegraphics[width=1.2\textwidth]{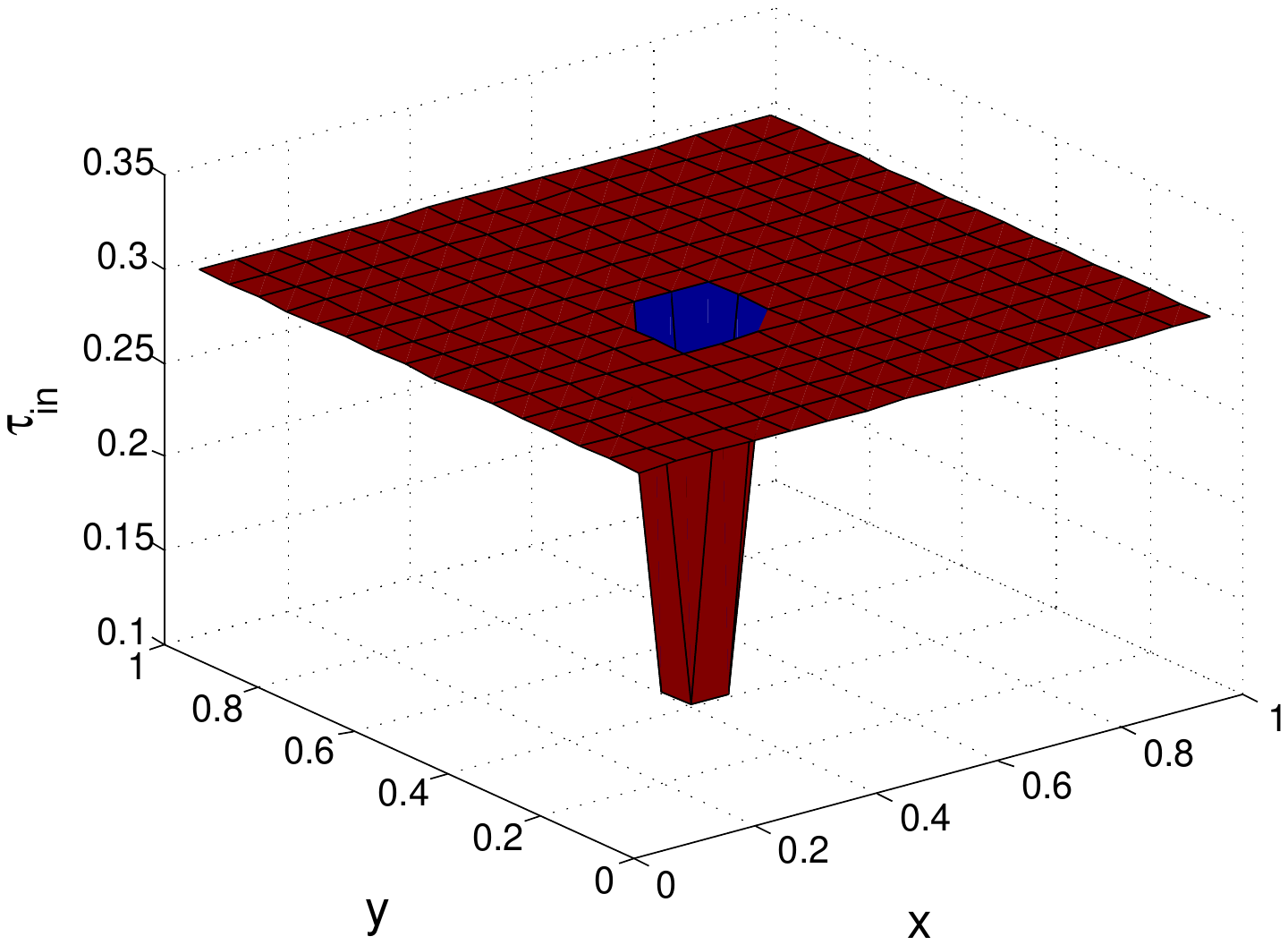}
		\vspace{-12em}
		\caption{}
		\label{ms_tin}
	\end{subfigure}\hfill
	\begin{subfigure}[t]{0.5\textwidth}
		\centering
		\includegraphics[width=1.2\textwidth]{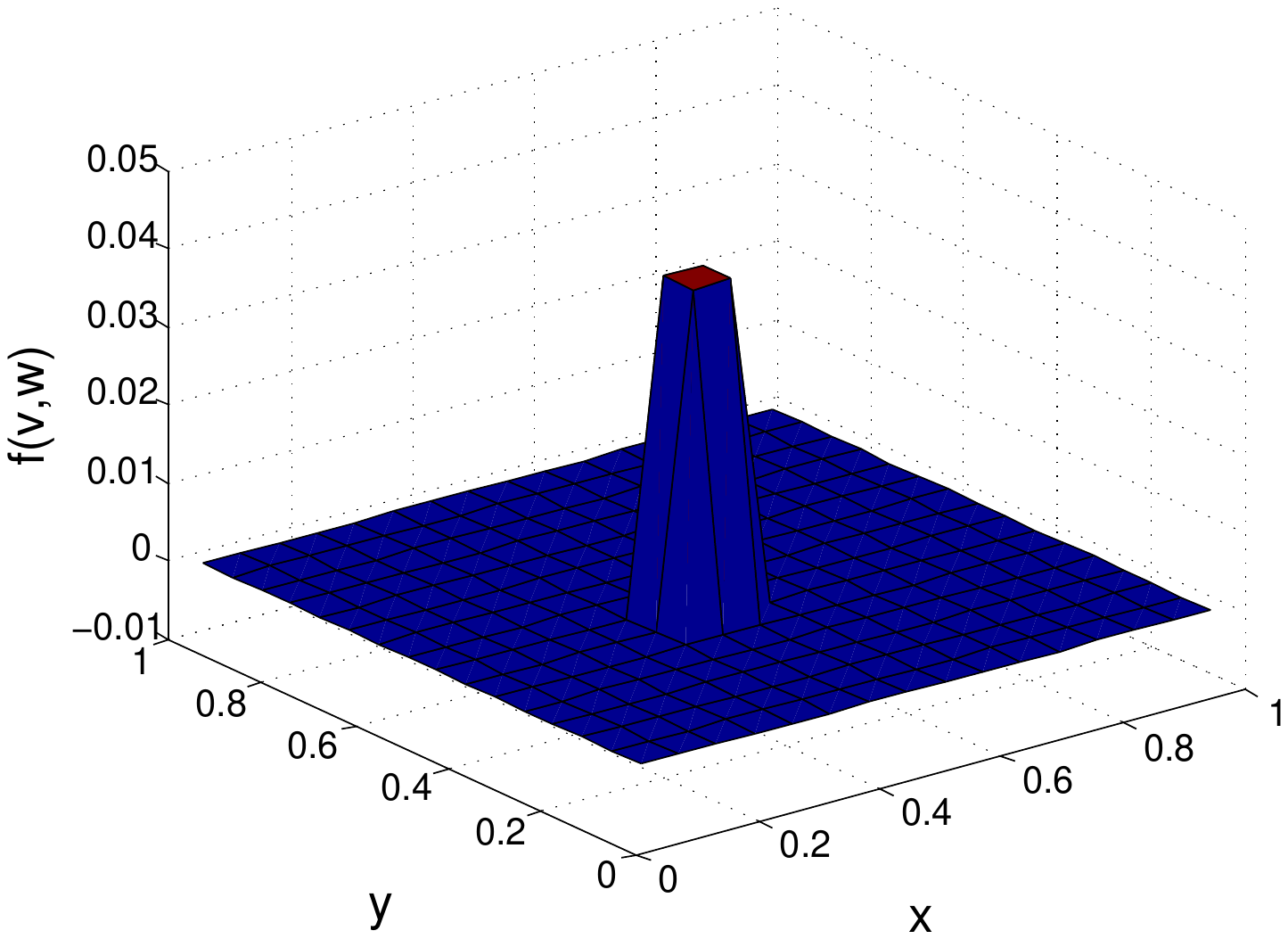}
		\vspace{-12em}
		\caption{}
		\label{ms_Ion}
	\end{subfigure}
	\vspace{-10em}
	\caption{\textbf{(a)$\tau_{in}$ parameter value, (b) Function $f(v,w)$.}}
	\label{MS}
\end{figure}

\begin{figure}
	\vspace{-17em}
	\centering
	\includegraphics[width=0.8\textwidth]{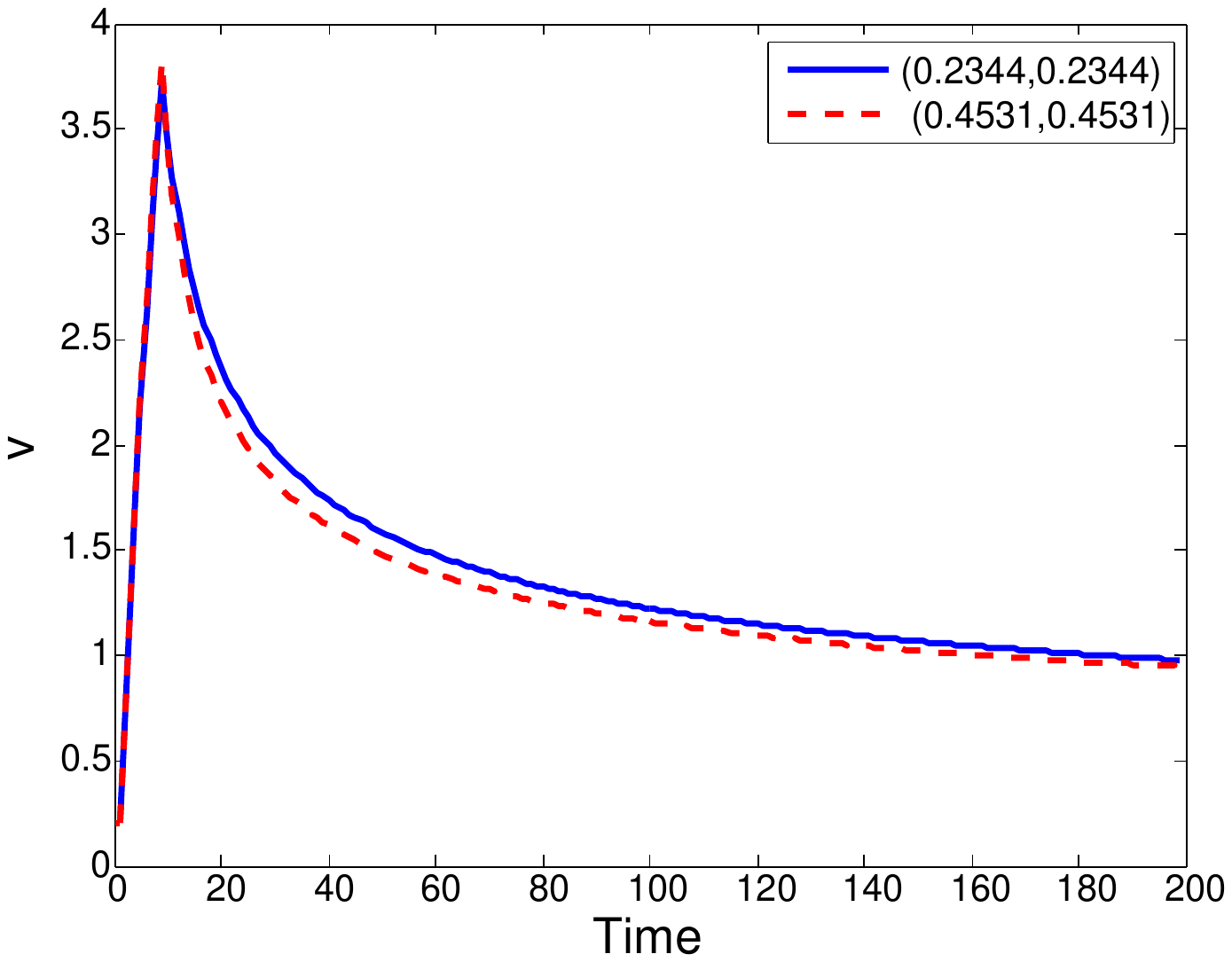}
	\centering
	\vspace{-16em}
	\caption{Solution $v$ with jump discontinuity in $\tau_{in}$ parameter.}
	\label{AP_ms}
\end{figure}

\begin{figure}
	\vspace{-13em}
	\centering
	\begin{subfigure}[t]{0.5\textwidth}
		\centering
		\includegraphics[width=1\textwidth]{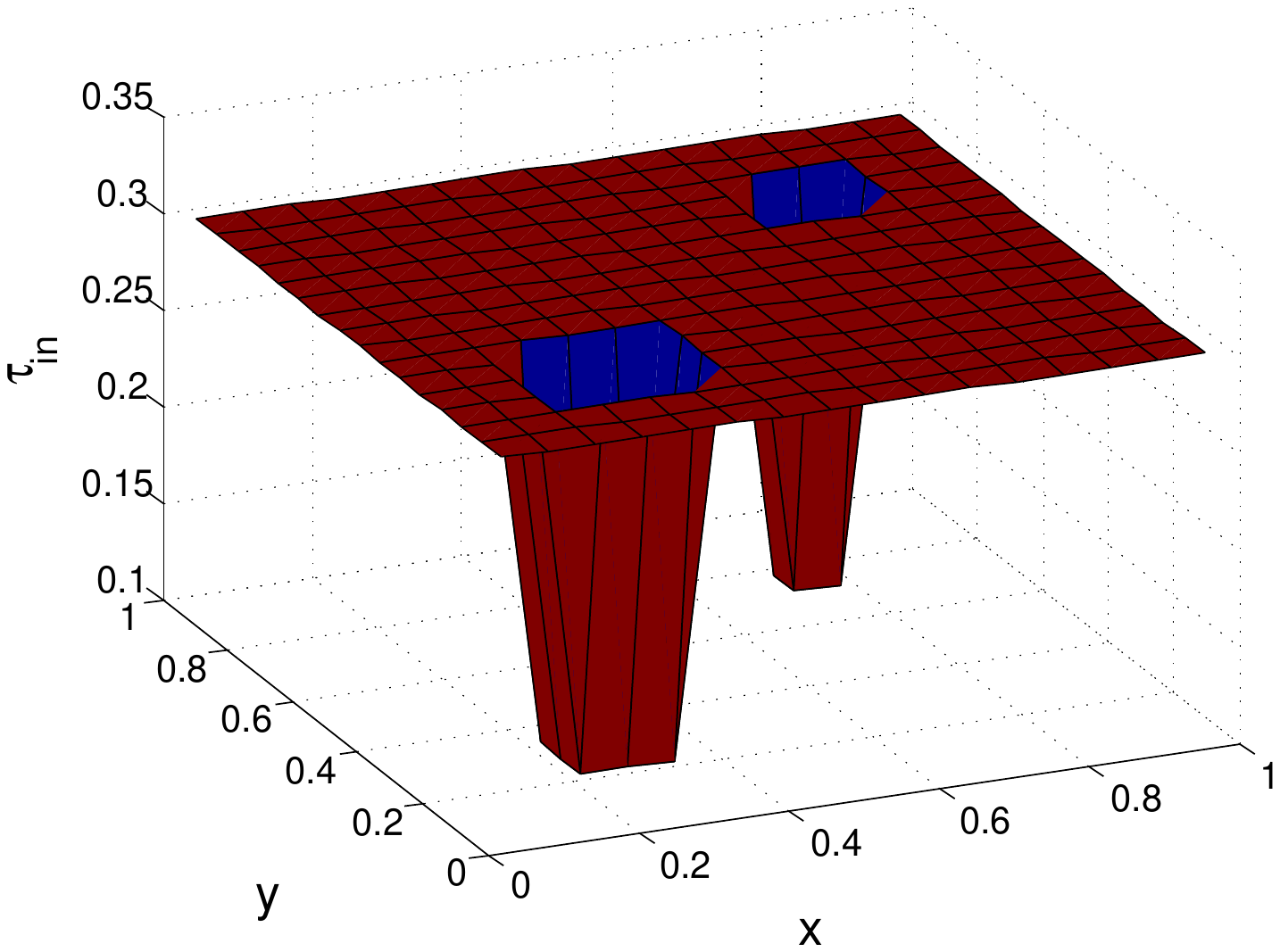}
		\vspace{-12em}
		\caption{}
		\label{ms_tin_2jump}
	\end{subfigure}\hfill
	\begin{subfigure}[t]{0.5\textwidth}
		\centering
		\includegraphics[width=1\textwidth]{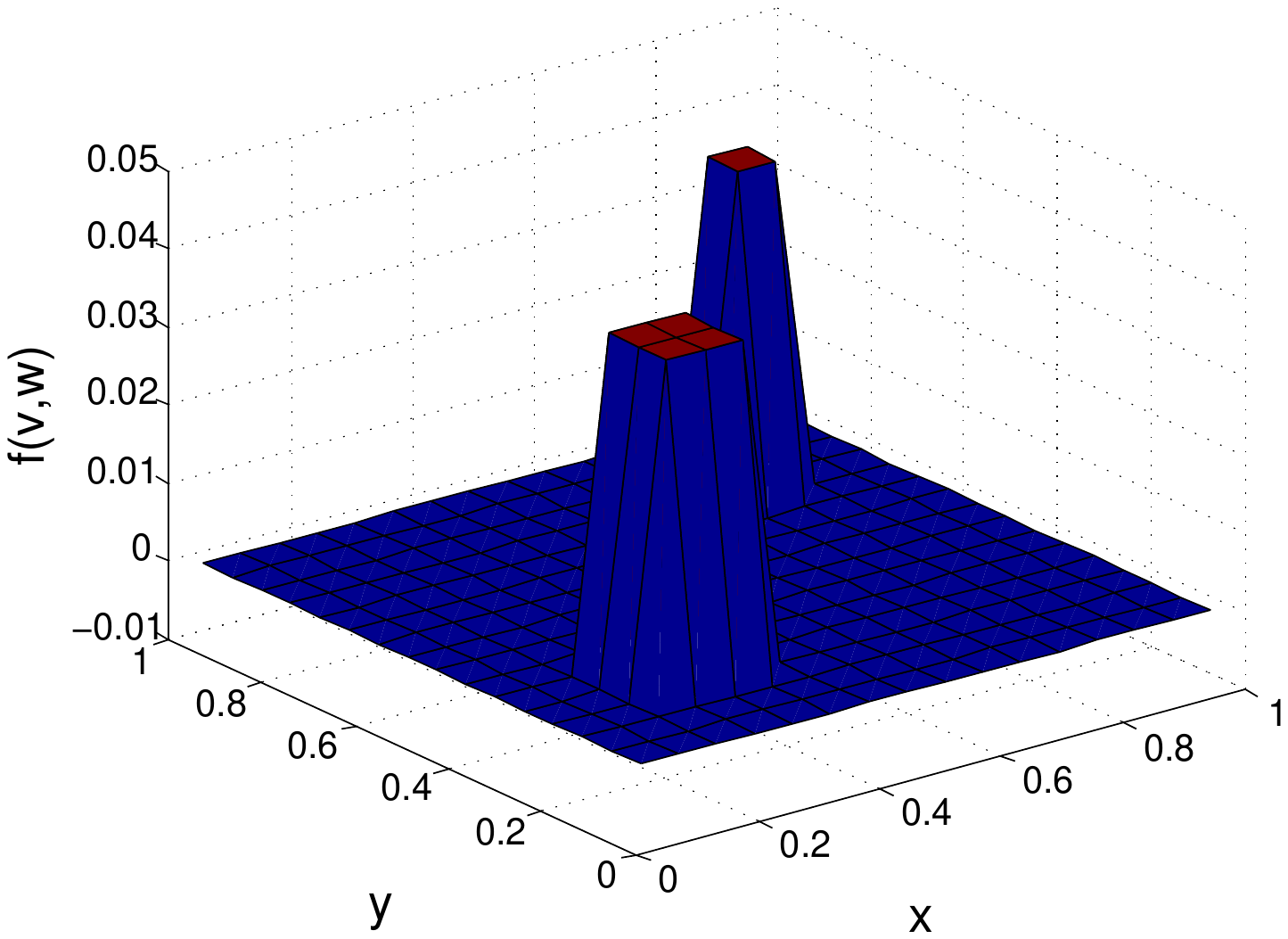}
		\vspace{-12em}
		\caption{}
		\label{ms_Ion_2jump}
	\end{subfigure}
	\vspace{-9em}
	\caption{\textbf{(a) $\tau_{in}$ parameter value with jump discontinuity at two regions, (b) Corresponding function $f(v,w)$.}}
	\label{MS_2jump}
\end{figure}

\begin{figure}   
	\vspace{-15em}
	\centering
	\includegraphics[width=0.7\textwidth]{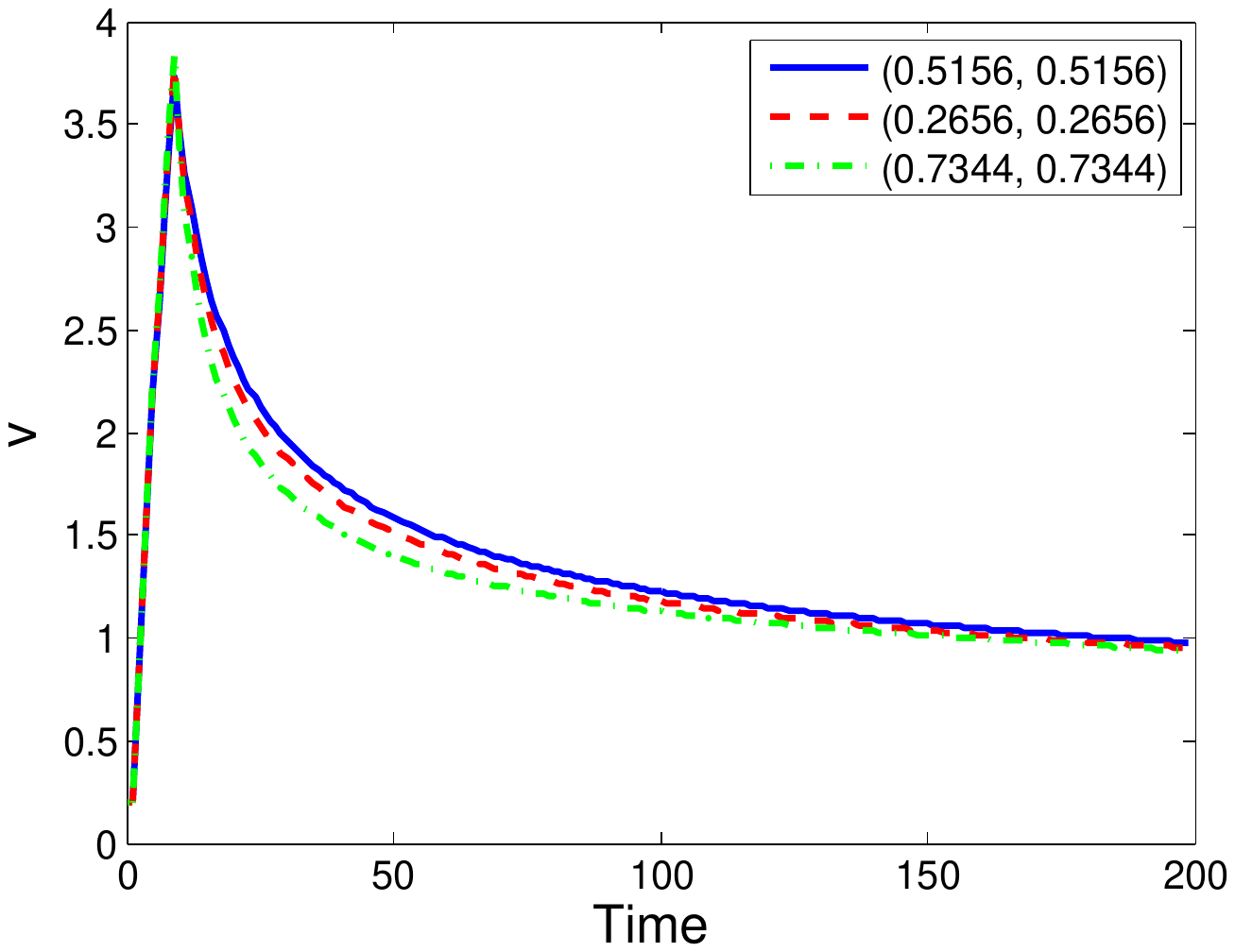}
	\centering
	\vspace{-14em}
	\caption{Solution $v$ with jump discontinuity at two places in $\tau_{in}$ parameter.}
	\label{AP_ms_2jump}
\end{figure}

Now, if the parameter $\tau_{in}$ in the above equation have jump discontinuity at some places in the domain, as shown in 
Fig. \ref{ms_tin} and Fig. \ref{ms_tin_2jump}, the corresponding function $I_{ion}(v,w)$ will also be discontinuous, shown in Fig. 
\ref{ms_Ion} and \ref{ms_Ion_2jump}. This discontinuity is automatically handled using Haar wavelets. The solution $v$ corresponding
to $\tau_{in}$ in Fig. \ref{ms_tin} at the points (0.2344,0.2344) and (0.4531,0.4531), lie outside and inside the jump region of
the domain are shown in Fig. \ref{AP_ms}. Now, 
the solution $v$ corresponding to $\tau_{in}$ in Fig. \ref{ms_tin_2jump} at the points (0.2656, 0.2656), (0.5156, 0.5156) and 
(0.7344, 0.7344) are presented in Fig. \ref{AP_ms_2jump}. From these Figures, we can clearly see that Haar wavelets automatically handle
the multiple jumps in the parameter $\tau_{in}$. This clearly describes the power of using Haar wavelets.

This type of single or multiple jump discontinuity problems are useful to calculate the effect of ischemia in the cardiac electrical activity in a tissue with single or multiple ischemic subregions. In such type of problems local ischemia is handled by changing the value of the corresponding parameters into the model at these ischemic subregions only. Thus, Haar wavelet method is a computationally effective tool to solve such types of pathological problems.

Comparison of the computation time of finite element solution and the Haar wavelet solution for 2D model in example 2 and example 
3 is presented in Table \ref{cputime2d}. Clearly, from the table, the haar wavelet method reduces the computational cost in 
comparison of FEM.\\

\begin{table}
	\vspace{-1em}
	\begin{center}
		\begin{tabular}{ | c |c | c | c | c| } 
			\hline
			& J & dt & FEM CPU time (sec) & HW CPU time (sec)\\
			\hline 
			Example 2 & 3 (256 grid) & $10^{-2}$ & 2.06464 & 1.690459 \\
			
			& 4 (1024 grid) & $10^{-2}$ & 7.30766 & 6.40874 \\
			%
			\hline	
			Example 3 & 3 (256 grid) & $10^{-2}$ & 2.94036 & 1.527416\\
			
			& 4 (1024 grid) & $10^{-2}$ & 9.84829 & 5.342718 \\
			\hline 
			
		\end{tabular}
		\vspace{0em}
		\caption{CPU time 2D}
		\label{cputime2d}
	\end{center}
\end{table}

\paragraph{Example 4} 
We will solve the following Monodomain cardiac tissue level model with the complex non-linear Hodgkin-Huxely ionic model \cite{HH} using 
the above developed Haar wavelet method.
\begin{align*}
\frac{\partial v}{\partial t}- div(D(x)\nabla v) +I_{ion}(v,w)   &= I_{app}, &  x \in \Omega,  0\leq t \leq T\\
\frac{\partial w}{\partial t}-g(v,w)&= 0, & x \in \Omega,  0\leq t \leq T\\
v(x,0)= v_0(x,0), \hspace{5mm} w(x,0)&=w_{i,0}(x,0), & x \in \Omega\\
n^T D(x) \nabla v &=0, & x \in \partial \Omega, 0\leq t \leq T
\end{align*}
where, $I_{app}$ is the applied stimulus and $D(x)$ is the conductivity tensor defined as: 
\begin{align*}
D(x)= {\sigma_t} I +({\sigma_l} -{\sigma_t} )b_l(x) {b_l}^t(x) \hspace{3mm} with \hspace{3mm}
\sigma_{s}=\frac{\widetilde{l}
	\sigma^i_{s}}{1+\widetilde{l}}, s=l,t,n.
\end{align*}	
$g(v,w)$ an $I_{ion}(v,w)$ are given by the models at the cell level, called ionic models. Here, we are choosing the Hodgkin
Huxely ionic model.

Hodking and Huxely in 1952 modeled an ionic model to calculate the action potential in a squid giant axon. In this
model, sodium and potassium ionic channels and hence corresponding two separate ionic currents are taken into account.
All other small ionic currents are lumped into the leakage current. The model equation is given as follows
\begin{align}
\label{Hudgkin}
I_{ion} = g_{Na}m^3h(v-E_{Na}) + g_{K}n^4(v-E_{K})+ g_{L}(v-E_{L}),
\end{align}
where, $g_{Na}$ and $g_{K}$ are the maximal conductances corresponding to their membrane which depends on the voltage
and time, $E_{Na}, E_{K}, E_{L}$ are the corresponding ion resting potentials. $I_{app}$ is the externally applied current.
$m$, $h$ and $n$ are the gating variables and the dynamics of each of the gating variables is described by the ordinary
differential equation:
\begin{align}
\label{gating}
\frac{dw}{dt}= \alpha_{w}(1-w)-\beta_{w}, \hspace{3mm} w=m,h,n,
\end{align}
where, $\alpha_{w}$ and $\beta_{w}$ are the rates depends upon the membrane voltage.

\begin{figure}
	\vspace{-15em}
	\centering
	\includegraphics[width=0.7\textwidth]{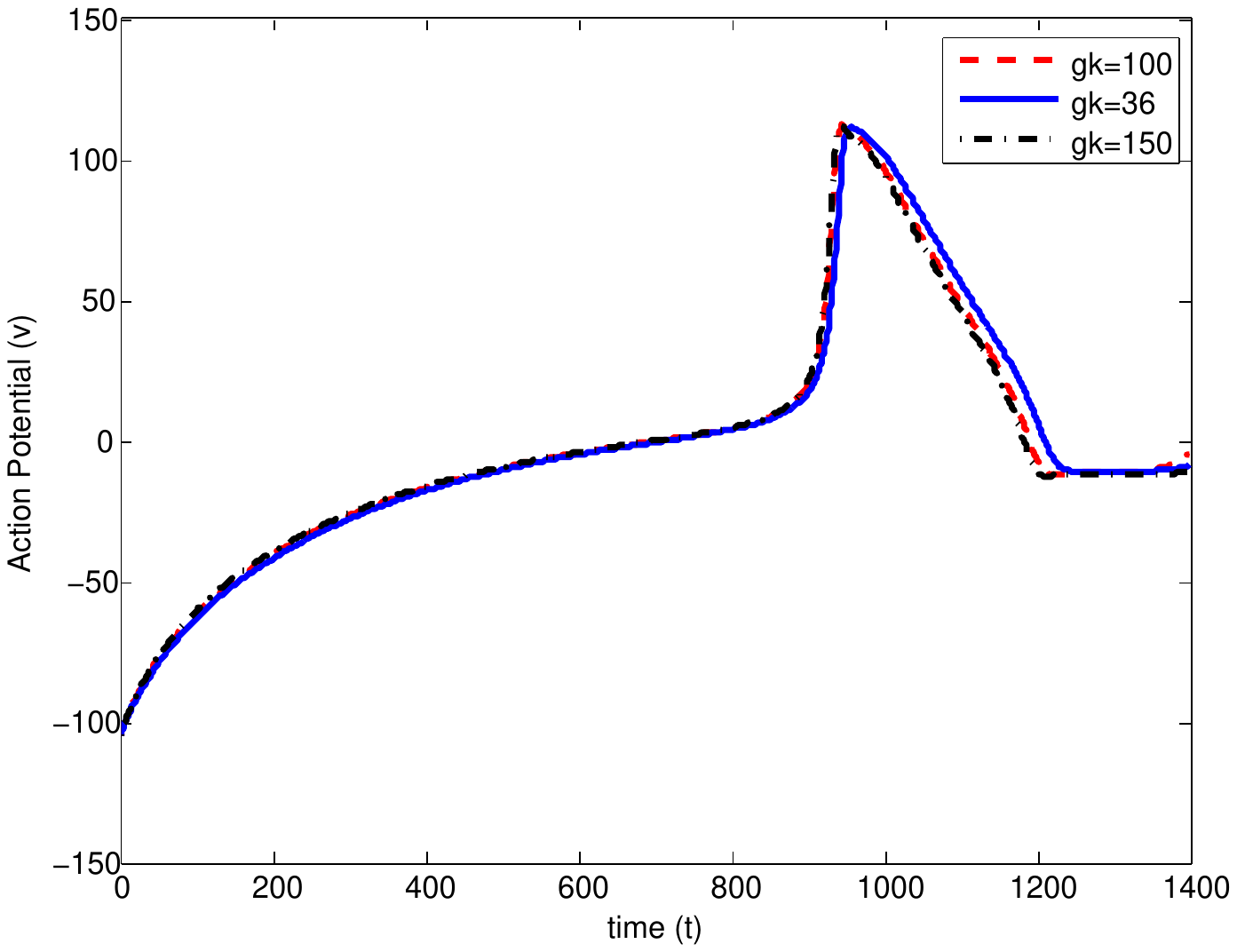}
	\centering
	\vspace{-14em}
	\caption{Action Potential for HH model}
	\label{AP_HH}
\end{figure}

\begin{figure}
	\vspace{-13em} 
	\centering
	\begin{subfigure}[t]{0.5\textwidth}
		\centering
		\includegraphics[width=1.\textwidth]{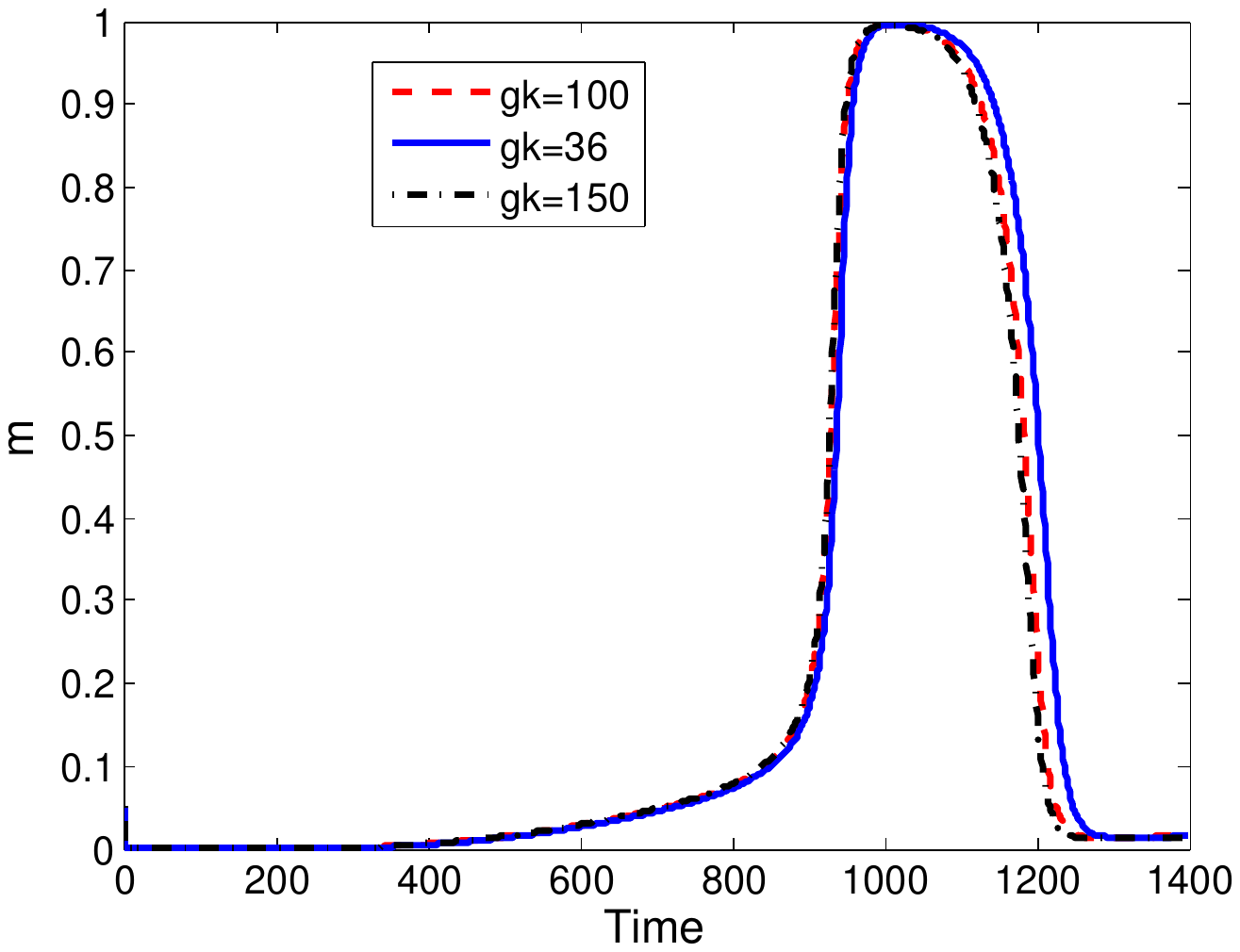}
		\vspace{-12em}
		\caption{}
		\label{HH_m}
	\end{subfigure}\hfill
	\begin{subfigure}[t]{0.5\textwidth}
		\centering
		\includegraphics[width=1.\textwidth]{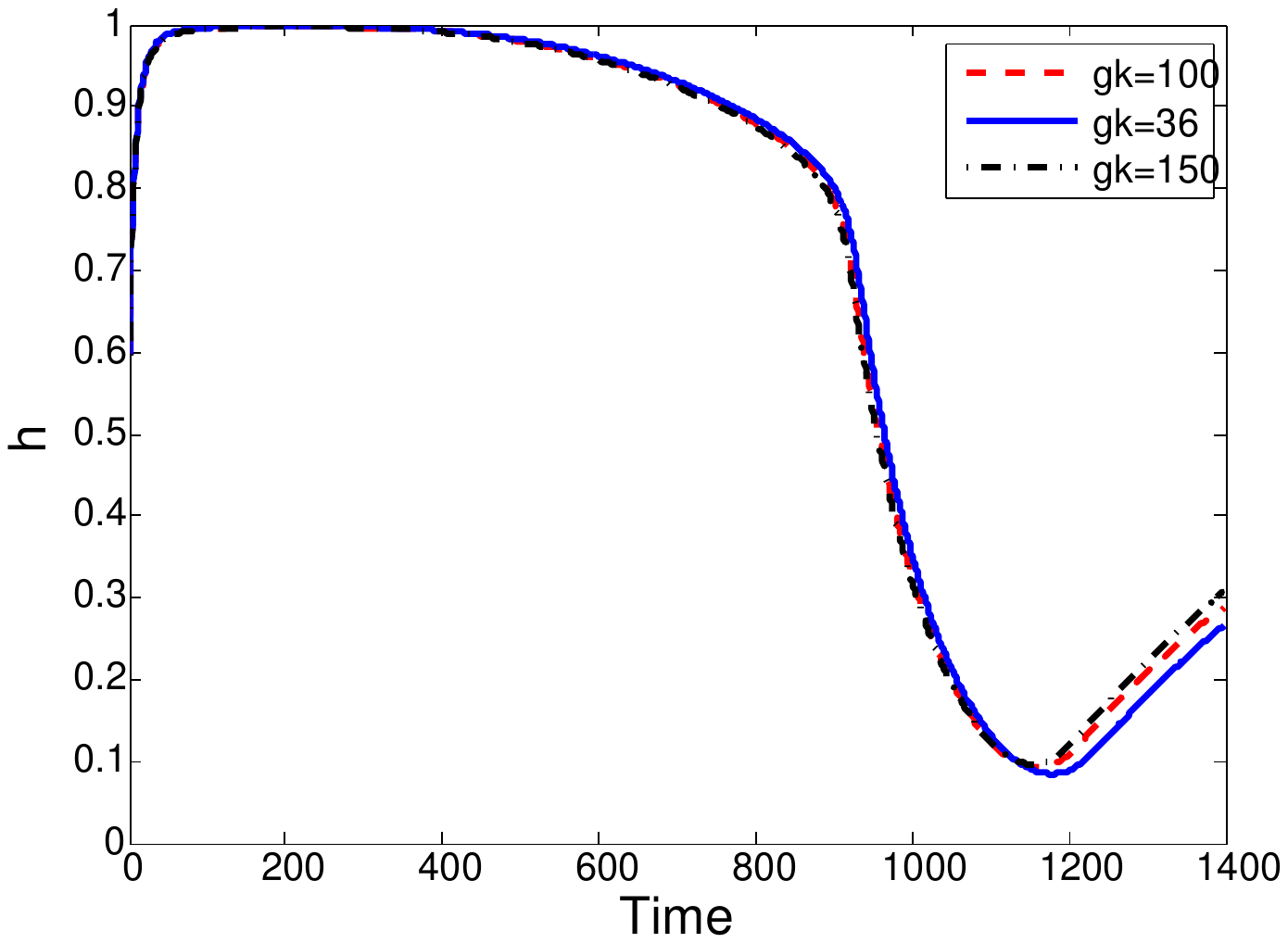}
		\vspace{-12em} 
		\caption{}
		\label{HH_h}
	\end{subfigure}
	\vspace{-10em}
	\caption{\textbf{(a)m gate, (b) h gate.}}
	\label{HH}
\end{figure}

Here, the conductivity of potassium ion $g_{K}$ depends on the voltage and time. If we change this parameter value membrane potential
or the action potential will also change in the tissue. Action potential and the gating variables for the different value of the parameter is presented in Fig. \ref{AP_HH}, \ref{HH_m} and \ref{HH_h}.

The jump discontinuity in the parameter $g_K$ and hence in $I_{ion}$ is handled by the Haar wavelets.

\begin{figure}
	\vspace{-17em}
	\centering
	\includegraphics[width=0.7\textwidth]{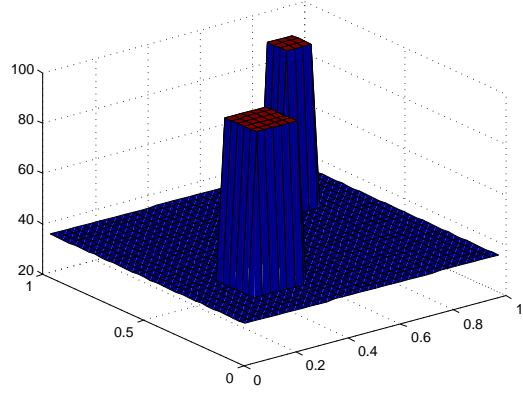}
	\centering
	\vspace{-14em}
	\caption{$g_K$ parameter value having jump discontinuity at two places,}
	\label{HH_gk_2jump}
\end{figure}

\begin{figure}
	\vspace{-17em}
	\centering
	\includegraphics[width=0.8\textwidth]{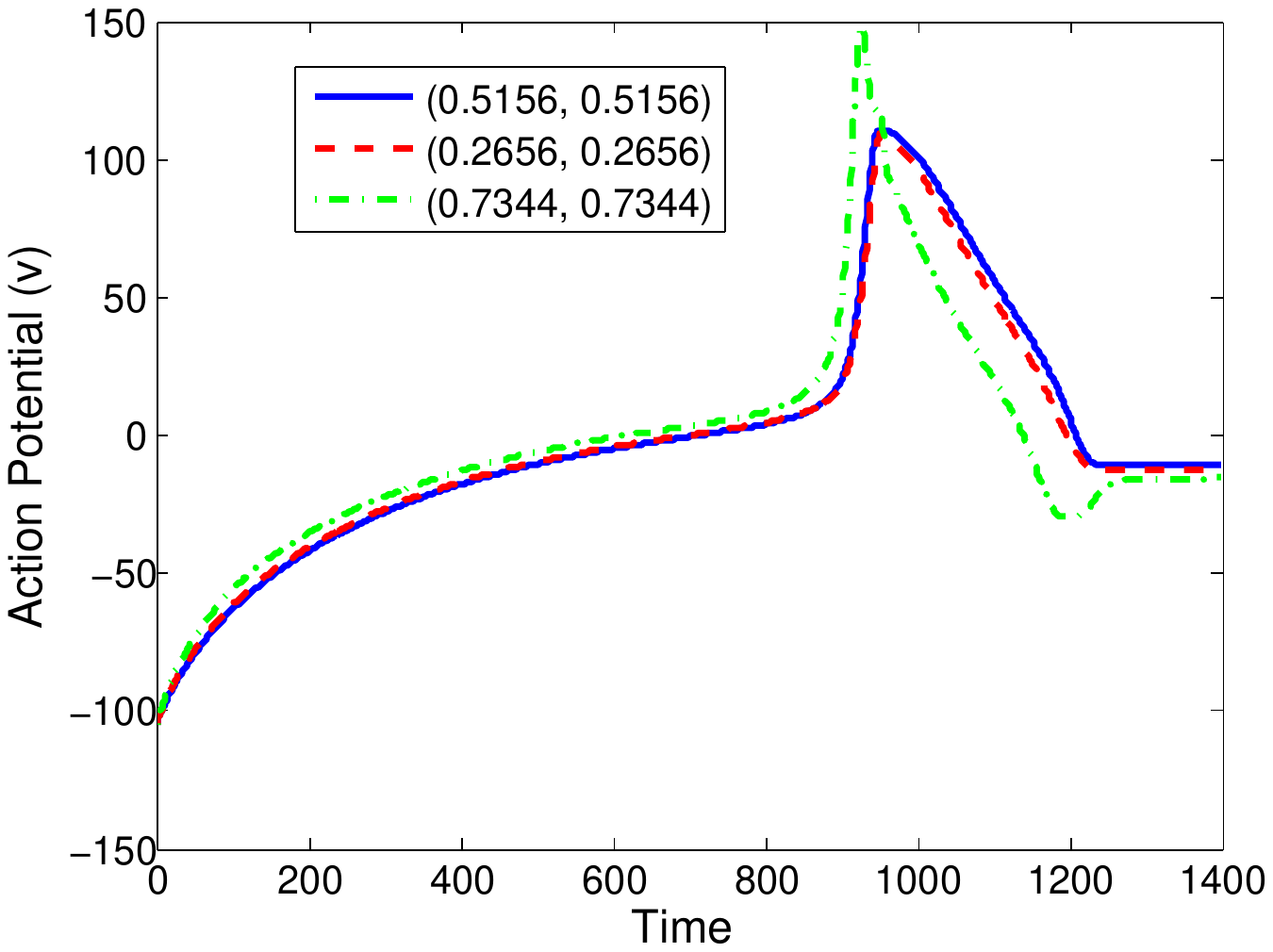}
	\centering
	\vspace{-17em}
	\caption{Action Potential for HH model $g_K$ parameter value having jump discontinuity at two places.}
	\label{AP_HH_2jump}
\end{figure}

Now, we are choosing the parameter $g_k$ value having jump discontinuity at the two subdomains of the domain, as shown 
in Fig. \ref{HH_gk_2jump}. The corresponding function $I_{ion}(v,w)$ will also be discontinuous in the corresponding subdomains. The solution $v$ corresponding to this $g_k$ at the points (0.2156, 0.2156), (0.5156, 0.5156) and 
(0.7344, 0.7344) is presented in Fig. \ref{AP_HH_2jump}.

\newpage

\paragraph{Example $5$} We will consider the three dimensional example which is important in the field of cardiac electrophysiology. This model includes only inward and outward current. It contains four time constants, for the four phases of cardiac action potential: initial, plateau, decay and recovery. We will solve the problem having jump discontinuity in any of these parameter using the Haar wavelets.
\begin{align*}
\frac{\partial v}{\partial t} - \bigg(d_{11} \frac{\partial ^2 v}{\partial x^2}+d_{22} 
\frac{\partial ^2 v}{\partial y^2}+d_{33} \frac{\partial ^2 v}{\partial z^2}\bigg) 
+I_{ion} &= I_{app}, & 0\leq x,y,z \leq 1,  0\leq t \leq T\\
\frac{\partial w}{\partial t}&=v-2w, & 0\leq x,y,z \leq 1,  0\leq t \leq T\\
v(x,y,z,0)= 0.2, \hspace{5mm} w(x,y,z,0)&=0.2, & 0\leq x,y,z \leq 1
\end{align*}
with Neumann boundary conditions in $v$. 
where,
\begin{align*}
&I^{ion}=-\frac{w}{\tau_{in}}u^2(u-1)-\frac{u}{\tau_{out}} \\		
&G(u,w)=  \begin{cases} 
\frac{1-w}{\tau_{open}} &  u\leq u_{gate},\\
\frac{-w}{\tau_{close}} & u>u_{gate}.\\
\end{cases}
\end{align*}

\begin{figure}[h]
	\centering
	\vspace{-11em}
	\includegraphics[width=0.6\textwidth]{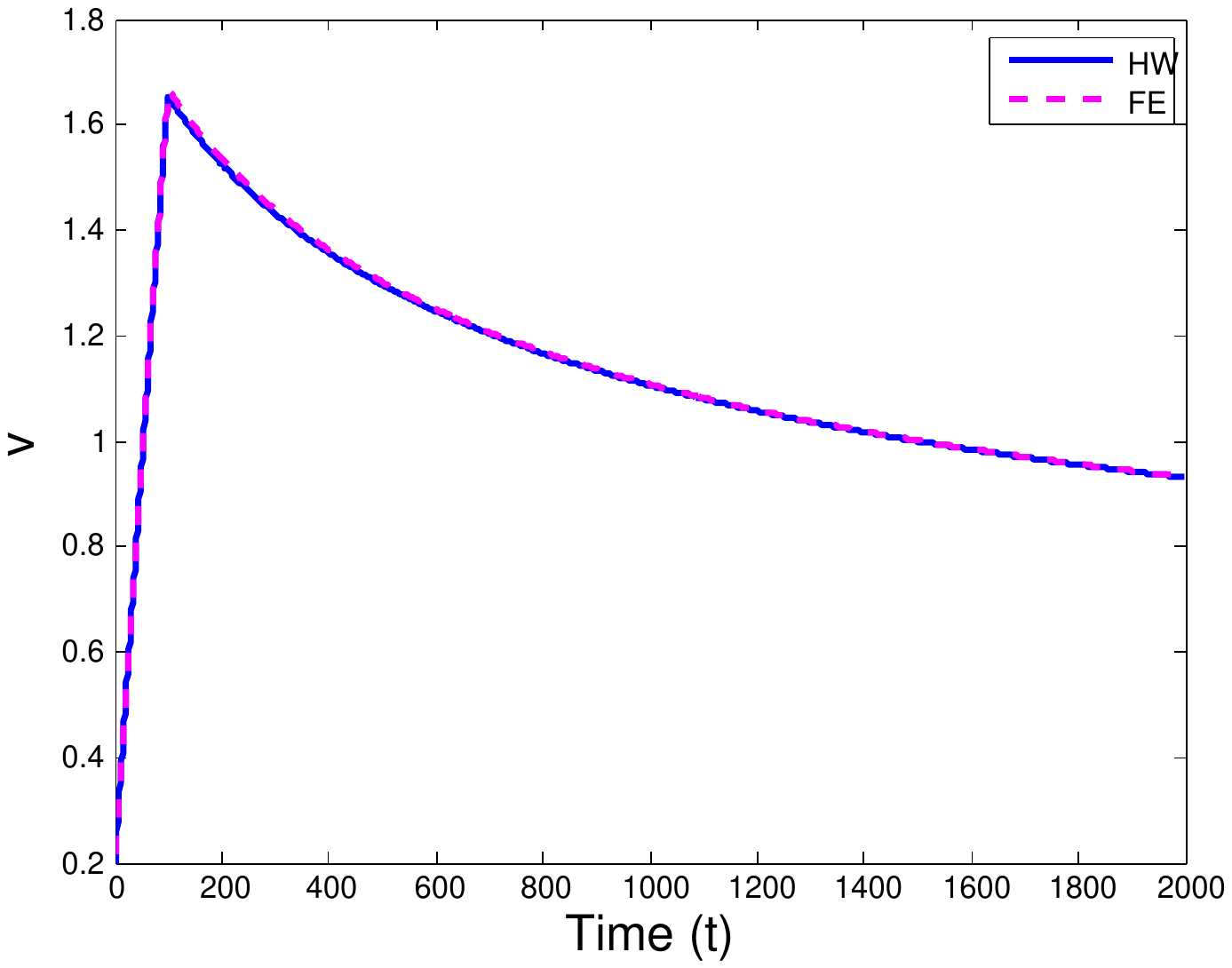}
	\vspace{-12em}
	\caption{\textbf{Comparison FE and HW solution $v$ with J=4 and T=1 and $dt=10^{-3}$.} at point (0.5312, 0.5312, 0.5312).}
	\label{ms_hw_3D_comparison}
\end{figure}

\begin{figure}[h]
	\vspace{-12em}
	\centering
	\begin{subfigure}[t]{0.5\textwidth}
		\centering
		\includegraphics[width=1.2\textwidth]{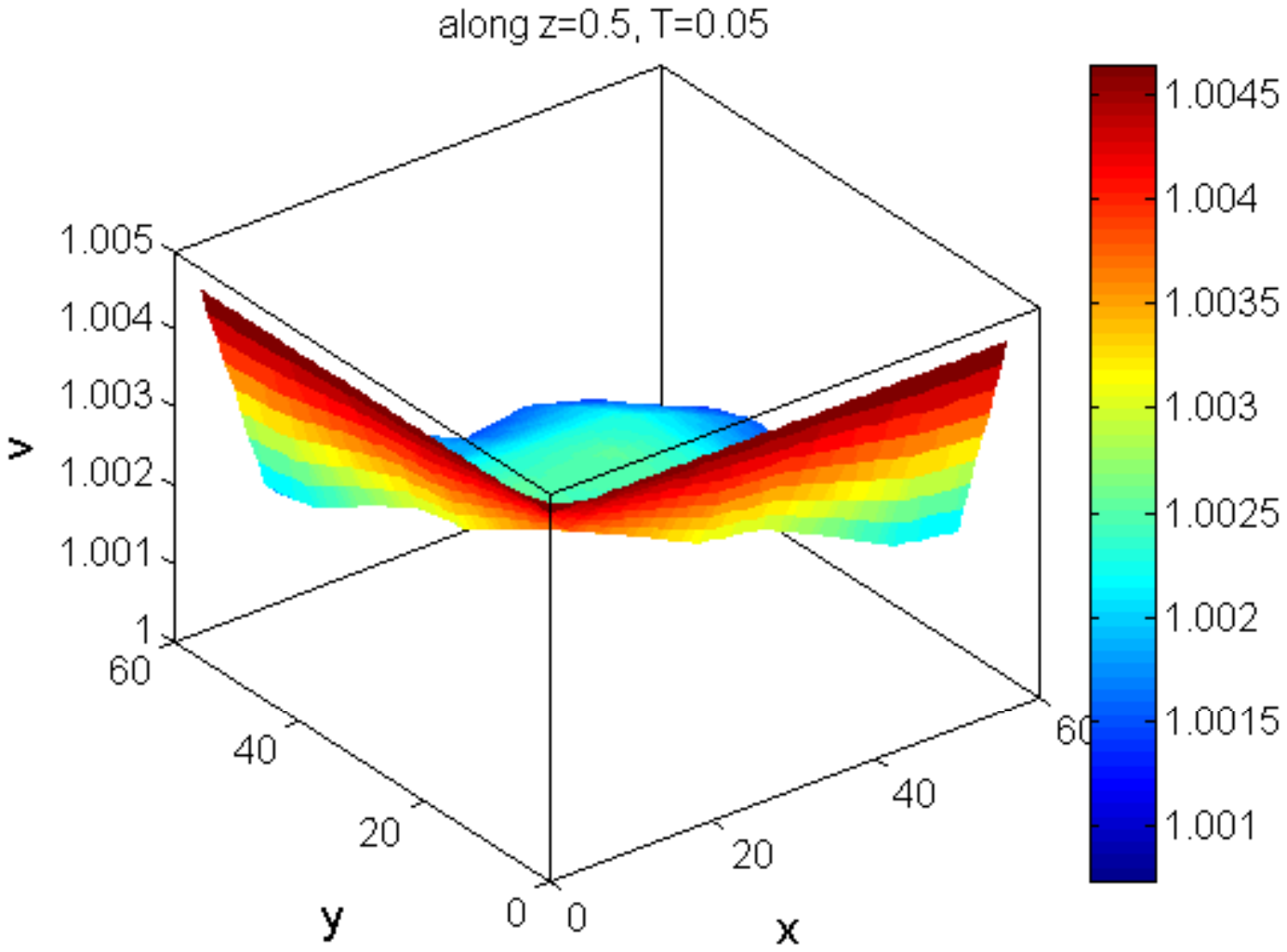}
		\vspace{-13em}
		\caption{}
		\label{ms3d_surf_xy_T50}
	\end{subfigure}\hfill
	\begin{subfigure}[t]{0.5\textwidth}
		\centering
		\includegraphics[width=1.2\textwidth]{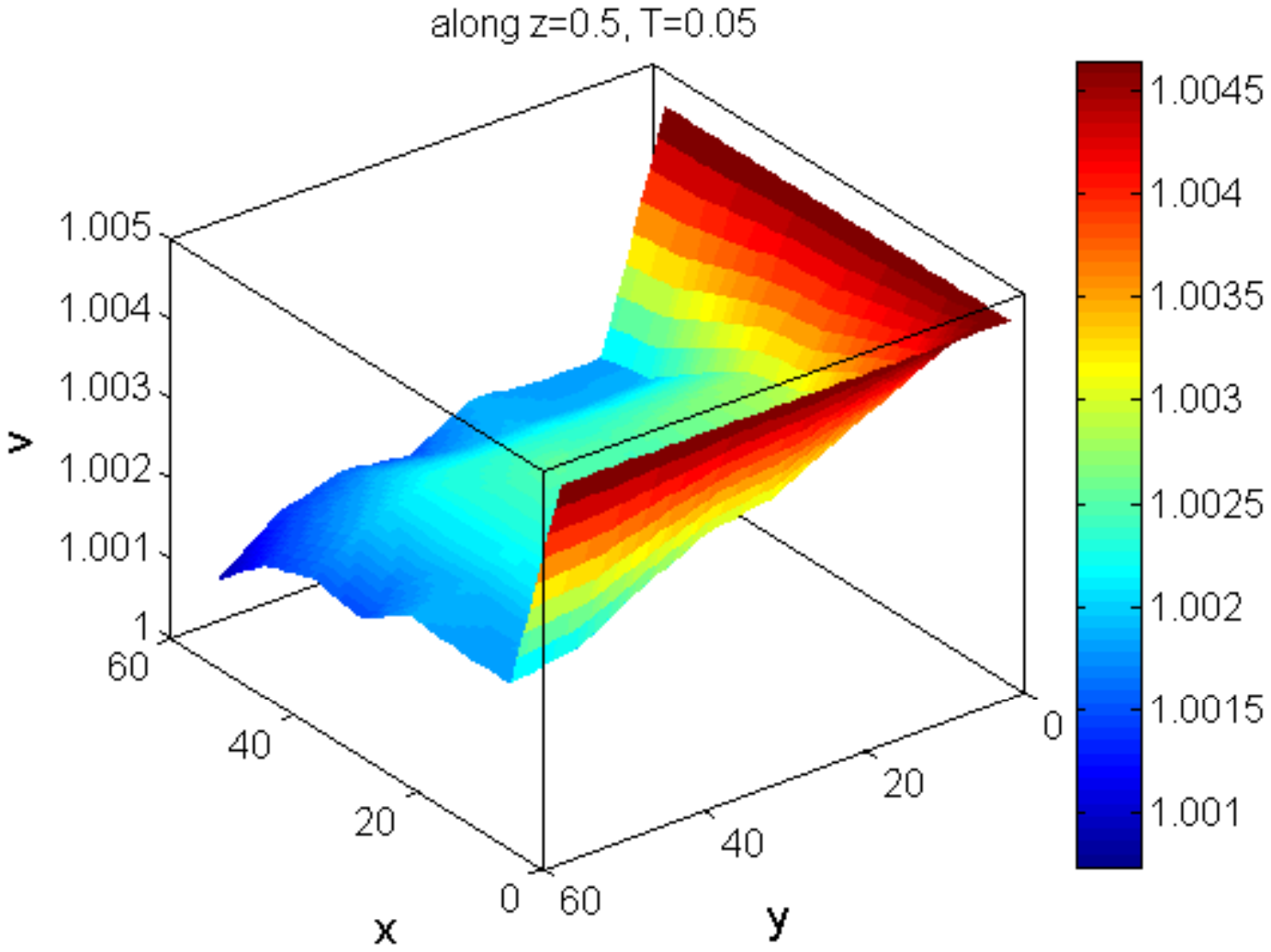}
		\vspace{-13em}
		\caption{}
		\label{ms3d_surf_xy_T50_rotate}
	\end{subfigure}\\
	\vspace{-16em}
	\begin{subfigure}[t]{0.5\textwidth}
		\centering
		\includegraphics[width=1.2\textwidth]{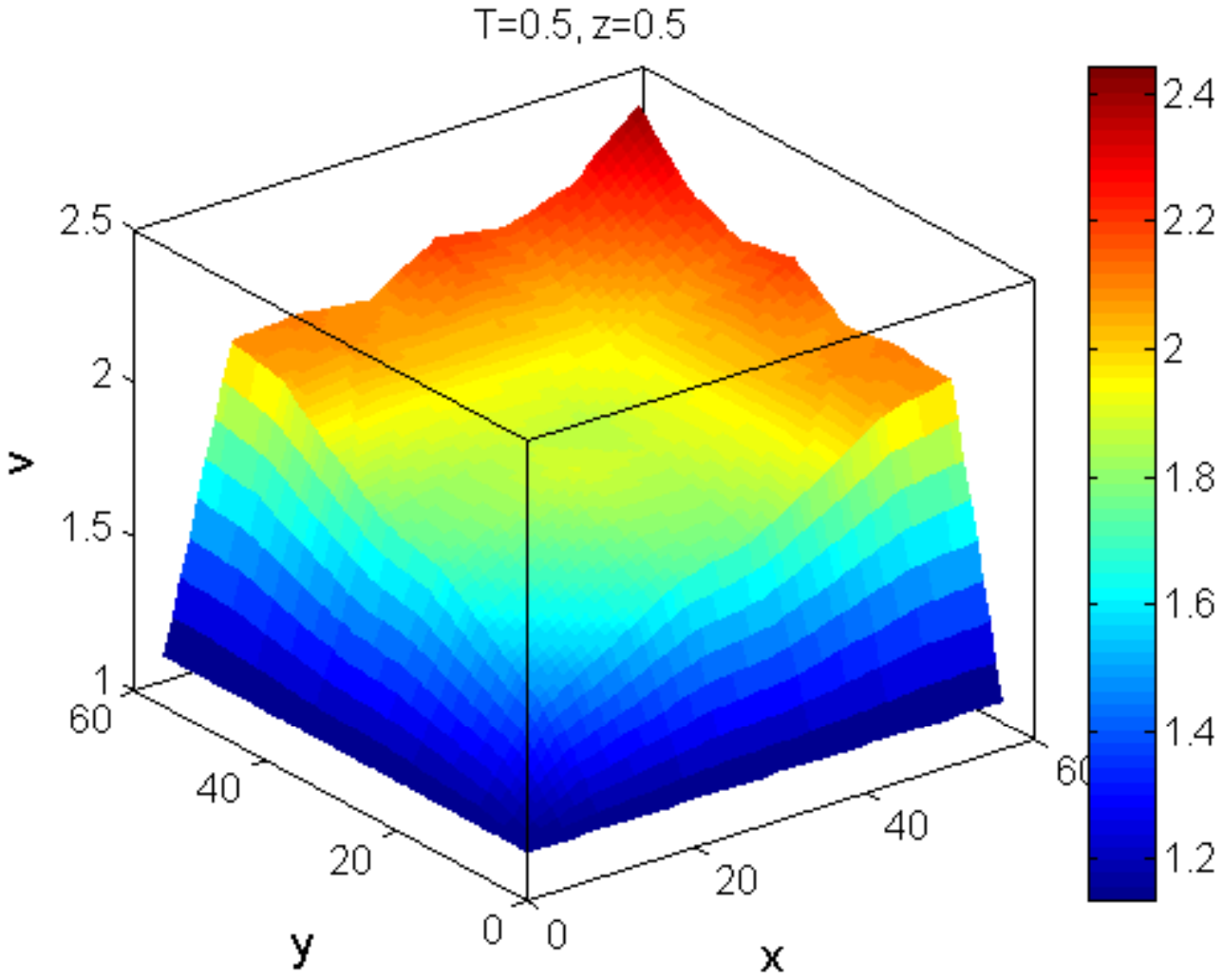}
		\vspace{-13em}
		\caption{}
		\label{ms3d_surf_xy_T100}
	\end{subfigure}
	\vspace{-10em}
	\caption{\textbf{HW solution along z=0.5 at time, (a) T=0.05, (b) T=0.05 (rotate view) (c)T=0.5.}}
	\label{ms3d_surf_xy}
\end{figure}

\begin{figure}[h]
	\vspace{-12em}
	\centering
	\begin{subfigure}[t]{0.5\textwidth}
		\centering
		\includegraphics[width=1.2\textwidth]{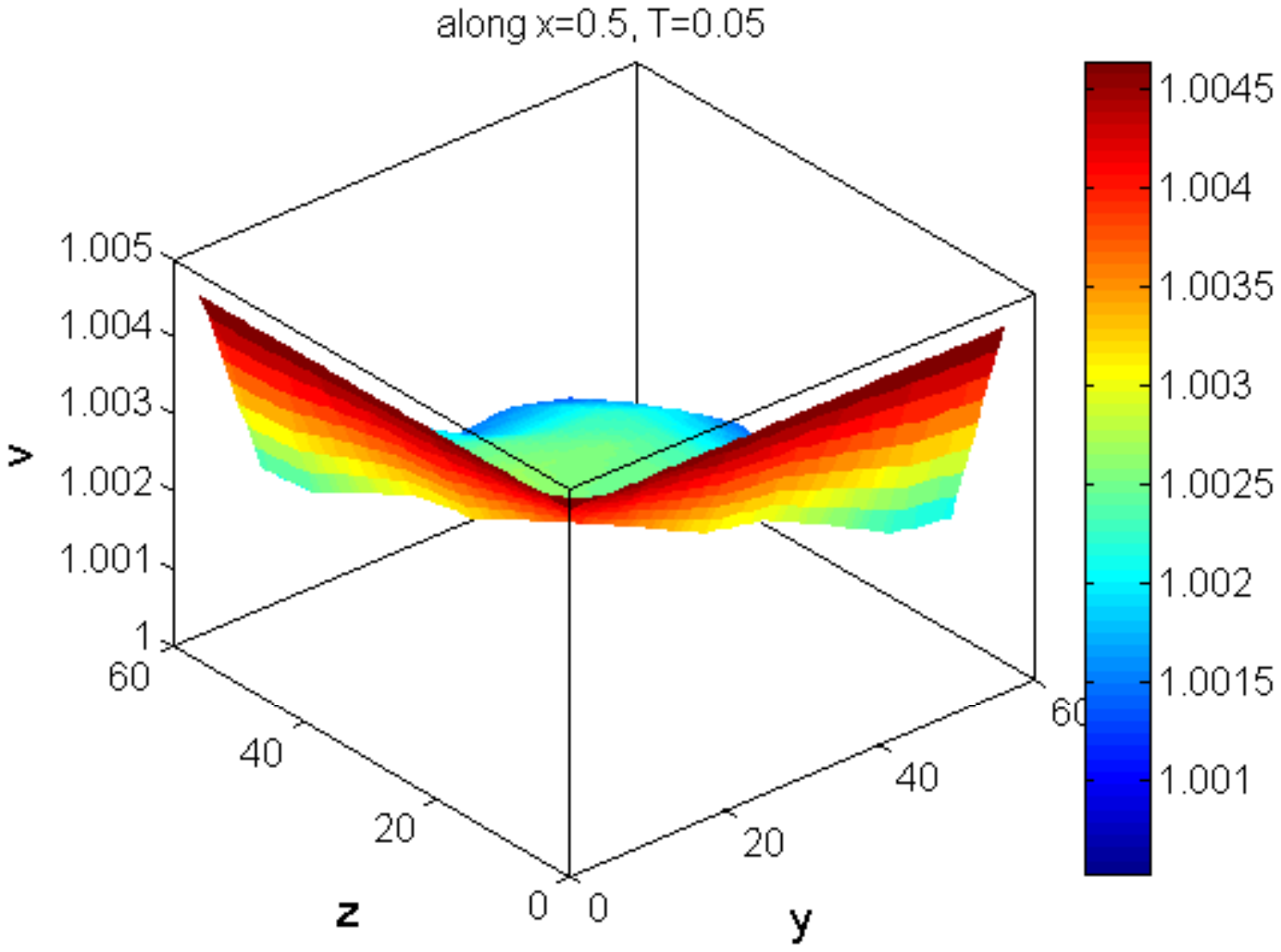}
		\vspace{-13em}
		\caption{}
		\label{ms3d_surf_yz_T50}
	\end{subfigure}\hfill
	\begin{subfigure}[t]{0.5\textwidth}
		\centering
		\includegraphics[width=1.2\textwidth]{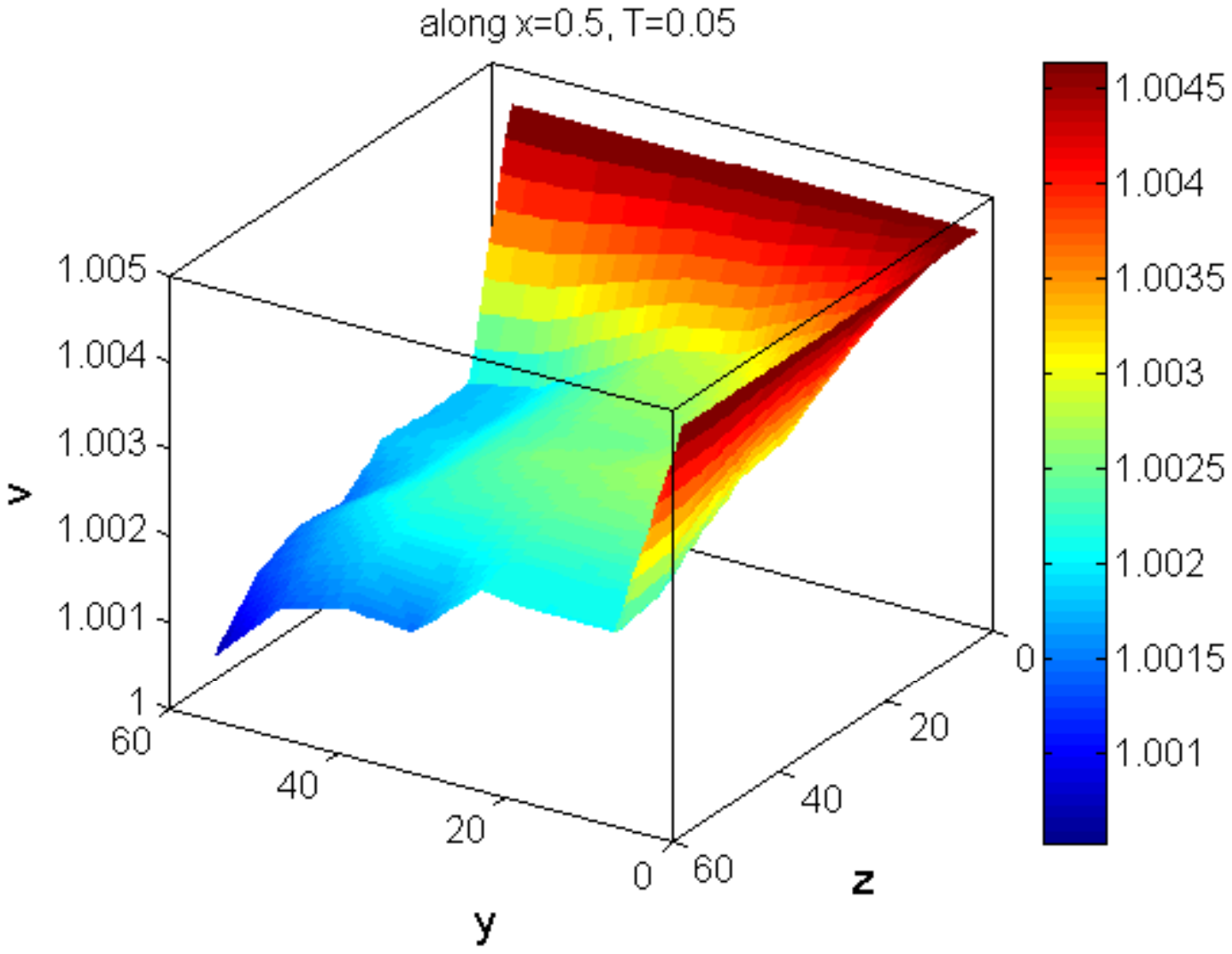}
		\vspace{-13em}
		\caption{}
		\label{ms3d_surf_yz_T50_rotate}
	\end{subfigure}\\
	\vspace{-16em}
	\begin{subfigure}[t]{0.5\textwidth}
		\centering
		\includegraphics[width=1.2\textwidth]{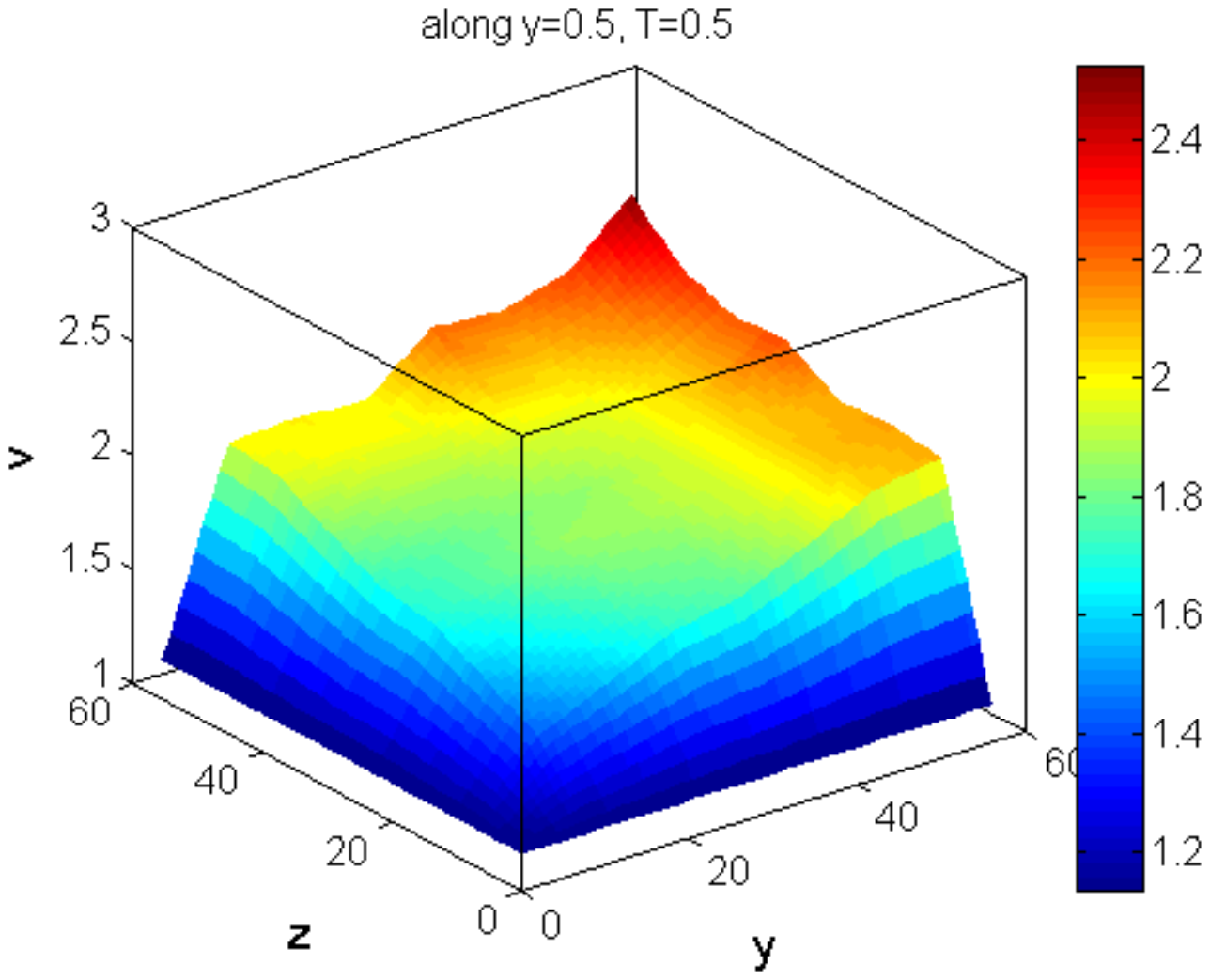}
		\vspace{-13em}
		\caption{}
		\label{ms3d_surf_yz_T100}
	\end{subfigure}
	\vspace{-10em}
	\caption{\textbf{HW solution along x=0.5 at time, (a) T=0.05, (b) T=0.05 (rotate view) (c)T=0.5}}
	\label{ms3d_surf_yz}
\end{figure}

\begin{figure}[h]
	\vspace{-12em}
	\centering
	\begin{subfigure}[t]{0.5\textwidth}
		\centering
		\includegraphics[width=1.2\textwidth]{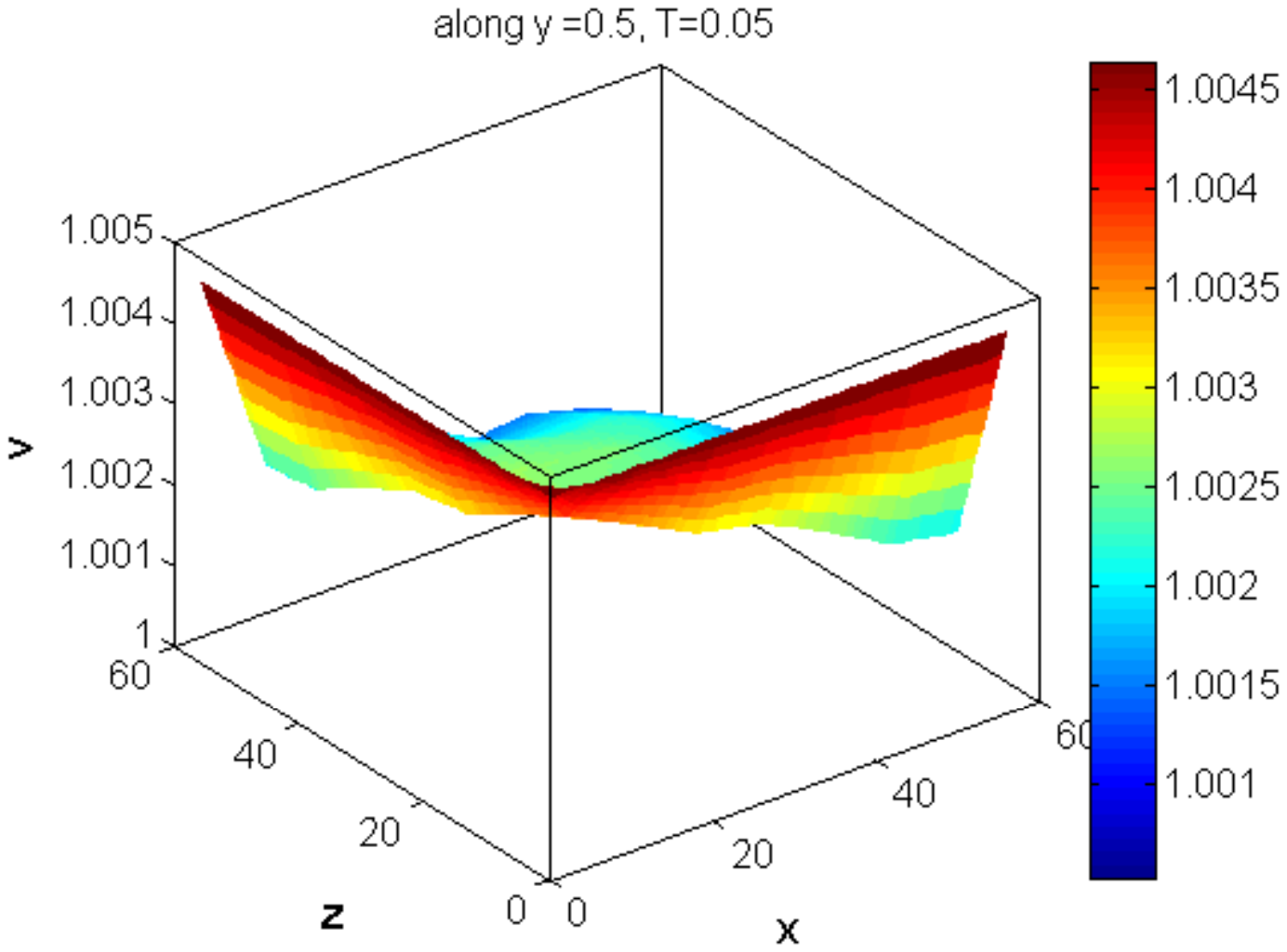}
		\vspace{-13em}
		\caption{}
		\label{ms3d_surf_zx_T50}
	\end{subfigure}\hfill
	\begin{subfigure}[t]{0.5\textwidth}
		\centering
		\includegraphics[width=1.2\textwidth]{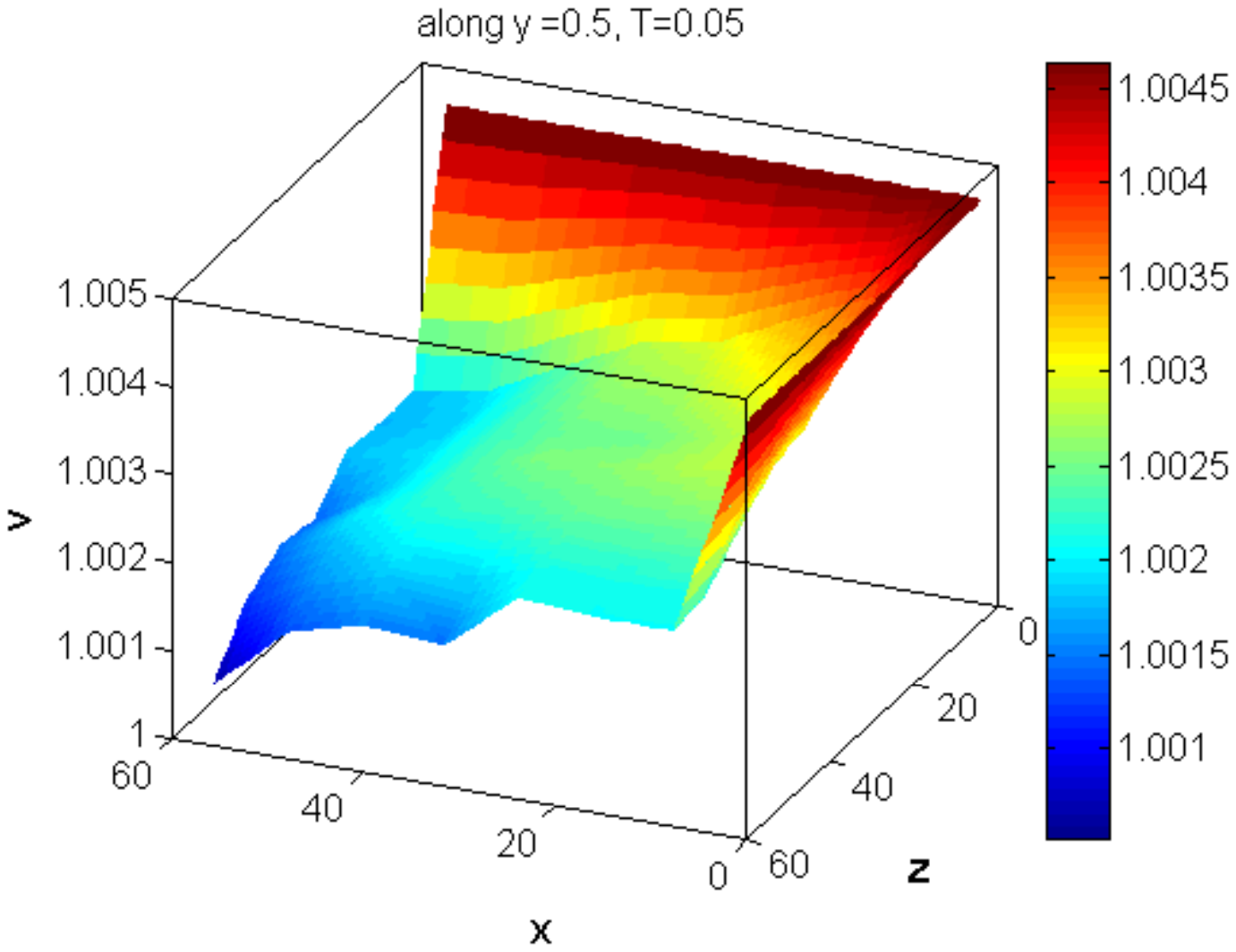}
		\vspace{-13em}
		\caption{}
		\label{ms3d_surf_zx_T50_rotate}
	\end{subfigure}\\
	\vspace{-16em}
	\begin{subfigure}[t]{0.5\textwidth}
		\centering
		\includegraphics[width=1.2\textwidth]{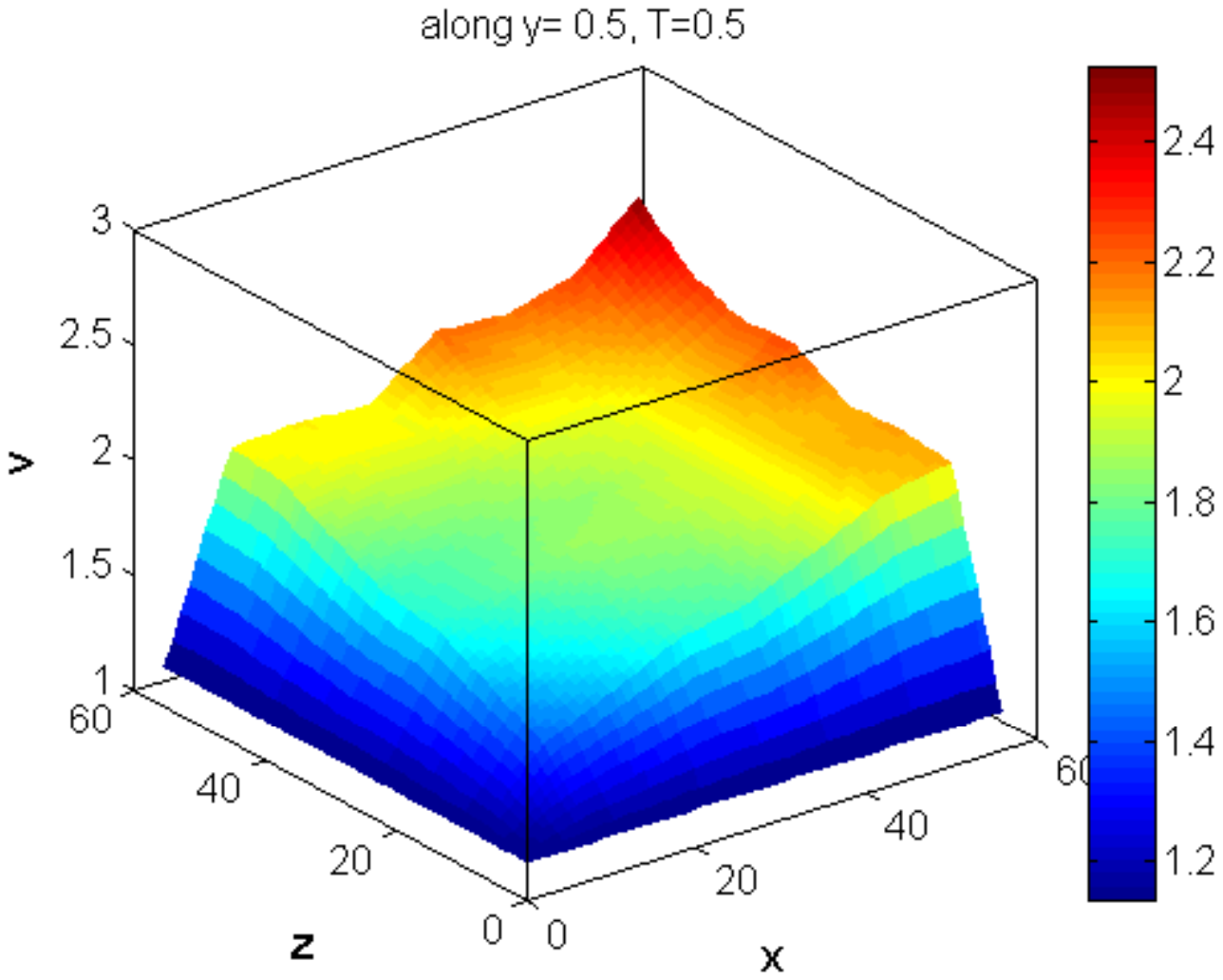}
		\vspace{-13em}
		\caption{}
		\label{ms3d_surf_zx_T100}
	\end{subfigure}
	\vspace{-10em}
	\caption{\textbf{HW solution along y=0.5 at time, (a) T=0.05, (b) T=0.05 (rotate view) (c)T=0.5.}}
	\label{ms3d_surf_zx}
\end{figure}

\begin{figure}[h]
	\vspace{-14em}
	\centering
	\begin{subfigure}[t]{0.5\textwidth}
		\centering
		\includegraphics[width=1.3\textwidth]{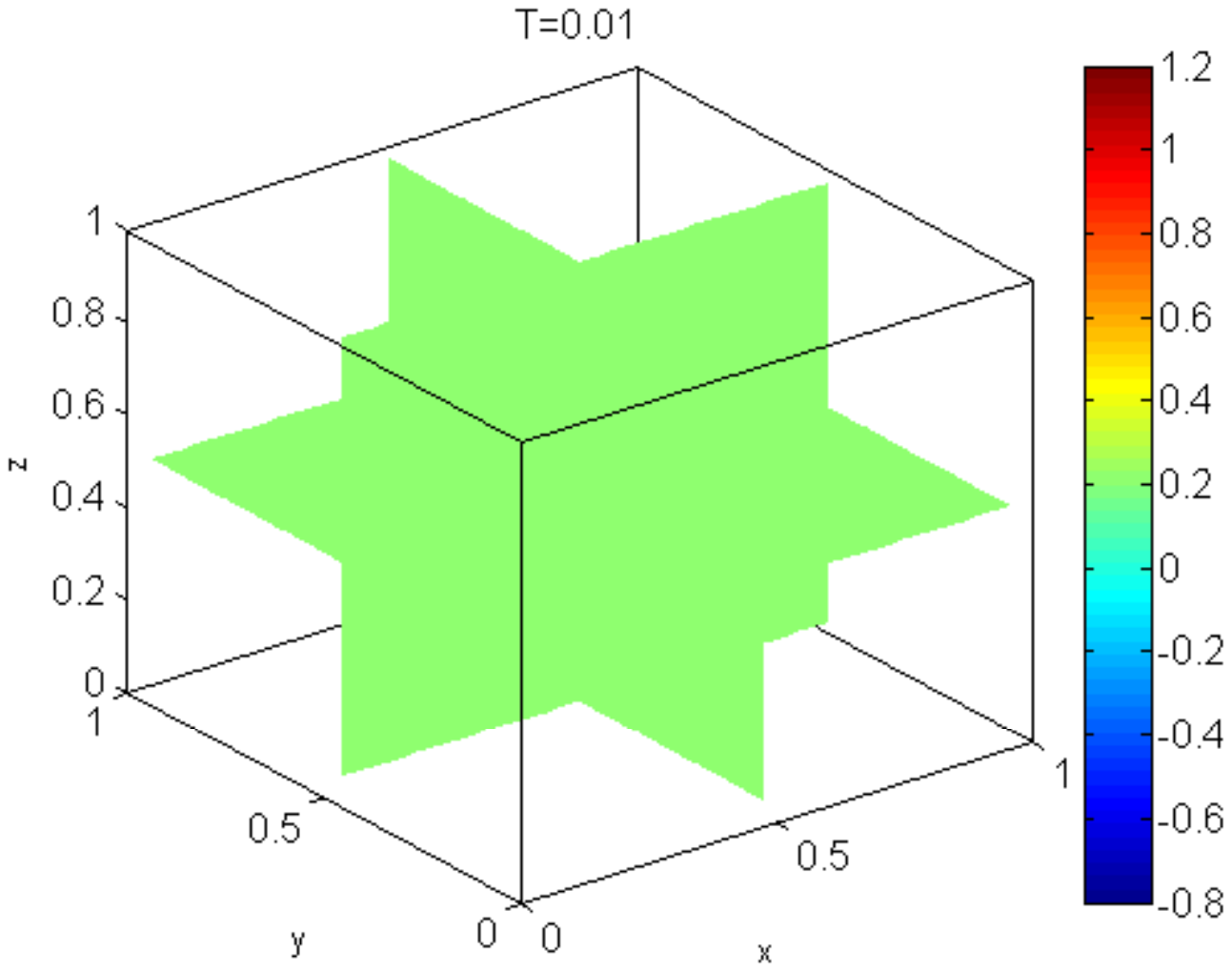}
		\vspace{-12em}
		\caption{}
		\label{T25}
	\end{subfigure}\hfill
	\begin{subfigure}[t]{0.5\textwidth}
		\centering
		\includegraphics[width=1.3\textwidth]{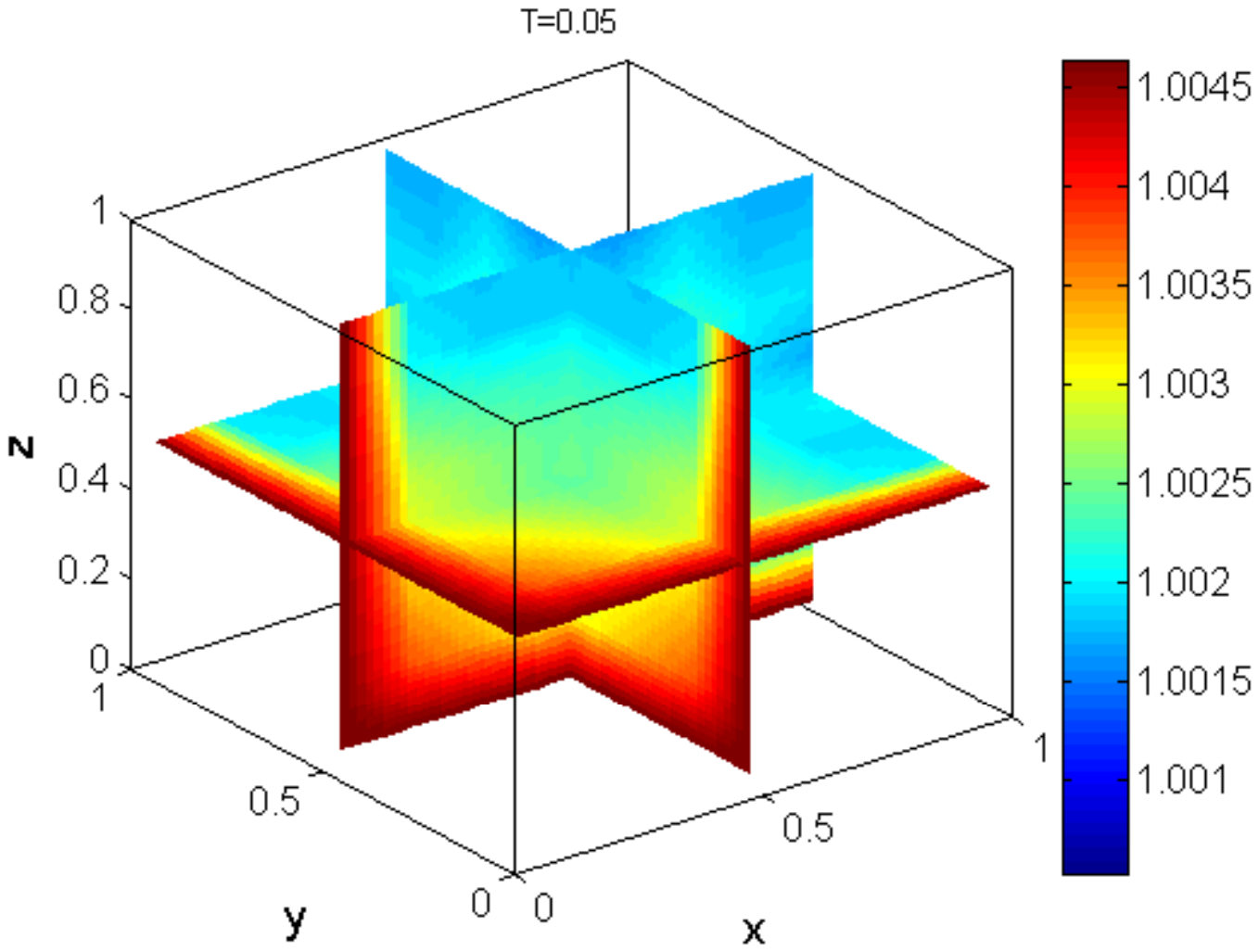}
		\vspace{-12em}
		\caption{}
		\label{T05}
	\end{subfigure}\\
	\vspace{-18em}
	\begin{subfigure}[t]{0.5\textwidth}
		\centering
		\includegraphics[width=1.3\textwidth]{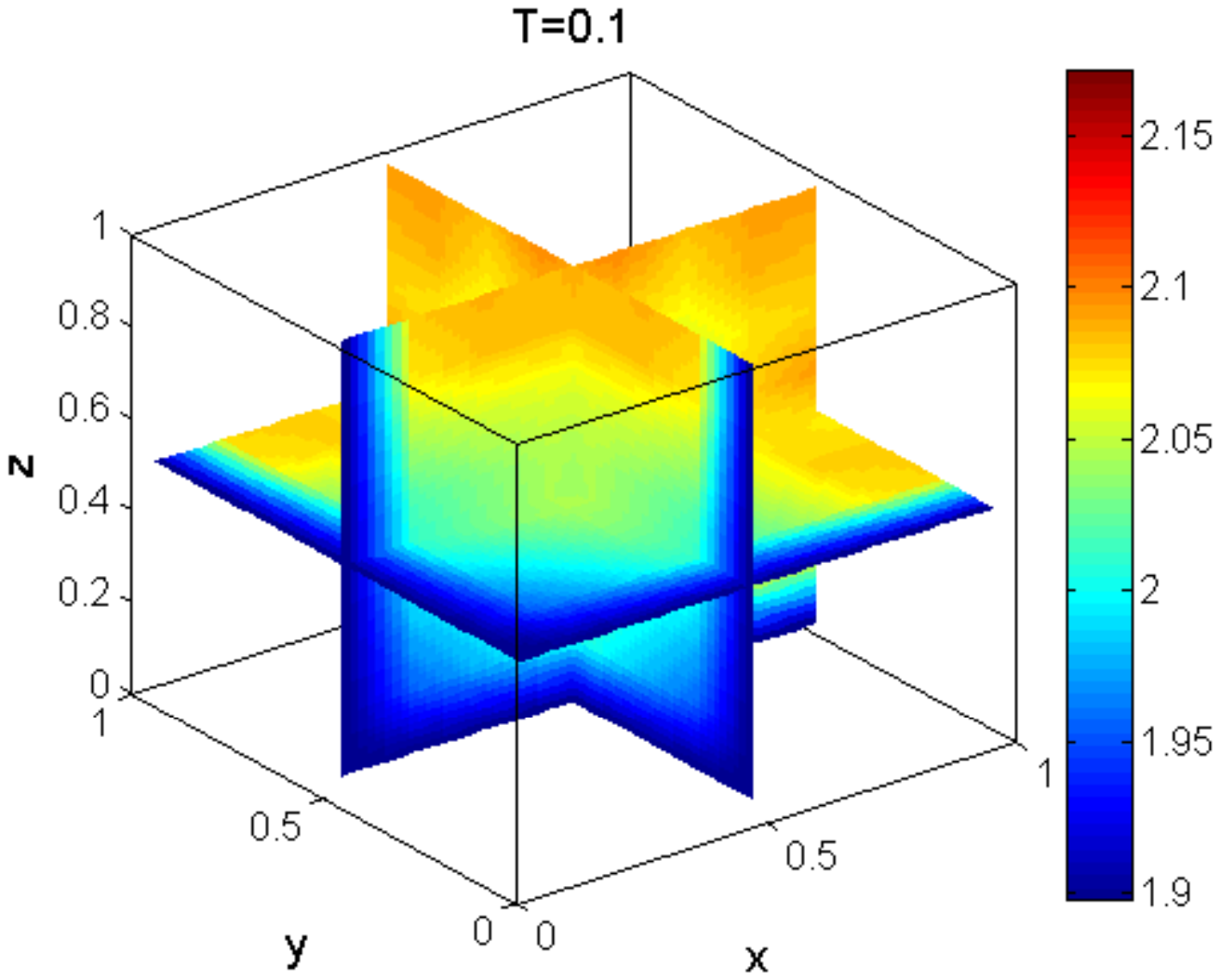}
		\vspace{-12em}
		\caption{}
		\label{T75}
	\end{subfigure}\hfill
	\begin{subfigure}[t]{0.5\textwidth}
		\centering
		\includegraphics[width=1.3\textwidth]{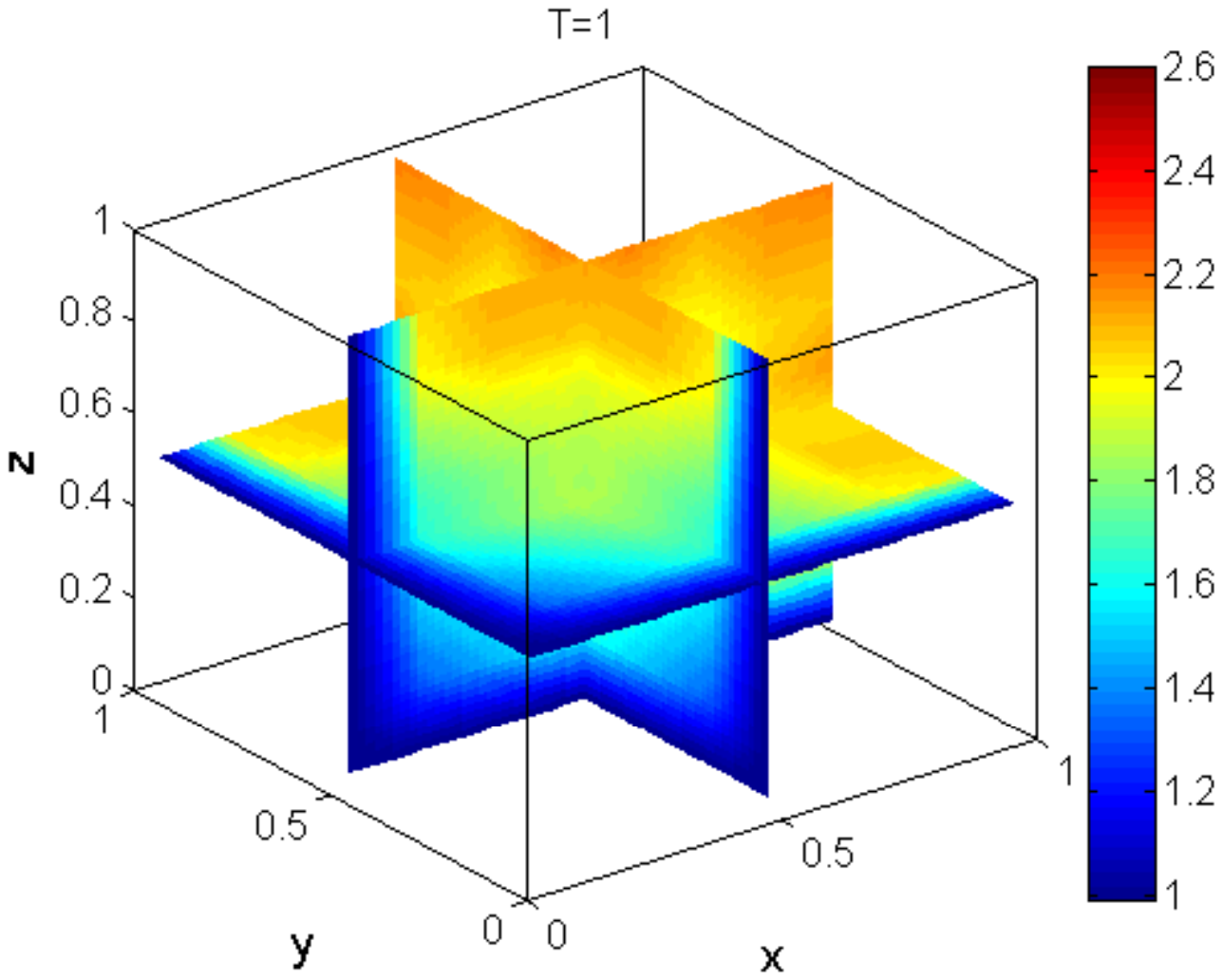}
		\vspace{-12em}
		\caption{}
		\label{T1}
	\end{subfigure}
	\vspace{-10em}
	\caption{\textbf{HW solution in the Three dimensional domain at J=4, $dt=10^{-2}$ at time, (a) T=0.01, (b) T=0.05,(c) T=0.1, (d) T=1.}}
	\label{3d_SlicePlot_ms}
\end{figure}

Comparison of the solution of this system using Haar wavelet method and FEM is presented in Fig. \ref{ms_hw_3D_comparison}. 
We have drawn the solution along the three cross section in the 3-dimensional domain in figures \eqref{ms3d_surf_xy}, \eqref{ms3d_surf_yz} and \eqref{ms3d_surf_zx} respectively at time T=0.05 and T=0.5. We have also shown the solution $v$ in 3-dimensional slice at different time levels in Fig \ref{3d_SlicePlot_ms}. 

%
%

\begin{figure}
	\vspace{-18em}
	\centering
	\includegraphics[width=0.8\textwidth]{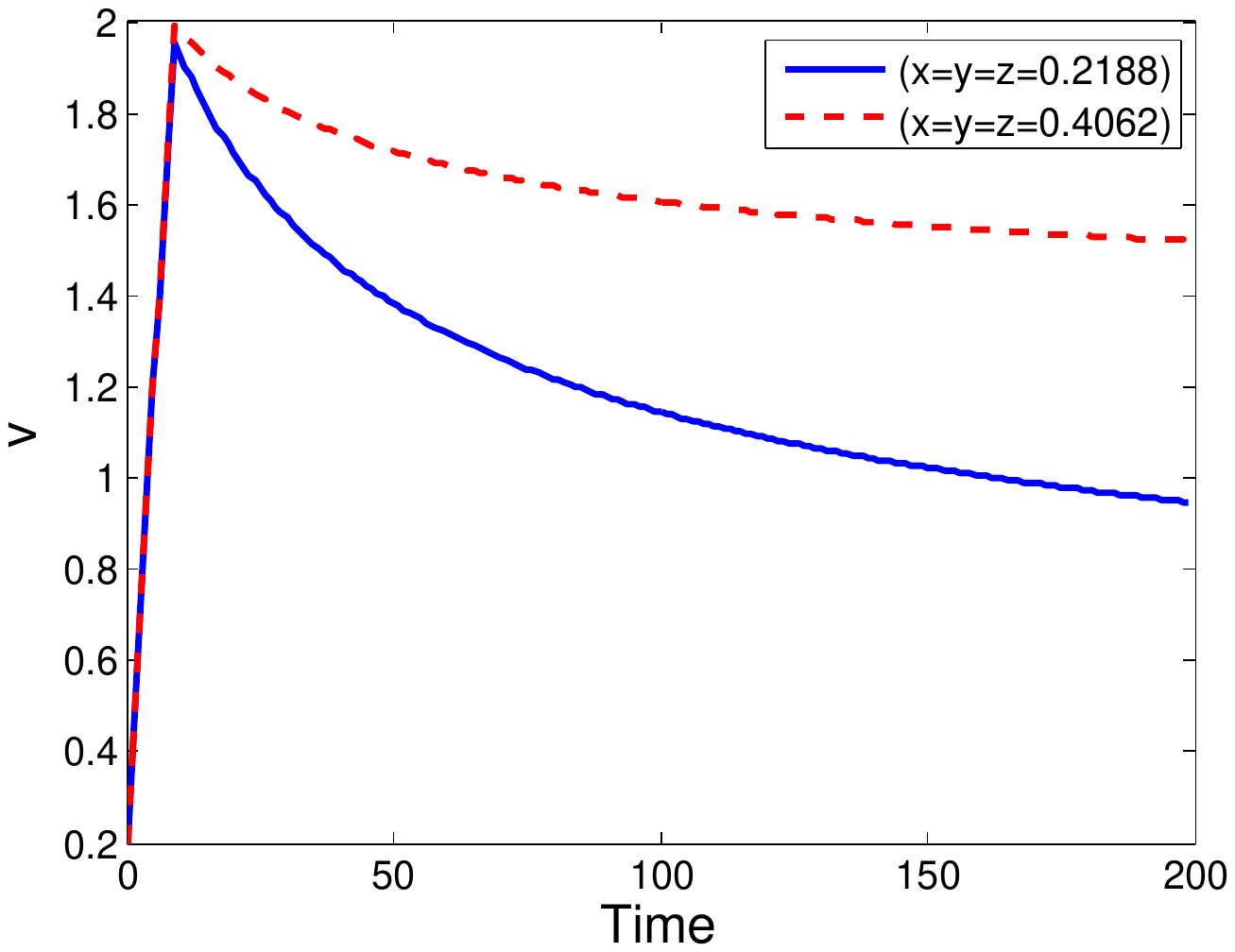}
	\centering
	\vspace{-16em}
	\caption{Solution $v$ with $\tau_{in}$ parameter value having jump discontinuity.}
	\label{AP_ms3d_jump}
\end{figure}

Next, we calculate the solution of the problem if the parameter $\tau_{in}$ have jump discontinuity at some places in the domain, 
the corresponding function $I_{ion}(v,w)$ will also be discontinuous, as shown in two dimensional case. The solution $v$ corresponding to this $\tau_{in}$ at the points (0.4062, 0.4062, 0.4062) and 
(0.2188, 0.2188, 0.2188), which are presented inside and outside the jump region respectively, of the domain is shown 
in Fig. \ref{AP_ms3d_jump}. The smoothness of the solution with jump clearly shows that the Haar wavelets handle the jump discontinuity in the value of parameter $\tau_{in}$.

We also solved this problem using the linear finite element in space and implicit-explicit (IE) Euler method in time. Comparison of the computation time of finite element solution and the Haar wavelet solution for 3D level model in example 5 is given in Table \ref{cputime_3d}. Clearly, from the table, the Haar wavelet method reduces the computational cost in comparison of FEM. Thus, it could be concluded that in higher dimension it is easy to use the Haar wavelet method to solve any problem even having the discontinuities.
\begin{table}[h]
	\begin{center}
		\begin{tabular}{ | c |c | c | c | c| } 
			\hline
			& J & dt & FEM CPU time (sec) & HW CPU time (sec)\\
			%
			\hline 
			
			Example 5 & 2 (256 grid) & $10^{-2}$ & 73.7888 & 40.9545\\
			
			& 3 (1024 grid) & $10^{-2}$ & 571.902 & 143.249 \\
			\hline 	
		\end{tabular}
		\caption{CPU time in 3D}
		\label{cputime_3d}
	\end{center}
\end{table}

\begin{table}[h]
	\begin{center}
		\begin{tabular}{ |c| c |c | c | c |  } 
			\hline
			J & cgs & bicg & bicgstab & gmres \\
			\hline 
			$J=5, dt=10^{-3}$(1D model) & 7.629 & 6.70 & 8.345 & 5.2011  \\
			\hline
			$J=4, dt=10^{-3}$(2D model)& 6.9015 & 6.8364 & 7.07158 & 6.4849  \\
			\hline
			$J=3, dt=10^{-2}$(3D model) & 85.14 & 81.1456 & 83.8982 & 72.875 \\
			\hline 
		\end{tabular}
		\caption{CPU time with different solvers for 1D, 2D, and, 3D problems}
		\label{cputime_pcg}
	\end{center}
\end{table}

\section{Conclusion}
A Haar Wavelet Method for a class of coupled non-linear PDE-ODEs system with jump discontinuities in model parameters or model 
coefficients has been proposed. The method is both simple and easy to implement in two and three dimensions. It is also 
computationally cost effective over linear finite element method. ILU- preconditioned GMRES Krylov solver accelerates the numerical solution 
convergence of the related linear system to required tolerances faster than other regular Krylov solvers. Convergence analysis has also been done to ensure the stability and accuracy. Numerical error reduces with the decrease in time step size or increase in resolution level. Model  problems having single or multiple jump discontinuities in some parameters have been successfully solved. Behavior of the solutions clearly indicated that Haar wavelets automatically handle such type of discontinuities. Problems with clinical relevance have also been successfully dealt with.


\begin{acknowledgements}
Authors would to acknowledge SERB. Govt of India for funding the project (SERB/F/10196/2017-2018) on WGM for PDEs. We would also like to thank the DST for support through Inspire Fellowship, ID no. is IF130906.
\end{acknowledgements}


\bibliographystyle{spbasic}      
\bibliography{reference}
%
%

\end{document}